\documentclass[graybox]{svmult}

\usepackage{type1cm}        
\usepackage{makeidx}         
\usepackage{graphicx}        
\usepackage{multicol}        
\usepackage[bottom]{footmisc}

\usepackage{newtxtext}       %
\usepackage[varvw]{newtxmath}       
\usepackage{amsmath}

\makeindex             

\let\Bbb\mathbb
\let\goth\mathfrak

\def\og{\leavevmode\raise.3ex\hbox{$\scriptscriptstyle\langle\!\langle\,$}}

\def\fgf{\/\leavevmode\raise.3ex\hbox{$\scriptscriptstyle\,\rangle\!\rangle$}}

\def\fg{\fgf\ }

\newcommand{\carre}{\qed}

\newcommand{\findem}{\ensuremath\blacksquare}

\let\stz\ss

\def\Demd#1|{\parindent=0pt\par{\sl D\'emonstration d#1}\pointir\parindent=20pt}

\let\petcap\sc

\def\finc{\vskip12pt}

\def\Defn{\parindent=0pt\par{\rm\bf D\'efinition.\ }\parindent=20pt}

\def\Dem{\parindent=0pt\par{\sl D\'emonstration}\pointir\parindent=20pt}

\def\Demo#1|{\parindent=0pt\par{\sl D\'emonstration #1}\pointir\parindent=20pt}

\def\Defns{\parindent=0pt\par{\rm\bf D\'efinitions.\ }\parindent=20pt}

\makeatletter
\newdimen\indentTh\indentTh=0pt
\def\p@int{{\rm .}}
\def\p@intir{\discretionary{\rm .}{}{\rm .\kern.35em---\kern.7em}}
\def\pointir{\afterassignment\pointir@\let\next=}
\def\pointir@{\ifx\next\par\p@int\else\p@intir\fi\next}
\long\def\Thc#1|#2\finc{\Th{}{#1}{\pointir}{#2}}
\long\def\Th#1#2#3#4{\parindent=\indentTh\par\vskip5pt
{#1}{\petcap #2}{\sl #3}\parindent=20pt{\sl #4\par}\vskip 5pt\parindent=20pt}

\newdimen\indentssec\indentssec=20pt
\newdimen\indentrem\indentssec=0pt

\def\Demdsp#1|{\parindent=0pt\par{\sl D\'emonstration d#1.}\parindent=20pt}
\def\oldstyle{}
\def\build#1_#2^#3{\mathrel{\mathop{\kern 0pt#1}\limits_{#2}^{#3}}}

\def\hdfl#1#2{\smash{\mathop{\hbox to 12mm{\rightarrowfill}}
\limits^{\scriptstyle#1}_{\scriptstyle#2}}}

\def\hdhfl#1#2{\smash{\mathop{\hbox to 12mm{\hookrightarrowfill}}
\limits^{\scriptstyle#1}_{\scriptstyle#2}}}

\def\hgfl#1#2{\smash{\mathop{\hbox to 12mm{\leftarrowfill}}
\limits^{\scriptstyle#1}_{\scriptstyle#2}}}

\def\hghfl#1#2{\smash{\mathop{\hbox to 12mm{\hookleftarrowfill}}
\limits^{\scriptstyle#1}_{\scriptstyle#2}}}

\newcount\secno \secno=0
\newcount\ssecno \ssecno=0
\newcount\sssecno \sssecno=0
\newcount\chapno \chapno=0
\newcount\notenumber \notenumber=1
\newcount\exino \exino=0
\newcount\expno \expno=0
\newdimen\indentsec\indentsec=20pt
\newdimen\indentssec\indentssec=20pt
\newdimen\indentsssec\indentsssec=20pt
\newdimen\indentrem\indentssec=0pt

\newdimen\indentTh\indentTh=0pt

\newdimen\indentth\indentssec=0pt
\def\sectiongen#1#2#3{\parindent=\indentsec\par\vskip .3cm
\vskip 0mm plus -20mm minus 1,5mm\penalty-50
{\bf #1}{\bf #2}{#3}\nobreak\parindent=20pt}%
\def\secc#1|{\sectiongen{}{#1}{\pointir}}
\def\nsecc#1|%
{\global\advance\secno by 1\global\ssecno=0\global\sssecno=0
\sectiongen{\the\secno\ }{#1}{\pointir}}
\def\secp#1|{\sectiongen{}{#1}{}\par}
\def\nsecp#1|%
{\global\advance\secno by 1\global\ssecno=0\global\sssecno=0
\sectiongen{\the\secno\ }{#1}{}\par}
\def\ssectiongen#1#2#3#4{\parindent=\indentssec\par\vskip .2cm
\vskip 0mm plus -20mm minus 1,5mm\penalty-50
{\bf #1}{\sl #2}{\sl #3}{#4}\nobreak\medskip\parindent=20pt}%
\def\ssecc#1|#2{\ssectiongen{}{#1}{\pointir}{#2}}
\def\nssecc#1|#2{\global\advance\ssecno by 1\global\sssecno=0
\ssectiongen{\the\secno.\the\ssecno\ }{#1}{\pointir}{#2}}
\def\ssecp#1|{\ssectiongen{}{#1}{}{}\par}
\def\nssecp#1|{\global\advance\ssecno by 1\global\sssecno=0
\ssectiongen{\the\secno.\the\ssecno\ }{#1}{}{}\par}
\def\sssectiongen#1#2#3#4{\parindent=\indentsssec\par\vskip .2cm
\vskip 0mm plus -20mm minus 1,5mm\penalty-50
{\bf #1}{\sl #2}{\sl #3}{#4}\nobreak\medskip\parindent=20pt}%
\def\sssecc#1|#2{\sssectiongen{}{#1}{\pointir}{#2}}
\def\nsssecc#1|#2{\global\advance\sssecno by 1
\ssectiongen{\the\secno.\the\ssecno.\the\sssecno\ }{#1}{\pointir}{#2}}
\def\sssecp#1|{\sssectiongen{}{#1}{}{}\par}
\def\nsssecp#1|{\global\advance\sssecno by 1
\sssectiongen{\the\secno.\the\ssecno.\the\sssecno\ }{#1}{}{}\par}
\long\def\Th#1#2#3#4{\parindent=\indentTh\par\vskip5pt
{#1}{\petcap #2}{\sl #3}\parindent=20pt{\sl #4\par}\vskip 5pt\parindent=20pt}
\long\def\pTh#1#2#3#4{\parindent=\indentth\par\vskip5pt
{#1}{\small \bf #2}{\small \sl #3}\parindent=20pt{\small \sl #4\par}\vskip 5pt\parindent=20pt}
\long\def\remarque#1#2#3#4{\parindent=\indentrem\par\vskip5pt
{#1}{\small \sl #2}{\small\sl #3}\parindent=20pt{\small#4\par}
\vskip 5pt\parindent=20pt}
\long\def\remarques#1#2#3#4{\parindent=\indentrem\par\vskip5pt
{#1}{\small \sl #2}{\small\sl #3}{\small#4\par}
\vskip 5pt\parindent=20pt}

\long\def\remarquesa#1#2#3#4#5{\parindent=\indentrem\par\vskip5pt
{#1}{\small \sl #2}{\small\sl #3}{\small#4}
\parindent=20pt{\small#5\par}
\vskip 5pt\parindent=20pt}

\long\def\Remarque#1#2#3#4{\parindent=\indentrem\par\vskip5pt
{#1}{ \sl #2}{\sl #3}\parindent=20pt{#4\par}
\vskip 5pt\parindent=20pt}
\long\def\remarquesn#1#2#3#4{\parindent=\indentrem\par\vskip5pt
{\small \sl #1}{\small #2}{\small\sl #3}{\small#4\par}
\vskip 5pt\parindent=20pt}

\long\def\Remarquen#1#2#3#4{\parindent=\indentrem\par\vskip5pt
{ \sl #1}{#2}{\sl #3}\parindent=20pt{#4\par}
\vskip 5pt\parindent=20pt}

\long\def\Thc#1|#2\finc{\Th{}{#1}{\pointir}{#2}}
\long\def\Thnc#1|#2|#3\finnc{\Th{#1}{#2}{\pointir}{#3}}
\long\def\Exic#1|#2\finc{{\global\advance\exino by 1}\Remarquen%
{Exercice }{\the\exino}{ #1\pointir}{#2}}
\long\def\Expc#1|#2\finc{{\global\advance\expno by 1}\Remarquen%
{Exemple }{\the\expno}{ #1\pointir}{#2}}
\long\def\exic#1|#2\finc{{\global\advance\exino by 1}\remarquesn%
{\bf Exercice }{\the\exino}{ #1\pointir}{#2}}
\long\def\expc#1|#2\finc{{\global\advance\expno by 1}\remarquesn%
{Exemple }{\the\expno}{ #1\pointir}{#2}}

\long\def\Ec#1\finc{\Th{}{}{}{#1}}
\long\def\thc#1|#2\finc{\pTh{}{#1}{\pointir}{#2}}%
\long\def\thsnc#1\finc{\pTh{}{}{}{#1}}

\long\def\Thp#1|#2\finp{\Th{}{#1}{\par}{#2}}
\long\def\thp#1|#2\finp{\pTh{}{#1}{\par}{#2}}
\long\def\rmc#1|#2\finc{\remarque{}{#1}{\pointir}{#2}}
\long\def\Rmc#1|#2\finc{\Remarque{}{#1}{\pointir}{#2}}

\long\def\rmp#1|#2\finp{\remarque{}{#1}{\par}{#2}}
\long\def\Rmp#1|#2\finp{\Remarque{}{#1}{\par}{#2}}

\long\def\parc#1\finc{\remarque{}{}{}{#1}}
\long\def\parcs#1\fincs{\remarques{}{}{}{#1}}

\long\def\parcsa#1\fins#2\fincsa{\remarquesa{}{}{}{#1}{#2}}

\def\Rm#1|{\parindent=0pt\par\vskip5pt{\sl #1}\pointir\parindent=20pt}
\global\ssecno=0\global\secno=0
 \def\page {\leaders\hbox to 2mm{\hfil.\hfil}\hfill
}
\def\npage {\vfill\eject \global\setcounter{footnote}{0}}

\begin{document}

 
\title*{Alg\`ebre lin\'eaire%
.}

\author{{\it par} Alexis Marin}
\maketitle

\vskip-45mm

\centerline{\small\sl R\'esum\'e}
{\small\sl
La premi\`ere partie de ce cours ramass\'e donne le langage des applications lin\'eaires entre modules \`a gauche et \`a droite, sous-modules et modules quotients,
noyau et image, produits et sommes directes, dualit\'e, transposition et application canonique vers le bidual.

La seconde pr\'esente le calcul matriciel et les d\'eterminants dans un anneau commutatif
comme application d'\og identit\'es remarquables\fg
dans l'anneaux de polyn\^omes \`a c\oe fficients entiers en les variables c\oe ficients et second membre de la
\og m\'ethode de Gau\stz\  g\'en\'erale\fgf,
ce qui requiert la factorialit\'e des anneaux de polyn\^omes \`a c\oe fficients entiers,
factorialit\'e  ici \'etablie en modifiant la d\'emonstration de Zermolo de la factorialit\'e de l'anneau des entiers.

La troisi\`eme, revenant au cas g\'en\'eral non n\'ecessairement commutatif,
pr\'esente les espaces vectoriels sur un corps et leur dimension. Exploitant la pr\'esentation,
aujourd'hui m\'econnue, que la premi\`ere \'edition de du trait\'e d'Alg\`ebre de Nicolas Bourbaki donne
des structures rationnelles sur un sous corps on obtient, sans th\'eorie des alg\`ebres centrales simples,
le baba de la th\'eorie des corps non commutatifs, le th\'eor\`eme de Wedderburn sur les corps finis
et celui d'Erd\"os-Kaplanski sur la dimension du dual d'un espace vectoriel de dimension infinie.
}
\vskip5mm
\centerline{\small\sl Abstract}
{\small\sl
After the language of module and theirs morphisms, this short course presents  matricial calculus and determinants in a commutative ring as appliction of ``remarquable identities'' in the ring of polynomials with integer coefficients with variable coefficients and second member of the ``general Gau\stz\ method'' requiring the factoriality of such polynomial rings,  here obtained by a modification of Zermolo's proof of the factoriality of the integers. As third part,
using the, today almost forgotten, presentation in the first edition of N. Bourbaki's Algebra
of rational structures for a subfield, one gets, without central simple algebras, the baba of skewfields,
and the theorems of Wederburn and Erd\"os-Kaplanski (commutativity of finite fields and dimension of the dual).
}

\vskip10mm
\parc
On utilisera les baba ensemblistes
(entiers naturels, ensemble, applications, injection, surjections
et les propri{\'e}t{\'e}s {\'e}l{\'e}mentaires des
cardinaux),  relire {\bf [Ld]\/} et
les \S 0 {\`a} \S 8 de
{\bf [Go]\/}.

Le {\bf 2.2\/} utilise qu'un anneau de polyn{\^o}mes
({\`a} nombre fini de variables) {\`a}
coefficients entiers relatifs est factoriel.
Ce fait tr{\`e}s {\'e}l{\'e}mentaire
(voir la preuve\footnote{\small
n'utilisant pas qu'un anneau de polyn{\^o}mes sur un corps est principal (donc factoriel).
}
de l'appendice {\bf 2.2.4\/}) adaptant celle de Zermolo 
[voir p. 11-12 de {\bf [Sa1]\/}] pour les entiers)
est  rarement expos{\'e} dans les manuels.
Pour l'approche classique  voir
les exercices  {\bf 21\/} du \S{\bf 31\/} et
{\bf 31\/} du \S{\bf 32\/} de  {\bf [Go]\/}, les Chap.
\uppercase\expandafter{\romannumeral1} et 
\uppercase\expandafter{\romannumeral2} de {\bf [Sa1]\/},
ou les \S18 puis \S30 et \S 31, ainsi que l'algorithme \S32 de {\bf [vW]\/},
ou encore le \S20 de {\bf [Wb]\/}.

Les {\bf 3.5.3\/}, {\bf 3.5.4\/} utilisent un peu
de  r{\'e}currence transfinie
(voir chap. \uppercase\expandafter{\romannumeral2}
de {\bf [Kr]\/} et/ou {\bf [Ha]\/}),
on consultera aussi la seconde partie (p. 17-30) du cours de magist{\`e}re {\bf [Gu]\/}.
\finc
\vfill

\npage
\nsecp Modules.|
\parc Soient
$f, g\!: M\!\rightarrow\!N$ 
deux homorphismes
 d'un groupe
ab{\'e}lien 
$M$ 
vers un groupe ab{\'e}lien 
$N$%
, leur somme 
$h=f+g$ 
est l'application qui {\`a} 
$m\in M$
 associe 
$h(m)=f(m)+g(m)\in N$%
. Comme l'addition
$+$
de
$N$
est commutative\footnote{\small %
$h(x\!+\!y)\!=\!f(x\!+\!y)+g(x\!+\!y)\!=\!f(x)\!+\!f(y)+g(x)+g(y)%
\!=\!f(x)+g(x)+f(y)+g(y)\!=\!h(x)+h(y)$
}%
,  c'est un homorphisme de groupe
ab{\'e}lien.
Muni de cette loi d'addition, l'ensemble
${\rm Hom\/}(M, N)$ 
des homomorphismes de 
$M$
 vers 
$N$
est un groupe ab{\'e}lien. 

Dans le cas o{\`u} 
$M=N$%
, source et but co{\"\i}ncidant,
deux homomorphismes 
$f, g : M\rightarrow M$%
, des {\it endomorphismes\/} de
 $M$%
,
se composent et
 $f\!\circ\!g : M\rightarrow M,\ 
x\mapsto f\!\circ\!g\,(x)=f\,(g\,(x))$
est un endomorphisme. 
Cette deuxi{\`e}me loi est associative, 
a l'identit{\'e}
${\rm Id\/}_M$
 de
 $M$
 comme {\'e}l{\'e}ment neutre et est distributive 
par rapport {\`a} l'addition, ainsi
 $({\rm End\/}\,(M), +, \circ)$ 
est un anneau. L'anneau oppos{\'e}
 ${\rm End\/}\,(M)^{\rm o\/}$ 
peut s'interpr{\'e}ter comme le m{\^e}me
ensemble muni des m{\^e}mes lois, mais en 
{\'e}crivant la variable {\`a} gauche
 de l'application : l'image par l'application
 $f$
 du point 
$m$
 est not{\'e}e
$(m)f$%
.
\finc

\Defns Soit 
$\Lambda$ 
un anneau, dit {\it anneau des scalaires\/}, un {\it
 $\Lambda$%
-module {\`a} gauche\/}

(respectivement {\it {\`a} droite\/})
est un groupe ab{\'e}lien 
$M$
 muni d'un morphisme d'anneau
$\gamma : \Lambda\rightarrow{\rm End\/}\,(M)$
(respectivement
 $\delta : 
\Lambda\rightarrow{\rm End\/}\,(M)^{\rm o\/}$%
).

{\small
Pour tout 
$\lambda\!\in\!\Lambda$
 et 
$m\!\in\!M$ 
on note alors
$\gamma(\lambda)(m)\!=\!\lambda\,m$
(resp. 
$\delta(\lambda)(m)$
 (ou 
$(m)\delta(\lambda))\!=\!\, m\lambda$%
. L'endomorphisme
 $\gamma(\lambda)$
(resp.
$\delta(\lambda)$%
) est l'{\it homoth{\'e}tie\/}
({\it {\`a} gauche\/}) (resp.
 ({\it {\`a} droite\/})) de rapport 
$\lambda$%
.
}

\Defn
Un
{\it
$\Lambda$%
-module bilat{\'e}ral\/}
est un groupe ab{\'e}lien
$M$
muni de structures de
$\Lambda$%
-modules {\`a} gauche et {\`a} droite
$\gamma, \delta : \Lambda\rightarrow{\rm End\/}\,(M)$
commutant~:\hfill\break
 pour tout
$\lambda, \mu\!\in\!\Lambda, m\!\in\! M$
on a
$(\lambda m)\mu\!=\!
\delta(\mu)\circ\gamma(\lambda)(m)\!=\!%
\gamma(\lambda)\circ\delta(\mu)(m)\!=\!%
\lambda(m \mu)$%
.
\rmc Exemples| 
a) Si
$M\!\!=\!\!\Lambda$
(ou plus g{\'e}n{\'e}ralement
$\Lambda^n, n\!\in\!{\bf N\/}$%
)
les multiplications {\`a} gauche et {\`a} droite
$((\lambda, \mu), m)\!\mapsto\!\lambda m\mu$
(resp.
$((\lambda, \mu), (m_i)_{1\leq i\leq n})\!\mapsto\!%
(\lambda m_i\mu)_{1\leq i\leq n})$%
)
d'un {\'e}l{\'e}ment
$m\!\in\!M\!=\!\Lambda$
(resp.
$\Lambda^n$%
)
par les scalaires
$\lambda, \mu\!\in\!\Lambda$
donnent {\`a}
$M$
sa {\it structure de
$\Lambda$%
-module bilat{\'e}ral usuelle\/}.

b) Si l'anneau
 $\Lambda$
 est commutatif les notions de module
$M$
{\`a} gauche et {\`a}
droite sont pratiquement identiques\footnote{\small
Si
$M$
est un
$\Lambda$%
-module {\`a} gauche alors
pour tout
$\lambda, \mu\!\in\!\Lambda$
et
$m\!\in\!M$
on a~:%
\hfill\break
$\lambda\cdot(\mu\cdot m)\!=\!(\lambda\mu)\cdot m\!=\!%
(\mu\lambda)\cdot m\!=\!\mu\cdot(\lambda\cdot m)$
donc
$(m, \lambda)\mapsto m\cdot \lambda
\build{=}_{}^{\rm Def}
\lambda\cdot m$
donne {\`a}
$M$
une structure de
$\Lambda$%
-module {\`a} droite commutant {\`a}
la structure initiale. De m{\^e}me un
 un 
$\Lambda$%
-module {\`a} droite 
$N$
est, par
$(\lambda, n)\mapsto \lambda\cdot n
\build{=}_{}^{\rm Def}
n\cdot\lambda$
un 
$\Lambda$%
-module {\`a} gauche et un
$\Lambda$%
-module bilat{\'e}ral.
}
et produisent une structure de
$\Lambda$%
-module bilat{\'e}rale dite {\it canonique\/}.

c) Tout groupe ab{\'e}lien 
$M$ 
a une unique structure de
${\bf Z\/}$%
-module bilat\'erale~:\hfill\break
si 
$k\geq0$%
,
$\gamma\,(k)(m)=
\delta\,(k)(m)=
\underbrace{m+\cdots+m}_{k\ {\hbox{\rm fois\/}}}$
 et si
 $k\leq0$%
, 
$-(\underbrace{m+\cdots+m}_{-k\ {\hbox{\rm fois\/}}})$%
.
\finc

\parc Si l'anneau
 $\Lambda$
n'est pas commutatif les notions de module
$M$
{\`a} gauche et {\`a}
droite sont distinctes et
quand un ensemble est muni de deux structures de 
$\Lambda$%
-module, l'une  {\`a} gauche l'autre {\`a} droite, 
comme par exemple l'anneau 
$\Lambda$%
, on pr{\'e}cisera
de quel c{\^o}t{\'e} on se place~:
les structure {\`a} gauche et {\`a}
droite de 
$\Lambda$
 sont not{\'e}es\footnote{\small
en suivant l'usage (issu du latin)
$s$
et
$d$
pour sinistra et dextra.
}
$\Lambda_s, \Lambda_d$%
.

Hormis pour la structure de 
$\Lambda$%
-module sur le dual servant {\`a} produire
 le morphisme canonique d'un 
$\Lambda$%
-module dans son bidual, et les applications
de la dimension aux corps gauches
cet opuscule ne s'occupe que de situations
 o{\`u} tous les modules sont du m{\^e}me c{\^o}t{\'e}.

 Conform{\'e}ment {\`a} un usage trop bien {\'e}tabli
pour {\^e}tre {\'e}branl{\'e},
il est {\'e}crit pour les
 $\Lambda$%
-modules {\`a} gauche 
 (mais, suivant l'usage,
le calcul matriciel est d{\'e}velopp{\'e} pour les
 $\Lambda$%
-modules {\`a} droite!).


\finc
\npage
\nssecp Applications lin{\'e}aires,
morphismes et isomorphismes.|
\Defns Une application {\it 
$\Lambda$%
-lin{\'e}aire\/} 
d'un 
$\Lambda$%
-module 
$M$
 vers un
 $\Lambda$%
-module 
$N$
est un homomorphisme 
$f : M\rightarrow N$ 
des groupes ab{\'e}liens sous-jacents
qui commute avec les homoth{\'e}ties : pour tout 
$m\in M$
 et 
$\lambda\in \Lambda$
on a 
$f\,(\lambda\,m)=\lambda\,f(m)$%
.

{\small
Ainsi une application
$f : M\rightarrow N$
entre deux 
$\Lambda$%
-modules est%
\footnote{\small
$f(x\!-\!y)\!=\!f(1 x+(\!-1)y)\!=\!%
f(x)-f(y)$
et
$f(\lambda x)\!=\!f(\lambda x\!+\!0 0)\!=\!%
\lambda f(x)\!+\!0f(0)\!=\!\lambda f(x)$%
.
} 
$\Lambda$%
-lin{\'e}aire si et seulement
si pour tout
$x, y\in M$
et
$\lambda, \mu\in\Lambda$
on a~:
$f(\lambda x+\mu y)=
\lambda f(x)+\mu f(u)$%
.
}

Un synonyme d'application lin{\'e}aire est
{\it morphisme de
 $\Lambda$%
-modules\/}
ou {\it 
$\Lambda$%
-mor\-phisme\/}. Si l'anneau 
$\Lambda$
est fix{\'e} par le contexte on parlera simplement
de {\it modules\/}, {\it applications lin{\'e}aires\/}
 ou {\it morphismes\/}.

L'identit{\'e} 
${\mathop{\rm Id}\nolimits}_M$
 d'un module 
$M$
est un morphisme  et le compos{\'e}
$g\circ f : M\rightarrow P$
 de deux morphismes
 $f : M\rightarrow N$ 
et
$g : N\rightarrow P$
est un morphisme.

L'ensemble des morphismes de
$\Lambda$%
-modules de 
$M$
 vers 
$N$
est not{\'e} 
${\mathop{\rm Hom\/}\nolimits}_\Lambda(M, N)$%
. C'est un sous-groupe du groupe
 ${\mathop{\rm Hom\/}\nolimits}\,(M, N)=
{\mathop{\rm Hom\/}\nolimits}_{\Bbb Z\/}\,(M, N)$
des homomorphismes de groupe ab{\'e}lien de 
$M$
 vers
 $N$%
.

Un morphisme 
$f :M\rightarrow N$ 
est un
{\it isomorphisme\/} si il y a un morphisme 
$g : N\rightarrow M$
tel que
$g\circ f={\mathop{\rm Id}\nolimits}_M$
et 
$f\circ g={\mathop{\rm Id}\nolimits}_N$%
.

{\small
En ce cas l'application
$f$
est bijective d'application r{\'e}ciproque
$f^{-1}=g$%
. En fait~:
}
\Thc Lemme|
Soit
$f : M\rightarrow N$
un morphisme bijectif. \hfill\break
Alors son application r{\'e}ciproque 
$g=f^{-1} : N\rightarrow M$
est un morphisme. 

Ainsi
un morphisme est un isomorphisme si et seuleument si
il est bijectif.
\finc
{\small
\Dem Soit
$u, v\in N, \lambda,\mu\in\Lambda, w\!=\!\lambda g(u)+\mu g(v)\in M$
alors comme~:\hfill\break
$u\!=\!f(g(u)), v\!=\!f(g(v)), g(f(w))=w$
\hfill\break
on a~:
$g(\lambda u+\mu v)=
g(\lambda f(g(u))+\mu f(g(v)))=
g(f(\lambda g(u)+\mu g(v)))=
\lambda g(u)+\mu g(v)$%
.\hfill\findem
}

Dans le cas o{\`u} le module
$N$
est muni d'une structure bilat{\'e}rale\hfill\break
{\small
[par exemple l'usuelle
 (resp. la canonique) sur
$N=\Lambda^n$ 
(resp., si l'anneau
$\Lambda$
est commutatif)]
}
pour tout
$f\in {\mathop{\rm Hom\/}\nolimits}_\Lambda\,(M, N)$
et tout 
$\lambda, \mu\in\Lambda$
et
$y\in N$%
, on a~:
$$f\,(\lambda\,x)\mu=(\lambda\,f\,(x))\,\mu=
\lambda\,(f\,(x)\mu)$$
ainsi l'application
$f\,\mu :M\rightarrow N,\
 x\mapsto f\,(x)\,\mu$
est
 $\Lambda$%
-lin{\'e}aire, c'est
l'homoth{\'e}tique de 
$f$
de rapport
 $\mu$ pour la
{\it structure de 
$\Lambda$%
-module {\`a} droite\/} de
 ${\mathop{\rm Hom\/}\nolimits}_\Lambda\,(M, N)$%
.

\vskip10mm
{\small
\centerline{\small\petcap
 Commentaires bibliographiques\/}
\vskip3mm

Plus de d{\'e}tails sur  modules 
morphismes sont dans\footnote%
{\small
Un manuel \og{\'e}l{\'e}mentaire mais correct\fgf,
sans la \og standard mistake\fg
 de l'alg{\`e}bre lin{\'e}aire~:
se \og m{\'e}langer les  c{\^o}t{\'e} des scalaires\fg
en passant de l'alg{\`e}bre lin{\'e}aire abstraite
au calcul matriciel.
}
les \S10 {\`a} \S13 de {\bf [Go]\/}, le {\bf 1.4\/} de {\bf [Sa2]\/}, 
ou  les plus exaustifs Kap. 12 de {\bf [vW]\/}, 
\S1 et \S2 de {\bf [Bo]\/} Chap
{\uppercase\expandafter{\romannumeral2}}, 
{\bf [Ja]\/} vol. \uppercase\expandafter{\romannumeral1}
et chap. \uppercase\expandafter{\romannumeral3} de {\bf [Lg]\/}.

}
\npage
\nssecp Sous-modules et modules quotients.|
\Defn Un {\it sous-module\/} d'un
 $\Lambda$%
-module 
$M$
est un sous-groupe
ab{\'e}lien
$U\subset M$
de
$M$
stable par homoth{\'e}tie [%
pour tout $\lambda\in \Lambda$
on a~:
$\gamma\,(\lambda)(U)=:\lambda\,U\subset U$%
].

C'est un 
$\Lambda$%
-module et l'inclusion
$i_U : U\hookrightarrow M$
 est un morphisme.
\rmc Exemples|
$M, \{0\}\subset M$
sont des sous-modules dits 
{\it triviaux, plein et nul\/}.
On note
$0\!=\!\{0\}$%
.
\finc
\Thc Proposition 1| Soit 
$U_i, i\!\in\!I$
une famille (index\'ee par un ensemble
$I$%
)
de sous-modules  d'un module 
$M$%
. Alors :

(i) L'intersection
 $\cap_{i\in I}U_i$
de cette famille
est un sous-module de
 $M$%
,\hfill\break
 c'est le plus grand sous-module de
 $M$ 
inclus dans tous les 
$U_i$%
.

(ii) L'intersection des sous-modules contenant 
chacun tous les
 $U_i$
est le plus petit sous-module
 les contenant tous.\hfill\break
Il est form{\'e} des sommes finies d'{\'e}l{\'e}ments
$u_i\in U_i$%
, c'est la {\it somme\/} des 
$U_i$%
.\hfill\break
 On le note
$\sum_{i\in I}U_i$
 ou, dans le cas o{\`u} 
$I=\{1,\ldots,n\}$
est fini, 
$U_1+U_2+\cdots+U_n$%
.
\finc
\Thc Compl{\'e}ment et d{\'e}finitions|
Soit
$\underline{m}=(m_i)_{i\in I}, m_i\in M$
(resp.
$E\subset M$%
) une famille (resp. partie) d'un
$\Lambda$%
-module
$M$%
. Il y a un plus petit sous-module de
$M$
contenant tous les
$m_i, i\!\in\!I$
(resp. la partie
$E$%
) c'est le {\it sous-module engendr{\'e} par\/}
la famille
$\underline{m}$
(resp. la partie
$E$%
) not{\'e}e
$<\underline{m}>=<\underline{m}>_\Lambda$
(resp.
$<E>=<E>_\Lambda$%
). 

Si
$M=<\underline{m}>_\Lambda$%
, le module
$M$
est dit {\it engendr{\'e} par\/} la famille
$\underline{m}$%
.

Un module
$M$
 est dit {\it de type fini\/}
({\it engendr{\'e} sur
$\Lambda$
par\/}
$m_1,\ldots,m_n$%
)
si\hfill\break
$M=<m_1,\ldots, m_n>_\Lambda$%
, il est engendr{\'e} par une famille
$\underline{m}=m_1,\ldots, m_n$
finie.
\finc

\Defns Deux {\'e}l{\'e}ments
 $x, y\in M$
sont {\it congrus\/} modulo un sous-module 
$U$\hfill\break 
si 
$x-y\in U$%
, on note alors
$x\equiv y\ {\rm mod\/}\,U$%
. 

{\small
La congruence modulo un sous-module est 
une relation d'{\'e}quivalence et l'ensemble
quotient est muni d'une unique structure de 
$\Lambda$%
-module, le
{\it 
$\Lambda$%
-module quotient\/}\
 $M/U$%
, pour laquelle
l'application quotient est un morphisme, le
{\it morphisme quotient\/}\ 
$\pi : M\rightarrow M/U$%
. De plus~:
}
\Thc Proposition 2| La pr{\'e}composition 
$g\mapsto g\circ \pi$ 
par 
$\pi$ 
est un isomorphisme du groupe ab{\'e}lien 
${\rm Hom\/}_\Lambda\,(M/U, N)$ 
sur le sous-groupe
$\{f\in {\rm Hom\/}_\Lambda\,(M, N)| f_{|U}=0\}$ 
des morphismes de 
$M$ 
vers 
$N$
s'annulant sur 
$U$%
.

Si l'anneau 
$\Lambda$ 
est commutatif ou plus g\'en\'eralement si le module
$N$
est muni d'une structure bilat\'erale, c'est un isomorphisme de 
$\Lambda$%
-module.
\finc
\vskip3mm
{\small
\centerline{\small\petcap
 Commentaires bibliographiques\/}

La notion de sous-objet et 
structure quotient en alg{\`e}bre,
pr{\'e}par{\'e}e en exercice dans
{\bf [Go]\/}, et le point de vue de la
Proposition 2 (insister plus sur morphismes 
qu'objets) est plus syst{\'e}ma\-tiquement
utilis{\'e}e et d{\'e}velopp{\'e}e dans 
{\bf [Bo]\/} et {\bf [Lg]\/}, et le 
Vocabulaire Math{\'e}matique de
{\bf [Co]\/}.

}
\npage
\nssecp Noyau, image, factorisation canonique
et suites exactes.|
\Defns Soit 
$f\in{\mathop{\rm Hom\/}\nolimits}_\Lambda(M, N)$ 
un morphisme de
 $M$
 vers 
$N$
son {\it noyau\/} et son {\it image\/} sont
${\mathop{\rm ker\/}\nolimits}(f)=
\{m\in M|\ f(m)=0\}$ 
et 
${\mathop{\rm Im\/}\nolimits}(f)=\{f(m)|\ m\in M\}$%
.\hfill\break
Ce sont des sous-modules de
 $M$ 
et 
$N$
 respectivement.
Le {\it conoyau\/} et la {\it coimage\/} de 
$f$ 
sont respectivement~:
${\mathop{\rm coker\/}\nolimits}(f)=
N/{\mathop{\rm Im\/}\nolimits}(f)$ 
et 
${\mathop{\rm coim\/}\nolimits}(f)= 
M/{\mathop{\rm ker\/}\nolimits}(f)$%
. Ainsi~:

{\small
{\sl Un morphisme est surjectif si et seulement si
son conoyau
${\rm coker\/}\,(f)\!=\!0$
est nul.}\hfill\break
Plus utile que ce pl{\'e}onasme est la fondamentale~:
}

\Thc Proposition 1| Un morphisme 
$f : M\rightarrow N$
 est injectif 
ssi
${\rm ker\/}\,(f)=0$%
.
\finc
\vskip-2mm
\Thc Proposition 2| Tout morphisme 
$f : M\rightarrow N$
a une unique factorisation
\vskip-4mm
$$f : 
M\build{\longrightarrow}_{}^{\pi}{\rm coim\/}\,(f)
\build{\longrightarrow}_{}^{\overline{f}}{\rm Im\/}\,(f)
\build{\longrightarrow}_{}^{i_{{\rm Im\/}\,(f)}}N$$
\vskip-2mm
\noindent{\small
o{\`u} 
$\pi$
 est 
le morphisme quotient,$\overline{f}$
un isomorphisme et 
$i_{\rm Im\/}\,(f)$
l'inclusion de l'image dans le but.
}
\finc
{\small
\Defn Une {\it suite exacte\/} ({\it courte\/})
de 
$\Lambda$%
-modules, not{\'e}e
$0\!\rightarrow\!U\!\build{\rightarrow}_{}^{i}\!%
V\!\build{\rightarrow}_{}^{j}\!W\!%
\rightarrow\!0$%
, est la donn{\'e}e de morpismes
$i : U\rightarrow V$
injectif et
$j : V\rightarrow W$
surjectif  tels que
${\mathop{\rm ker\/}\nolimits}(j)=
{\mathop{\rm Im\/}\nolimits}(i)$%
.
\rmc Exemples|
a) Si
$U\subset M$
est un sous-module d'un module
$M$%
, on a la suite exacte~:
$$0\rightarrow U
\rightarrow M%
\rightarrow M/U
\rightarrow0$$

b) Si
$f : M\rightarrow N$
est un morphisme de modules, on a la suite exacte~:
$$0\rightarrow{\mathop{\rm ker\/}\nolimits}(f)%
\rightarrow M%
\rightarrow{\mathop{\rm Im\/}\nolimits}(f)%
\rightarrow0$$
\finc
\vskip-4mm
{\small
\Defn
Un {\it morphisme\/} d'une suite exacte
$0\rightarrow U\build{\rightarrow}_{}^{i}%
V\build{\rightarrow}_{}^{j}W%
\rightarrow0$
vers une suite exacte\break
$0\!\rightarrow\!U'\!\build{\rightarrow}_{}^{i'}\!%
V'\!\build{\rightarrow}_{}^{j'}\!W'\!\!%
\rightarrow\!0$
est la donn{\'e}e, not{\'e}e
$(f, g, h)$%
, de trois morphismes
$f\!: U\!\rightarrow\!U', g\!: V\!\rightarrow\!V'$
et
$h\!: W\!\rightarrow\!W'$
tels que
$g\circ i\!=\!i'\circ f$
et
$h\circ j\!=\!j'\circ g$%
, conditions que l'on exprime en disant que~:
$$\begin{matrix}
0\rightarrow & U & \build{\rightarrow}_{}^{i} & V & \build{\rightarrow}_{}^{j} & W & \rightarrow0  \\
\null  & f\downarrow\ \ & \null & g\downarrow\ \ & \null & h\downarrow\ \ & \null \\
0\rightarrow & U' & \build{\rightarrow}_{}^{i'} & V' & \build{\rightarrow}_{}^{j'} & W' & \rightarrow0
\end{matrix}$$
est un {\it diagramme commutatif.\/} On a alors
le {\it (petit\footnote{\small
Pour le vrai Lemme des cinq,
faire l'exercice {\bf 25\/} du {\bf \S7\/} de {\bf [Go]\/}.
}) Lemme des cinq\/}~:
}
\Thc Lemme| Si les deux morphismes
lat{\'e}raux 
$f$
et
$h$
sont des isomorphismes,\hfill\break
alors le morphisme central 
$g$
est un isomorphisme.
\finc
{\small
\Dem D'apr{\`e}s le lemme de {\bf 1.1\/} il
suffit de montrer que 
$g$
est bijectif. 

Soit
$v'\!\in\!V'$.
Comme
$h, j$
sont surjectifs, il y a
$w\!\in\!W, v_1\!\in\!V$
tel que
$h(w)\!=\!j'(v'), j(v_1)\!=\!w$%
. Ainsi
$j'(v'\!-\!g(v_1))\!=\!j'(v')\!-\!j'(g(v_1))\!=\!%
j'(v')\!-\!h(j(v_1))\!=\!0$%
, donc
$v'\!-\!g(v_1)\!\in\!{\mathop{\rm ker\/}\nolimits}(j')\!=\!%
{\mathop{\rm Im\/}\nolimits}(i')$
et il y a
$u'\!\in\!U'$
avec
$i'(u')\!=\!v'-g(v_1)$%
. Comme
$f$ est surjectif, il y a
$u\!\in\!U$
avec
$f(u)\!=\!u'$%
, d'o{\`u}
$g(i(u)+v_1)\!=\!g(i(u))+g(v_1)\!=\!i'(f(u))+g(v_1)\!=\!
v'-g(v_1)+g(v_1)\!=\!v'$
et
$g$
est surjectif.\hfill\carre

Soit
$v\!\in\!{\mathop{\rm ker\/}\nolimits}(g)$%
. On a
$h(j(v))\!=\!j'(g(v))\!=\!j'(0)\!=\!0$%
, donc
$j(v)\!\in\!{\mathop{\rm ker\/}\nolimits}(h)$%
, qui est nul d'apr{\`e}s la Proposition 1,
puisque
$h$ 
est injectif.
Ainsi
$v\!\in\!{\mathop{\rm ker\/}\nolimits}(j)\!=\!%
{\mathop{\rm Im\/}\nolimits}(i)$
et il y a
$u\!\in\!U$
tel que
$i(u)\!=\!v$%
. Comme
$i'\circ f(u)=g\circ i(u)=g(i(u))=g(v)=0$
et
$i', f$
sont injectifs on a~:
$u=0$%
, d'o{\`u}
$v=0$%
. Ainsi
${\mathop{\rm ker\/}\nolimits}(g)=0$
et la Proposition 1 assure que
$g$ 
est aussi injective, donc bijective.\hfill\carre\findem
}
}
\npage
\nssecp Produit, somme directe et 
module libre sur un ensemble.|
\Defn Soit 
$N_i, i\!\in\!I$ 
une famille de 
$\Lambda$%
-modules index{\'e}e par un ensemble
$I$%
.\hfill\break
Le
{\it 
$\Lambda$%
-module produit\/} 
$\prod_{i\in I}\!N_i$ 
est l'ensemble des suites
$(n_i\!)_{i\in I}$ 
indic{\'e}es par 
$I$ 
dont les 
${\imath}^{\hbox{\small i\`emes}}$
{\'e}l{\'e}ments sont des {\'e}l{\'e}ments de 
$N_i$%
. Cet ensemble de suites est
muni de l'addition et de l'action de 
$\Lambda$ 
composante par composante :
$$(n_i)_{i\in I}+(n'_i)_{i\in I}=%
(n_i+n'_i)_{i\in I}\ \hbox{\rm et\/}\
 \lambda\,(n_i)_{i\in I}=(\lambda\,n_i)_{i\in I}$$
Pour tout 
$k\in I$ 
le morphisme d{\'e}fini par
${\rm pr\/}_k((n_i)_{i\in I})=n_k$
est la
{\it 
$k^{\hbox{\small i\`eme}}$
projection\/}~:
$${\rm pr\/}_k :\prod_{i\in I}N_i\rightarrow N_k$$

Si 
$f_i : X\rightarrow N_i, i\!\in\!I$ 
est une famille index{\'e}e par $I$
de morphismes d'un
$\Lambda$%
-mo\-dule 
$X$ 
vers les 
$N_i$ 
alors il y a un unique morphisme
$f:X\rightarrow \prod_{i\in I}N_i$ 
tel que pour tout 
$i\in I$
on ait 
${\rm pr\/}_i\!\circ\!f=f_i$%
, il est donn{\'e} par 
$f(x)=(f_i(x))_{i\in I}$%
~:
\Thc Proposition 1|L'application
$f\mapsto ({\rm pr\/}_i\!\circ\!f)_{i\in I}$ 
est un isomorphime de groupe ab{\'e}lien de 
${\mathop{\rm Hom\/}\nolimits}_\Lambda(X, \prod_{i\in I}N_i)$ 
sur 
$\prod_{i\in I}%
{\mathop{\rm Hom\/}\nolimits}_\Lambda(X, N_i)$%
.

Si les modules
$N_i$
sont bilat{\'e}raux\footnote{\small
par exemple si l'anneau 
$\Lambda$ 
est commutatif et les
$N_i$
sont munis de leur structure canonique.
}
 c'est un isomorphisme de 
$\Lambda$%
-modules.
\finc
\Defn Soit 
$M_j, j\!\in\!J$ 
une famille de 
$\Lambda$%
-modules d'ensemble d'indices
$J$%
.\hfill\break
 Les suites 
$(m_j)_{j\in J}\!\in\!\prod_{j\in J}M_j$
dont {\it presque toutes\/} les composantes sont
nulles\footnote{\small 
{\it c. a d.\/} la suite 
$(m_j)_{j\in J}$
est de 
{\it support\/}%
, l'ensemble
$\{j\!\in\!J, m_j\!\ne\!0\}$%
,
fini.
}\hfill\break
forment un sous-%
$\Lambda$%
-module de 
$\prod_{j\in J}M_j$ 
appel{\'e}
{\it somme directe des\/}\ 
$M_j$%
, et not{\'e} 
$$\oplus_{j\in J}M_j$$
L'inclusion
${\rm i\/}_J : 
\oplus_{j\in J}M_j\rightarrow \prod_{j\in I}M_j$
 de la somme directe dans le produit est surjective
si et seulement si l'ensemble 
$\{j\in J|\ M_j\ne0\}$ 
est fini.

Pour tout 
$k\in J$ 
on a la 
{\it 
$k^{\hbox{i\`eme}}$
 inclusion\/}
$${\mathop{\rm i\/}\nolimits}_k:
M_k\rightarrow \oplus_{j\in J}M_j$$
le morphisme d{\'e}fini par 
${\mathop{\rm i\/}\nolimits}_k(m)=(m_j)_{j\in J}$ 
o{\`u} 
$m_k=m$ 
et 
$m_j=0$
pour 
$j\ne k$%
.

Si 
$g_j : M_j\rightarrow Y, j\in J$ 
est une famille, index{\'e}e par 
$J$%
, de morphismes de
 $M_j$
 vers un 
$\Lambda$%
-module 
$Y$ 
alors il y a un unique morphisme~:
$$g : \oplus_{j\in J}M_j\rightarrow Y$$
tel que pour tout
 $j\in J$ 
on ait~: 
$$g\circ {\mathop{\rm i\/}\nolimits}_j=g_j$$
il est donn{\'e} par
$g((m_j)_{j\in J})=\sum_{j\in J}g_j(m_j)$%
, cette somme {\'e}tant bien d{\'e}finie car, 
sauf pour un nombre fini d'entre eux, les 
$m_j$%
, donc les 
$g_j(m_j)$%
, sont nuls.

\Thc Proposition 2| L'application
$g\mapsto 
(g\circ{\mathop{\rm i\/}\nolimits_j})_{j\in J}$
est un isomorphisme du groupe ab{\'e}lien 
${\mathop{\rm Hom\/}\nolimits}_\Lambda%
(\oplus_{j\in J}M_j, Y)$
sur le groupe ab{\'e}lien 
$\prod_{j\in J}%
{\mathop{\rm Hom\/}\nolimits}_\Lambda(M_j, Y)$%
.

{\small
Si le module
$Y$
est muni d'une structure bilat{\'e}rale\footnote{\small
{\it p. e.\/} l'usuelle (resp. canonique)
si
$Y=\Lambda^n$
(resp.
$\Lambda$ 
est commutatif).
} 
c'est un isomorphisme de 
$\Lambda$%
-modules.
}
\finc
\Thc Lemme et d{\'e}finitions| Soit 
$P\subset M$ 
un sous-module de 
$M$ 
tel qu'il y ait un sous-module 
$Q\subset M$ 
avec 
$P\cap Q= 0$
et 
$P+Q=M$%
. Alors l'application 
$$P\oplus Q\rightarrow M$$ 
d{\'e}finie par les inclusions de 
$P$ 
et 
$Q$ 
dans 
$M$ 
est un isomorphisme. 

On dit alors que
 $P$ 
est un 
{\it facteur direct\/}
 de 
$M$%
, que 
$Q$ 
est un 
{\it suppl{\'e}mentaire\/}
de 
$P$ 
dans 
$M$ 
et que 
$M$ 
admet la 
{\it d{\'e}composition en somme directe\/}
$M=P\oplus Q$%
.

{\small
De m{\^e}me si une famille
$M_j\!\subset\!M,\, j\!\in\!J\/$ 
de sous-modules de 
$M\/$ 
v{\'e}rifie~:
$$\sum_{j\in J}M_j\!=\!M\
\hbox{\rm et pour chaque \/}\
j_0\!\in\!J,\ 
M_{j_0}\cap\bigl(\!\!\sum_{j\in J\setminus\{j_0\}}M_j\bigr)\!=\!0$$
alors le mor\-phis\-me 
$\oplus_{i\in J}M_j\rightarrow M$
d{\'e}fini par les inclusions des 
$M_j$ 
dans 
$M$ 
est un isomorphisme.
 
On dit alors que 
 $M$ 
admet la 
{\it d{\'e}composition en somme directe\/}
$M=\oplus_{i\in J}M_j$%
.
}
\finc

Ainsi si 
$M\!=\!\oplus_{j\in J}\!M_j$ 
est  d{\'e}composition en somme directe d'un 
$\Lambda$%
-mo\-dule 
$M$ 
et 
$N\!=\!\prod_{i\in I}N_i$ 
est d{\'e}composition
en produit d'un 
$\Lambda$%
-module 
$N$ 
l'application
$$f\mapsto
({\mathop{\rm pr\/}\nolimits}_i\circ f\circ%
{\mathop{\rm i\/}\nolimits}_j)_{(i,j)\in I\times J}$$
est un isomorphisme de groupe ab{\'e}lien de
${\mathop{\rm Hom\/}\nolimits}_\Lambda(M, N)$
sur 
$\displaystyle\!\!\!\prod_{(i,j)\in I\times J}\!\!\!\!%
{\mathop{\rm Hom\/}\nolimits}_\Lambda\,(M_j, N_i)$%
.
\Defn la 
{\it matrice d'endomorphismes\/} 
de 
$f$ 
associ{\'e}e {\`a} ces d{\'e}compositions est~:
$$(f_{i,j}\,(f) =
{\mathop{\rm pr\/}\nolimits}_i\circ f\circ%
{\mathop{\rm i\/}\nolimits}_j)_{(i,j)\in I\times J}
\in\prod_{(i,j)\in I\times J}%
{\mathop{\rm Hom\/}\nolimits}_\Lambda\,(M_j, N_i)$$

Si un troisi{\`e}me 
$\Lambda$%
-mo\-dule 
$L\!=\!\oplus_{k\in K}L_k$ 
est d{\'e}compos{\'e }en somme directe
et 
$g : L\rightarrow M$
une application lin{\'e}aire telle que la  matrice de 
${\mathop{\rm i\/}\nolimits}_J\!\circ\!g$
est 
$(g_{j,k})_{(j,k)\in J\times K}$%
. Alors la matrice 
$(h_{i,k})_{(i,k)\in I\times K}$
de l'application compos{\'e}e
$h = f\!\circ\!g : L\rightarrow M\rightarrow N$ 
est donn{\'e}e par la
formule du calcul matriciel des endomorphismes~: 
$$h_{i,k}=\sum_{j\in J}f_{i,j}\!\circ\!g_{j,k}$$
qui  a un sens m{\^e}me si 
$J$ 
est infini car le but de 
$g$ 
{\'e}tant la somme directe 
$M=\oplus_{j\in J} M_j$ 
pour chaque 
$k\in K$
et 
$z\in L$ 
il n'y a qu'un nombre fini de 
$j$ 
tels que 
$g_{j,k}\,(z)\ne0$
soit non nul et la somme 
$\sum_{j\in J}g_{j,k}\,(z)$
est en fait finie.

{\small
Soit 
$I$ 
un ensemble. 
On choisit pour chaque 
$i\in I$ 
un exemplaire 
$\Lambda_i$ 
de l'anneau 
$\Lambda$%
, on note son unit{\'e} 
$1_i$
et on le consid{\`e}re comme 
$\Lambda$%
-module {\`a} gauche, ainsi
$\Lambda_i=\{\lambda\,1_i|\ \lambda\in \Lambda\}$%
. 
}
\vfill\eject
\Defns Le 
{\it 
$\Lambda$%
-module libre sur $I$\/} 
est la somme directe
$$L_I=\oplus_{i\in I}\Lambda_i$$
{\small
(s'identifiant au module 
$\Lambda^{(I)}$
des lignes 
dans 
$\Lambda$
indic{\'e}es par 
$I$ 
d'{\'e}l{\'e}ments presques tous nuls),
}
munie de l'application
${\rm e\/}_I : I\rightarrow L_I$ 
d{\'e}finie par~:
$${\rm e\/}_I\,(i) = e_i=(e_{i_j})_{j\in I}\
\hbox{\rm o{\`u}\/}\
e_{i_j}=0\
\hbox{\rm si\/}\ j\ne i\ \hbox{\rm et\/}\ %
e_{i_i}=1_i$$

L'{\'e}l{\'e}ment
$e_i$
est dit
{\it {\'e}l{\'e}ment canonique d'indice
$i$
de base de\/}
$\Lambda_I$%
.

Cette application 
${\rm e\/}_I\!:\!I\rightarrow L_I$
a la propri{\'e}t{\'e} universelle suivante :
\Thc Proposition 3|Pour toute application
$f : I\!\rightarrow\!M$ 
d'un ensemble 
$I$
vers un 
$\Lambda$%
-module 
$M$
il y a un unique morphisme
$$L(f) : L_I\rightarrow M\quad
\hbox{\rm tel que\/}\quad 
L(f)\circ{\rm e\/}_I = f$$
c'est le morphisme de la somme directe 
$\oplus_{i\in I}\Lambda_i$ 
vers 
$M$ 
d{\'e}termin{\'e} par les morphismes
$f_i : \Lambda_i\rightarrow M,\ 
\lambda\,1_i\mapsto \lambda\,f(i)$%
.
\finc
De m{\^e}me que l'on a identifi{\'e} 
$L_I$ 
aux lignes d'{\'e}l{\'e}ments de 
$\Lambda$ 
indic{\'e}es par 
$I$%
, on utilisera l'abus de notation 
$f$ 
pour
$L(f)$%
.
\Defn Un 
{\it 
$\Lambda$%
-module libre de base 
$(b_i)_{i\in I}$
}
est un 
$\Lambda$%
-module 
$L$
muni d'un isomorphisme 
$\varphi : L_I\rightarrow L$ 
tel que pour tout
$i\in I$ 
on ait 
$\varphi\,(e_i) = b_i$%
.

 L'isomorphisme r{\'e}ciproque
$\psi : L\rightarrow L_I$ 
associe {\`a} chaque {\'e}l{\'e}ment 
$x\in L$
la ligne 
$(x_i)_{i\in I}$
de ses 
{\it coordonn{\'e}es dans la base\/} 
$(b_i)_{i\in I}$%
.

On a alors l'{\'e}criture unique 
$x=\sum_{i\in I}x_i\,b_i$%
. De plus l'inclusion 
$\{b_i|\ i\in I\}\hookrightarrow L$ 
poss{\`e}de la propri{\'e}t{\'e}
universelle de la {\petcap Proposition 3\/}.
\Thc Corollaire| Soit 
$f : M\!\rightarrow\!N$ 
une application 
$\Lambda$%
-lin{\'e}aire surjective et 
$g : L\!\rightarrow\!N$ 
un morphisme de source un 
$\Lambda$%
-module libre 
$L$ 
de base 
$(b_i)_{i\in I}$%
. 

Alors, si les 
$c_i\in M$ 
sont tels que 
$f\,(c_i) = g(b_i)$%
, l'application 
$\Lambda$%
-lin{\'e}aire 
$\tilde{g} : L\rightarrow M$
d{\'e}termin{\'e}e par 
$\tilde{g}\,(b_i) = c_i$  
v{\'e}rifie la relation 
$g= f\circ\tilde{g}$%
.
\finc

\vskip10mm
{\small
\centerline{\small\petcap
 Commentaires bibliographiques\/}
\vskip3mm

Les notions de somme et produit en alg{\`e}bre,
plus g{\'e}n{\'e}ralement de probl{\`e}me universel,
ainsi que le proc{\'e}d{\'e} de construction qui
{\`a} un ensemble 
$I$
associe le module libre sur
$I$
sont\footnote{\small
apr{\`e}s la notion sous-objet et de quotient de
{\bf 1.2\/}.
}
de nouvelles illustrations du
{\it point de vue des cat{\'e}gories et foncteurs\/}
qui est implicite dans  {\bf [vW]\/}, {\bf [Bo]\/} et est
plus amplement d{\'e}velopp{\'e} dans {\bf [Lg]\/}
et dans le Vocabulaire math{\'e}matique 
de {\bf [Co]\/}.
}
\npage
\nssecp Dualit{\'e}, transposition et
application canonique 
vers le bidual.|
{\small
\Defns Une 
{\it forme lin{\'e}aire\/} 
sur un 
$\Lambda$%
-module {\`a} gauche
$M$ 
est une application 
$\Lambda$%
-lin{\'e}aire {\`a} gauche
de
$M$ 
vers l'anneau 
$\Lambda$%
.

Le 
$\Lambda$%
-module {\`a} droite 
{\it dual\/} 
de
$M$ 
est le  groupe ab{\'e}lien 
${\mathop{\rm Hom\/}\nolimits}_\Lambda\,(M, \Lambda)$
des formes lin{\'e}aires sur 
$M$ 
muni de la structure  usuelle%
\footnote{\small
dont l'homoth{\'e}tie de rapport
$\lambda$ 
est
$\delta\,(\lambda)\!: 
{\mathop{\rm Hom\/}\nolimits}_\Lambda(M, \Lambda)\!\rightarrow\!%
{\mathop{\rm Hom\/}\nolimits}_\Lambda\,(M, \Lambda),
\varphi\!\mapsto\!\delta\,(\lambda)\,(\varphi)\!=\!%
\varphi\lambda$%
, o{\`u} 
$\varphi\,\lambda$
est la forme lin{\'e}aire 
$x\mapsto\varphi\,(x)\,\lambda$%
.
}.
On le note 
$M^\ast$%
.
\rmc Exemple| L'application 
$\varphi \mapsto \varphi\,(1)$
est un isomorphisme du dual 
$\Lambda^\ast=\bigl(\Lambda_s\bigr)^\ast$ 
de l'anneau 
$\Lambda$
vu comme 
$\Lambda$%
-module {\`a} gauche  sur 
$\Lambda_d$%
, l'anneau 
$\Lambda$ 
vu comme 
$\Lambda$%
-module {\`a} droite\footnote{\small
L'isomorphisme r{\'e}ciproque {\'e}tant 
$\mu \mapsto [\varphi_\mu :
\lambda\mapsto \lambda\,\mu]$%
.
}.
\finc

Si 
$M$ 
est libre de base 
$(b_i)_{i\in I}$ 
la
{\it forme lin{\'e}aire duale\/} 
du 
$k^{\hbox{\small i\`eme}}$
{\'e}l{\'e}ment de base 
$b_k$ 
est
$b_k^\ast : M\rightarrow \Lambda$ 
d{\'e}termin{\'e}e par 
$b_k^\ast\,(b_k) = 1$
et 
$b_k^\ast\,(b_i) = 0$ 
pour 
$i\ne k$%
.

La  Proposition 2 de {\bf 1.4} se traduit
en :

\thc Proposition 1| Le dual 
$(\oplus_{i\in I}M_i)^\ast$ 
de la somme directe des modules 
$M_i$ 
pour 
$i\in I$ 
s'identifie au produit
$\prod_{i\in I}M_i^\ast$ 
des duaux de ces modules.

En particulier le dual d'un
$\Lambda$%
-module {\`a} gauche libre de base finie
$(b_1,\,\ldots,\,b_n)$ 
est le 
$\Lambda$%
-module {\`a} droite libre de base la base duale 
$(b_1^\ast,\,\ldots,\,b_n^\ast)$%
.
\finc

%
\Defns La {\it transpos{\'e}e\/} d'une application 
$\Lambda$%
-lin{\'e}aire
$f\!:\!M\!\rightarrow\!N$
est l'application 
$f^\ast\!:\!N^\ast\!\rightarrow\!M^\ast$ 
qui {\`a} la forme lin{\'e}aire 
$\psi : N\rightarrow \Lambda$ 
associe la forme lin{\'e}aire
$f^\ast\,(\psi) = \psi\!\circ\!f : M\rightarrow N\rightarrow \Lambda$%
.

L'application
${\mathop{\rm t\/}\nolimits}: {\mathop{\rm Hom\/}\nolimits}_\Lambda(M, N)\!\rightarrow\!%
{\mathop{\rm Hom\/}\nolimits}_\Lambda%
(N^\ast, M^\ast),\,
f\mapsto {\rm t\/}\,(f)\!=\!f^\ast$
est un morphisme de groupe ab{\'e}lien et,
si les modules
$M$
et
$N$
sont munis de structure bilat\'erales, par exemple la canonique si l'anneau 
$\Lambda$ 
est commutatif, un morphisme de
$\Lambda$%
-module. C'est la 
{\it transposition\/}.

\thc Proposition 2| La transpos{\'e}e de 
l'identit{\'e} d'un module 
$M$
est l'identit{\'e} de son dual
$M^\ast$%
~:
$${\mathop{\rm t\/}\nolimits}({\rm Id\/}_M) =
{\mathop{\rm Id\/}\nolimits}_{M^\ast}$$

Si 
$g\in {\mathop{\rm Hom\/}\nolimits}_\Lambda\,(L, M)$ 
et 
$f\in {\mathop{\rm Hom\/}\nolimits}_\Lambda\,(M, N)$
sont composables alors~:
$${\mathop{\rm t\/}\nolimits}(f\!\circ\!g) = {\mathop{\rm t\/}\nolimits}(g)\!\circ\!{\mathop{\rm t\/}\nolimits}\,(f) \in {\mathop{\rm Hom\/}\nolimits}_\Lambda(N^\ast, L^\ast)$$
\finc

\thc Lemme et d{\'e}finitions| Soit
$m\in M$
un \'el\'ement du
$\Lambda$%
-module {\`a} gauche $M$%
.

L'application
$e_m\!: \!M^{\ast}\!\rightarrow\!\Lambda, \varphi\!\mapsto\!\varphi\,(m)$
est  un \'el\'ement de
$(M^{\ast})^{\ast}$%
, dit
{\it {\'e}valuation en\/} 
$m$%
, et l'application
$M\rightarrow(M^{\ast})^{\ast}, m\mapsto e_m$
est
$\Lambda$%
-lin{\'e}aire {\`a} gauche.

C'est l'%
{\it application canonique\/}
de 
$M$ 
dans son 
{\it bidual\/} 
$M^{\ast\ast}=(M^\ast)^\ast$%
, not{\'e}e
$\iota_M : M\rightarrow M^{\ast\ast}$%
.
\finc
La  Proposition 1 et l'%
{\sl Exemple\/} 
donnent :

\thc Proposition 3| Si 
$P$ 
est  facteur direct d'un module libre 
$L=P\oplus Q$
alors~:

$({\romannumeral1})$ l'application canonique 
$\iota_P : P\rightarrow P^{\ast\ast}$
est injective.

$({\romannumeral2})$ si 
$L$
est de type fini,
$\iota_P$
est un isomorphisme de 
$P$
sur son bidual 
$P^{\ast \ast}$%
.
\finc
\rmc Remarque|
Le {\petcap Corollaire\/} de la {\petcap Proposition 3\/} de {\bf 1.4\/}
vaut\footnote{\small
Si
$g_{P}: P\rightarrow N$
il suffit de d\'efinir
$g: L\!=\!P\oplus Q\!\rightarrow\!M, g(p, q)\!=\!(g_{P}(p), 0)$%
.
\finc
}pour un tel
$P$%
.
\finc
}
\vfill\eject
\npage
\nsecp Calcul matriciel et d{\'e}terminants.|
\nssecp Calcul matriciel.|
Pour un entier positif
$k\/$ 
on note 
$(x_1,\,\ldots,\,x_k)^t$ 
la colonne ayant l'{\'e}l{\'e}ment 
$x_i$
de l'anneau 
$\Lambda\/$ 
dans sa 
$i^{\hbox{\small i\`eme}}$
 ligne, c'est la \og transpos{\'e}e\fg 
de la ligne 
$(x_1,\,\ldots,\,x_k)$%
. 

Soit
$\Lambda^k\/$
le 
$\Lambda\/$%
-module {\`a} droite de ces colonnes de hauteur 
$k\/$ 
{\`a} coefficients dans l'anneau 
$\Lambda\/$
que l'on consid{\`e}re comme somme directe de 
$k\/$ 
exemplaires 
$\Lambda_i=e_i\,\Lambda\/$ 
de l'anneau
$\Lambda\/$
o{\`u} 
$e_1,\,\ldots,\,e_k$ 
est\footnote{\small
l'{\'e}l{\'e}ment
$e_j=(\delta_{j,\,i})^t$ 
est la colonne ayant 
$1$
dans la 
$j^{\hbox{i\`eme}}$
ligne et 
$0\/$ 
dans les autres.
} la 
{\it base canonique\/} 
de 
$\Lambda^k\/$%
.

Pour 
$f,g\!\in\!{\mathop{\rm End\/}\nolimits}%
_\Lambda\,(\Lambda_d)\!=\!%
{\mathop{\rm Hom\/}\nolimits}_\Lambda%
(\Lambda_d, \Lambda_d)$ 
deux endomorphismes de l'anneau 
$\Lambda\/$%
, vu comme 
$\Lambda\/$%
-module {\`a} droite, on a 
$f\!\circ g\!\,(1)\!=\!f\,(g\,(1))\!=\!%
f\,(1\,g\,(1))\!=\!f\,(1)\,g\,(1)$%
.

 Ainsi l'application 
$f\!\mapsto\!f(1)$ 
est un isomorphisme de l'anneau des endomorphismes 
$\Lambda\/$%
-lin{\'e}aires {\`a} droite de 
$\Lambda\/$ 
sur l'anneau 
$\Lambda\/$%
, (alors que l'anneau des endomorphismes 
$\Lambda\/$%
-lin{\'e}aires {\`a} gauche de 
$\Lambda$
est isomorphe
{\`a} l'anneau 
$\Lambda^o\/$ 
oppos{\'e} {\`a} 
$\Lambda\/$%
!).

 Ainsi, si 
$f : \Lambda^m\!\rightarrow\!\Lambda^n$
est une application 
$\Lambda\/$%
-lin{\'e}aire {\`a} droite, sa matrice 
d{\'e}crite en {\bf 1.4\/},
s'identifie {\`a} une matrice
$A\,(f)=A=
(a_{i,\,j\/})_{\build{1\leq i\leq n}_{1\leq j\leq m}^{}}$
{\`a} 
$n\/$ 
lignes et 
$m\/$ 
colonnes, on dit aussi 
$n\times m\/$%
, d'{\'e}l{\'e}ments 
$a_{i,\,j}={\rm pr\/}_i\,(f\,(e_j))\/$ 
de l'anneau
 $\Lambda\/$%
. 

Comme 
$\Lambda^n\/$ 
est aussi, par multiplication {\`a} gauche 
d'une colonne par un {\'e}l{\'e}ment de l'anneau,
$(\lambda,\, (x_1,\ldots,\,x_n)^t)\mapsto
(\lambda\, x_1,\ldots,\,\lambda\,x_n)^t)\/$%
, un 
$\Lambda\/$
module {\`a} gauche, 
${\mathop{\rm Hom\/}\nolimits}%
(\Lambda^m,\,\Lambda^n)\/)$ 
est, par multiplication {\`a} gauche au but
$(\lambda,\, f)\mapsto 
\lambda\,f : x\mapsto \lambda\,f\,(x)\/$%
, un
$\Lambda\/$%
-module {\`a} gauche et la matrice de 
$\lambda\,f\/$ 
est
$A\,(\lambda\,f)=\lambda\,A=
(\lambda\,a_{i,\,j\/})_%
{\build{1\leq i\leq n}_{1\leq j\leq m}^{}}\/$.
\Defns La 
{\it matrice \/}
$\!\!A(f)$ 
d'une application 
$\Lambda\/$%
-lin{\'e}aire {\`a} droite 
$f : \Lambda^m\rightarrow \Lambda^n$ 
est la matrice 
$n\!\times\!m$
{\`a} coefficients dans 
$\Lambda\/$ 
dont les 
$m\/$ 
colonnes sont les images par
$f\/$ 
des {\'e}l{\'e}ments 
$e_1,\,\ldots,\,e_m\/$
de la base canonique  de 
$\Lambda^m\/$%
.

Ainsi l'image par 
$f\/$ 
d'une colonne 
$X=(x_1,\,\ldots,\,x_m)^t$ 
de hauteur 
$m\/$
est la colonne 
$A(f)\cdot X=Y=(y_1,\,\ldots,\,y_n)^t$ 
de hauteur $n\/$ o{\`u}, pour 
$i=1,\ldots,n\/$ 
:
$$y_i=\sum_{j=1}^{m}a_{i,\,j}x_j$$

Et si 
$g : \Lambda^l\rightarrow \Lambda^m$
est une application 
$\Lambda\/$%
-lin{\'e}aire {\`a} droite de matrice
$B=
(b_{j,\,k\/})_%
{\buildrel{1\leq j\leq m}\over{1\leq k\leq l}}$%
, l'application compos{\'e}e 
$h= f\!\circ\!g :\Lambda^l\rightarrow \Lambda^n$
est de matrice la 
{\it matrice produit\/}
$A\cdot B=C=
(c_{i,\,k\/})_%
{\buildrel{1\leq i\leq n}\over{1\leq k\leq l}}$
o{\`u}, pour 
$i= 1,\dots,n$ 
et 
$k=1,\ldots,l\/$ :
$$c_{i,\,k}=\sum_{j=1}^m a_{i,\,j}\,b_{j,\,k}$$

On note 
${\rm M\/}_{n,\,m}\!(\Lambda)\/$ 
le 
$\Lambda\/$%
-module {\`a} gauche des matrices 
$n\!\times\!m\/$
{\`a} coefficients dans 
$\Lambda\/$ 
et
${\rm M\/}_{n}(\Lambda)\!=\!%
{\rm M\/}_{n,\,n}(\Lambda)\/$
celui des matrices carr{\'e}es 
$n\!\times\!n$ 
{\`a} coefficients dans 
$\Lambda\/$%
.

Le dual 
$(\Lambda^k)^\ast\/$ 
du 
$\Lambda\/$%
-module {\`a} droite
$\Lambda^k\/$ 
des colonnes
de hauteur 
$k\/$ 
est le module {\`a} gauche des matrices 
lignes de longueur 
$k\/$%
. Si une application 
$\lambda\/$%
-lin{\'e}aire {\`a} droite
$f : \Lambda^m\rightarrow \Lambda^n\/$ 
est de matrice 
$A\/$ 
l'application transpos{\'e}e 
$f^\ast : 
(\Lambda^n)^\ast\rightarrow (\Lambda^m)^\ast\/$
est la multiplication {\`a} droite 
des matrices lignes de longueur
$n\/$
par la matrice 
$A\/$%
.

Ces notions se traduisent dans le cas 
de modules libres munis de base :
\Defn Si 
$M\/$ 
et 
$N\/$ 
sont des 
$\Lambda\/$%
-modules {\`a} droite libres de bases respectives
${\cal C\/}=
(c_j)_{j\in J}\/$ et ${\cal D\/}=(d_i)_{i\in I}$ 
et
$f : M\rightarrow N$ 
est une application
$\Lambda\/$%
-lin{\'e}aire la 
{\it matrice de\/} 
$f$ 
{\it dans\/}
les bases 
${\cal C\/}$ et ${\cal D\/}$
est la matrice\footnote{\small
L'{\'e}criture dans cette notation de la base
${\cal C\/}$
de la source {\`a} droite est
pour des raisons typographiques~:
elle est du m{\^e}me c{\^o}t{\'e}
que la variable de source
$x$
dans l'{\'e}criture 
$f(x)$
de l'{\'e}valuation
de l'application lin{\'e}aire
$f$
en
$x$%
, et elle rendra plus mn{\'e}mothechnique la
formule pour la matrice d'une compos{\'e}
d'applications lin{\'e}aire ci-dessous.
} 
$A^{{\cal D\/},\,{\cal C\/}}\,(f)=
A\,(f) = (a\,(f)_{i,\,j})$
de 
$\varphi_N\!\circ\!f\!\circ\!\varphi_M^{-1} :
L_J\rightarrow L_I\/$
o{\`u} 
$\varphi_M : L_J\rightarrow M\/$ 
et 
$\varphi_N : L_I\rightarrow N\/$
sont les
isomorphismes d{\'e}finis par 
$\varphi_M\,(e_j)=c_j\/$ 
et 
$\varphi_N\,(e_i)=d_i\/$.
Pour tout 
$j\in J$ de $J$ 
on a la relation :
$$f\,(c_j)=\sum_{i\in I}d_i\, a\,(f)_{i,\,j}$$
qui caract{\'e}rise la matrice 
$A\,(f)=A^{{\cal D\/},\,{\cal C\/}}$
de l'application 
$f\/$
dans les bases 
${\cal C\/}$ et ${\cal D\/}$%
.

Si 
$g : L\rightarrow M$ 
est une application 
$\Lambda\/$%
-lin{\'e}aire {\`a} droite d'un 
$\Lambda\/$%
-module libre 
$L\/$ 
de base 
${\cal B\/}=(b_k)_{k\in K}$ 
alors on a la relation ~:
$$A^{{\cal D\/},\,{\cal B\/}}\,(f\!\circ\!g)=
A^{{\cal D\/},\,{\cal C\/}}\,(f)\cdot%
 A^{{\cal C\/},\,{\cal B\/}}\,(g)$$

Si l'anneau 
$\Lambda\/$ 
est commutatif le
$\Lambda\/$%
-module {\`a} gauche 
$(\Lambda^k)^\ast={\rm M\/}_{1\times k}\/$ 
des lignes de longueur 
$k\/$ 
est isomorphe par \og transposition\fg
$(\lambda_1,\ldots ,\lambda_k)\mapsto 
(\lambda_1,\ldots ,\lambda_k)^t\/$
au 
$\Lambda\/$%
-module {\`a} droite
$\Lambda^k={\rm M\/}_{k\times 1}\/$ 
des colonnes de hauteur 
$k\/$
d'o{\`u} la :
\Thc Proposition| Si l'anneau 
$\Lambda\/$ 
est commutatif et
$f : M\rightarrow N\/$ 
est une application 
$\Lambda\/$%
-lin{\'e}aire d'un 
$\Lambda\/$%
-module libre 
$M\/$ 
de base finie
${\cal B\/}\!=\!(b_1,\ldots, b_m)\/$
sur un 
$\Lambda\/$%
-module libre 
$N\/$ 
de base finie
${\cal C\/}\!\!=\!\!(c_1,\!\ldots\!,c_n)\/$
alors la transpos{\'e}e 
$f^\ast\!\!: N^\ast\!\!\rightarrow\!\!M^\ast\/$
est de matrice,
dans les bases duales 
${\cal C\/}^\ast\!=\!(c_1^\ast,\ldots ,\!c_n^\ast)\/$ 
et 
${\cal B\/}^\ast\!=\!(b_1^\ast,\ldots ,b_m^\ast)\/$%
,  la matrice obtenue en {\'e}changeant lignes et colonnes de
$A$%
, la
{\it matrice transpos{\'e}e\/}~:
$$A^t=(a_{j, i})_%
{\build{1\leq i\leq n}_{1\leq j\leq m}^{}}\/$$
de la matrice
$A=(a_{i, j})_%
{\build{1\leq i\leq n}_{1\leq j\leq m}^{}}\/$
de l'application 
$f\/$ 
dans les bases 
${\cal B\/}$ et ${\cal C\/}$%
.

Ainsi si 
$A\in {\rm M\/}_{n,\,m}\/$ 
et 
$B\in {\rm M\/}_{m,\,l}\/$ 
sont deux matrices composables alors dans 
${\rm M\/}_{l,\, n}\/$
on a la relation~:
$$(A\,B)^t=B^t\cdot A^t$$
\finc
\vfill\eject
\npage
\nssecp Quelques identit{\'e}s polynomiales.|
\nsssecc 
M{\'e}thode de Gauss%
\footnote{\small
appel{\'e}e fang cheng 
(mod{\`e}le rectangulaire) par  Chang Ts'ang
(
$\uppercase\expandafter{\romannumeral2}^{\hbox{i\`eme}}$
sci{\`e}cle avant notre {\`e}re)
{\`a} qui, selon P. Gabriel ({\bf [Ga]\/}), il faut attribuer
cet algorithme de r{\'e}solution des syst{\`e}mes
lin{\'e}aires
}
pour le 
\og syst{\`e}me 
$n\!\times\!n\/$
 g{\'e}n{\'e}ral\fgf|
{Soit}
\vskip-3mm
$$\widetilde{\Lambda}=\widetilde{\Lambda}_n=
{\Bbb Z\/}\,[a_{i,\,j}]_{1\leq i,\, j\leq n}\/$$
l'anneau  des polyn{\^o}mes 
{\`a} coefficients entiers en 
$n^2\/$ 
variables
$a_{i,\,j}, 1\leq i,\, j\leq n\/$%
. Cet anneau est factoriel, de corps des fractions
$\widetilde{K}=
{\Bbb Z\/}\,(a_{i,\,j})_{1\leq i,\, j\leq n}\/$ 
le corps des fractions rationnelles en ces 
$n^2\/$ 
variables 
$a_{i,\,j}\/$%
. On notera
$\widetilde{\Lambda}[Y]\/$ et $\widetilde{K}[Y]$ 
les anneaux des polyn{\^o}mes {\`a} coefficients dans 
$\widetilde{\Lambda}\/$ 
et 
$\widetilde{K}\/$ 
respectivement en 
$n\/$ 
variables 
$Y_1,\,\ldots,\, Y_n\/$
que l'on consid{\`e}re comme une colonne 
$Y=(Y_1,\,\ldots,\,Y_n)^t\/$ 
de variables.

Si 
$p\!\in\!\widetilde{\Lambda}\/$
est un polyn{\^o}me en les 
$n^2\/$ 
variables 
$a_{i,\,j}\/$%
, le polyn{\^o}me obtenu en rempla\c cant dans 
$p\/$ 
la 
$k^{\hbox{i\`eme}}$
 colonne de variables par la colonne
$Y=(Y_1,\ldots,Y_n)^t \/$%
, c'est {\`a} dire en substituant dans 
$p\/$%
, pour 
$i=1,\ldots, n\/$
la variable 
$Y_i\/$ 
{\`a} la variable 
$a_{i,\,k}\/$ 
est not{\'e}
$p\,(A_{\bullet,\, k}=Y)\/$%
.

La 
{\it matrice g{\'e}n{\'e}rale\/} 
$n\!\times\!n\/$ 
est la matrice
$A=(a_{i,\,j})_{\build{1\leq i\leq n}_{1\leq j\leq n}^{}}
\in 
{\mathop{\rm M}\nolimits}_n\,(\widetilde{\Lambda})\/$%
. 

Le 
{\it syst{\`e}me g{\'e}n{\'e}ral\/} 
$n\!\times\!n\/$ 
est le syst{\`e}me lin{\'e}aire 
$A\cdot X=Y\/$
dans 
$\widetilde{\Lambda}[Y]\/$ :
$$\left\lbrace
  \begin{matrix}
a_{1,\,1}\,X_1 & \ +\ & a_{1,\,2}\,X_2 & \ldots\,\ +\ & %
a_{1,\,i}\,X_i & \ldots\,\ +\  & a _{1,\,n}\,X_n & \ =\ & Y_1\\
a_{2,\,1}\,X_1 & \ +\ & a_{2,\,2}\,X_2 & \ldots\,\ +\ & %
a_{2,\,i}\,X_i & \ldots\,\ +\ & a_{2,\,n}\,X_n & \ =\ & Y_2\\
\ \vdots\hfill & & \ \vdots\hfill & & \ \vdots\hfill & %
\vdots & \ \vdots\hfill & \vdots & \ \vdots\hfill\\
a_{n,\,1}\,X_1 &\ +\ & a_{n,\,2}\,X_2 & \ldots\,\ +\ & %
a_{n,\,i}\,X_i & \ldots\,\ +\ & a_{n,\,n}\,X_n & \ =\ & Y_n
  \end{matrix}
  \right.
  \leqno(\Sigma)$$

On d{\'e}finit, pour 
$0\leq g\leq n$
les matrices
$A^g\!=\!%
(a^g_{i,\,j})_{\build{1\leq i\leq n}_{1\leq j\leq n}^{}}
\in 
{\mathop{\rm M\/}\nolimits}_n\,(\widetilde{\Lambda})\/$
par
$A^0\!=\!A$
 et, en posant
$a^0_{0,\,0}\!=\!1$
et utilisant le symbole de Kronecker
$\delta_{g,\,g}\!=\!1, \delta_{i, g}\!=\!0$
si
$i\!\ne\!g$%
,
 les relations de r{\'e}currence si 
$0<g$%
~:
$$a^{g}_{i,\,j}=
a^{g-1}_{g,\,g}\,a^{g-1}_{i,\,j}-
a^{g-1}_{i,\,g}\,a^{g-1}_{g,\,j}+
\delta_{i,\, g}a^{g-1}_{g-1,\,g-1}a^{g-1}_{i,\, j}$$
ainsi
$a^{g}_{g,\,j}=a^{g-1}_{g-1,\,g-1}a^{g-1}_{g,\, j}$
et, si
$i\ne g,
a^{g}_{i,\,j}=
a^{g-1}_{g,\,g}\,a^{g-1}_{i,\,j}-
a^{g-1}_{i,\,g}\,a^{g-1}_{g,\,j}$
donc
$a^{g}_{i,\,g}=0$%
~:

\Thc Lemme|({\romannumeral1})  Pour 
$j\!\leq\!g$
et
$i\!\ne\!j$
on a~:
$a^{g}_{j,\,j}\!=\!a^{g}_{g,\, g}\!=\!a^{g-1}_{g-1,\,g-1}a^{g-1}_{g,\, g}$
et
$a^{g}_{i,\,j}\!=\!0$

({\romannumeral2}) Les polyn{\^o}mes 
$a^g_{i,\,j}=B^g_{i,\,\bullet}\cdot A_{\bullet,\,j}\/$ 
sont lin{\'e}aires\footnote{\small
o{\^u} la ligne
$B^g_{i,\,\bullet}\in\tilde{\Lambda}^n$
(ne contenant aucun
$a_{i,\,j}$
pour
$1\leq i\leq n$%
)
exprime cette lin{\'e}arit{\'e}.
} en les variables de la 
$j^{\hbox{i\`eme}}$
 colonne, ne d{\'e}pendent que des variables 
$a_{r,\,s}\/$ 
avec 
$r\leq i,\,s\leq j\/$  
et sont homog{\`e}nes\footnote{\small
de degr{\'e}, \`a ordre
$g$
et indice de ligne
$i$
fix\'es, ind\'ependants de l'indice
$j$
de colonne.
}%
.

({\romannumeral3}) Pour 
$g=0,\,\ldots,\,n\/$ 
on a 
$A^g({\mathop{\rm Id}\nolimits})=
{\mathop{\rm Id}\nolimits}\/$%
.

\finc 
Pour
$0\leq g\leq n$
d{\'e}finissons la colonne
$Y^g\in \widetilde{\Lambda}[Y]^n$
 de composantes
$$Y^g_i=a^g_{i,\,n}(A_{\bullet,\,n}=Y)=
B^g_{i,\,\bullet}\cdot Y$$
et notons
$B^g$
la matrice dont la
$i^{\hbox{i\`eme}}$%
 ligne est
$B^g_{i,\,\bullet}$
donn{\'e}e par 
({\romannumeral2}) du Lemme.

{\small
La matrice
$A^g$
s'obtient {\`a} partir de
$A^{g-1}$
en multipliant sa
$g^{\hbox{i\`eme}}$
ligne par
$a^{g-1}_{g-1,\,g-1}$
et en retranchant aux autres 
$i^{\hbox{i\`eme}}$
lignes multipli{\'e}s par
$a^{g-1}_{g,\,g}$%
, leur
$g^{\hbox{i\`eme}}$
 ligne multipli{\'e} par
$a^{g-1}_{i,\,g}$%
. Comme\footnote{\small
selon
({\romannumeral1}) et ({\romannumeral3})
du Lemme.
}
$a^{g-1}_{g-1,\,g-1}\ne0\ne a^{g-1}_{g,\,g}\in%
 \widetilde{\Lambda}\subset\widetilde{K}$%
, le syst{\`e}me
$A^0X=AX=Y$
est, dans 
$\widetilde{K}[Y]$%
, {\'e}quivalent aux
 syst{\`e}mes interm{\'e}diaires
$A^{g-1}X=Y^{g-1}, g=2,\ldots, n$
qui, en posant
$d_g=a^g_{g,\,g}\in \widetilde{\Lambda}$%
, sont~:

$$\left\lbrace
\begin{matrix}
d_{g-1}\,X_1&\cdots+\,0\cdots&%
+\,0&\!\!\cdots+\,0\cdots&%
+a^{g-1}_{1,\,g}\,X_{g}&\!\!+\cdots&
+a^{g-1}_{1,\,n}\,X_n&=\ Y^{g-1}_1\\
\noalign{\vskip 2mm}
\build{\hfill 0 \hfill}_{\vdots}^{\vdots}\quad\ &
\!\!\build{\ \ \ \hfill\, +\, 0\, \hfill \ \ \ \ \ \ }_%
{^{\ \ \ \ \ \ \ddots}}^{\ddots \ \ \ \ \ \ \ \ }&
\build{+\, d_{g-1}\,X_k}_%
{\build{\hfill  \hfill}_{\vdots}^{}}^
{\build{\hfill  \hfill}_{}^{\vdots}}&
\!\!\build{\ \ \ \hfill\, +\, 0\, \hfill \ \ \ \ \ \ }_%
{^{\ \ \ \ \ \ \ddots}}^{\ddots \ \ \ \ \ \ \ \ }&
\build{+ a^{g-1}_{k,\,g}\,X_g}_{\vdots}^{\vdots}&%
\!\!+\cdots&
\build{a^{g-1}_{k,\,n}\,X_n}_{\vdots}^{\vdots}&
=\ \build{Y^{g-1}_k}_{\vdots}^{\vdots}\\
\noalign{\vskip 3mm}
\hfill 0 \hfill&\!\!\cdots+\,0\cdots&+\,0&%
\!\!\cdots+\,0\cdots&
+a^{g-1}_{g,\,g}\,X_{g}&\!\!+\cdots&
+a^{g-1}_{g,\,n}\,X_n&=\ Y^{g-1}_g\\
\vdots&\vdots &\hfill\vdots\ &
\vdots&
\vdots&\vdots&\vdots&\vdots\cr
0&\!\!\cdots+0\cdots&+0&\!\!\cdots+0\cdots&
+a^{g-1}_{n,\,g}\,X_{g}&\!\!+\cdots&
+a^{g-1}_{n,\,n}\,X_n&=\ Y^{g-1}_n
\end{matrix}
\right.
\leqno(\Sigma_{g-1})$$
et le syst{\`e}me 
$AX=Y$
est aussi {\'e}quivalent au dernier syst{\`e}me
$A^nX=Y^n$%
, qui, {\'e}tant  diagonal~:
$$d_nX_i=Y^n_i,\quad i=1,\ldots, n 
\leqno(\Sigma_{n})$$
et comme 
$Y^n_i=B_i^n\cdot Y$
est lin{\'e}aire en
$Y$%
,  a dans
$\widetilde{K}[Y]$
les solutions lin{\'e}aires en
$Y$%
~:
$$X_i=(d_n^{-1}B^n_i)\cdot Y,\quad i=1,\ldots, n
\leqno(S)$$
 Comme
$Y=A\cdot X$%
, on a
$Y=A\cdot\bigl((d_n^{-1}B^n)\cdot Y\bigr)=
\bigl(A\cdot(d_n^{-1}B^n)\bigr)\cdot Y$
et
$X=(d_n^{-1}B^n)\cdot Y=
(d_n^{-1}B^n)\cdot (A\cdot X))=
((d_n^{-1}B^n)\cdot A)\cdot X$
d'o{\`u}, dans
${\mathop{\rm M}\nolimits}_n\,(\widetilde{K})\/$%
, les relations~:
$$A\cdot (d_n^{-1}B^n)=
{\mathop{\rm Id}\nolimits}_n=
(d_n^{-1}B^n)\cdot A$$
traduisant que la matrice g{\'e}n{\'e}rale
$A$
est inversible dans 
${\mathop{\rm M}\nolimits}_n\,(\widetilde{K})\/$%
, d'inverse
$(d_n^{-1}B^n)$%
~:
}

\Thc Corollaire|La matrice g{\'e}n{\'e}rale 
est inversible dans 
${\mathop{\rm M\/}\nolimits}_n\!(\widetilde{K})\/$%
.\hfill\break
Plus pr{\'e}cis{\'e}ment~:
Il y a une matrice 
$B^n\!\in\! %
M_n(\widetilde{\Lambda})$
{\`a} coefficients homog{\`e}nes et un polyn{\^o}me 
$d_n\/$ 
homog{\`e}ne et lin{\'e}aire
en la derni{\`e}re colonne, tels que 
$d_n\,({\rm Id\/})\!=\!1\/$ 
et
$$B^n\cdot A=
d_n{\mathop{\rm Id}\nolimits}_n=
A\cdot B^n$$
De plus la derni{\`e}re
coordonn{\'e}e $x_n\/$ de la solution 
$x\in \widetilde{K}[Y]^n\/$ de $A\, X=Y\/$ 
v{\'e}rifie :
$$d_n x_n=d_n(A_{\bullet,\,n}=Y)$$

\finc
\npage
\nsssecp %
Les identit{\'es} remarquables de Crammer et %
du d{\'e}terminant.|
L'anneau 
$\widetilde{\Lambda}\/$ 
{\'e}tant factoriel d'unit{\'e}s 
$\{-1,\,1\}\/$%
, tout 
$M\!\in\!{\mathop{\rm M\/}\nolimits}_n\,(\widetilde{K})\/$ 
d{\'e}finie\footnote{\small
{\it c. a. d.} en
substi\-tuant
aux coefficients de la matrice g{\'e}n{\'e}rale 
$A\/$ 
ceux de la matrice identit{\'e}.
} 
en 
${\mathop{\rm Id\/}\nolimits}\/$
a une unique expression r{\'e}duite
$M\!=\!\frac{N}{d}\/$ 
de num{\'e}rateur 
$N\!\in\!{\mathop{\rm M\/}\nolimits}_n\,(\widetilde{\Lambda})\/$
une matrice {\`a} coeficients dans l'anneau 
$\widetilde{\Lambda}\/$
et d{\'e}nominateur
$d\!\in\!\widetilde{\Lambda}\/$ 
 un polyn{\^o}me {\`a} coefficients entiers
v{\'e}rifiant 
$d\,({\mathop{\rm Id\/}\nolimits})>0\/$ 
et tels que 
$N\/$ 
et 
$d\/$ 
sont sans facteurs communs.
\Defns Le 
{\it d{\'e}terminant g{\'e}n{\'e}ral\/} 
${\mathop{\rm det\/}\nolimits}(A)\/$ 
et la
{\it comatrice g{\'e}n{\'e}rale\/}
$\widetilde{A}\/$ 
de la matrice g{\'e}n{\'e}rale 
$A\/$ 
sont les d{\'e}nominateurs
et num{\'e}rateurs de l'expression r{\'e}duite
$\frac{\widetilde{A}}{{\rm det\/}\,A}=
\frac{B^n}{d_n}=A^{-1}\/$ 
de l'inverse 
$A^{-1}\/$
de la matrice g{\'e}n{\'e}rale 
$A\/$%
. D'apr{\`e}s  le {
Corollaire\/}
 de {\bf 2.2.1\/}
les coefficients de la matrice
$\widetilde{A}\/$ 
et 
${\mathop{\rm det\/}\nolimits}A\/$ 
sont  homog{\`e}nes et~:
$${\mathop{\rm det\/}\nolimits}%
({\mathop{\rm Id\/}\nolimits})=1
\leqno{(0)}$$
$$\widetilde{A}\cdot A={\rm det\/}(A)\,{\rm Id\/}=
A\cdot \widetilde{A}
\leqno{(1)}$$

\Thc Identit{\'e}s remarquables du d{\'e}terminant| 
Le d{\'e}terminant g{\'e}n{\'e}ral 
est multiplicatif pour composition et 
forme triangulaire par blocs :

$({\romannumeral 1})$
 Si 
$B\!\in\!{\rm M\/}_n\,(\widetilde{\Lambda}'_n),\
C\!\in\!{\rm M\/}_n\,(\widetilde{\Lambda}''_n)\/$ 
sont deux matrices g{\'e}n{\'e}rales 
$n\!\times\!n$
on a~:

$${\rm det\/}\,(B\cdot C)=
{\rm det\/}\,(C)\,{\rm det\/}\,(B)
\leqno{(2)}$$

$({\romannumeral2})$
Si 
$U\!\in\!{\rm M\/}_p\,(\widetilde{\Lambda}_p), %
V\!\in\!{\rm M\/}_q\,(\widetilde{\Lambda}'_q), %
W\!\in\!{\rm M\/}_{p,\, q}\,%
(\widetilde{\Lambda}_{p,\,q})%
\/$
sont trois matrices g{\'e}n{\'e}rales, respectivement de taille 
$ p\!\times\!p,\ q\!\times\!q\/$
et 
$p\!\times\!q\/$
on a~:

$${\rm det\/}\,(%
\begin{pmatrix}
  U &W\cr
  0 &V
\end{pmatrix}%
)=
{\mathop{\rm det\/}\nolimits}(U)\,
{\mathop{\rm det\/}\nolimits}(V),\quad
\widetilde{\begin{pmatrix}
U &W\\
0 &V
\end{pmatrix}
}=
\begin{pmatrix}
  {\rm det\/}(V)\widetilde{U} &%
-\widetilde{U}\cdot W\cdot\widetilde{V}\\
0 &{\rm det\/}(U)\widetilde{V}
\end{pmatrix}
\leqno{(3)}$$
\finc
{\small
\Dem
Soit 
$\widetilde{\Lambda}(2)=\widetilde{\Lambda}_n(2)=
{\Bbb Z\/}\,[b_{i,\,j}\,c_{r,\,s}]_%
{1\leq i,\,j,\,r,\,s\leq n}\/$
l'anneau des polyn{\^o}mes {\`a}
coeficients entiers en les 
$2\,n^2\/$ 
variables
$b_{i,\,j}, 1\leq i,\, j\leq n\/$ 
et 
$c_{r,\,s}, 1\leq r,\,s\leq n$ 
et
$\widetilde{K}(2)={\Bbb Z\/}\,
(b_{i,\,j},\,c_{r,\,s})_{1\leq i,\, j,\,r,\,s\leq n}\/$
le corps des fractions rationnelles en ces 
$2\,n^2\/$ 
variables. 
L'anneau 
$\widetilde{\Lambda}_2\/$ 
est factoriel,
contient deux exemplaires 
$\widetilde{\Lambda}'={\Bbb Z\/}\,[b_{i,\,j}]\/$ 
et
$\widetilde{\Lambda}''={\Bbb Z\/}\,[c_{r,\,s}]\/$
de l'anneau 
$\widetilde{\Lambda}\/$ 
tels que des {\'e}l{\'e}ments
$x'\in \widetilde{\Lambda}'\/$ 
et 
$x''\in \widetilde{\Lambda}''\/$
ne peuvent avoir que des constantes comme 
diviseurs communs.

Les deux matrices g{\'e}n{\'e}rales 
$B=(b_{i,\,j})_{1\leq i,\, j\leq n}\/$ 
et
$C=(c_{r,\,s})_{1\leq r,\,s\leq n}$ 
de 
${\rm M\/}_n\,(\widetilde{\Lambda}(2))\/$ 
{\'e}tant inversibles dans
${\rm M\/}_n\,(\widetilde{K}(2))\/$%
, leur produit l'est aussi et on a
$(B\cdot C)^{-1}=C^{-1}\cdot B^{-1}$%
. Ainsi
$${\mathop{\rm det\/}\nolimits}(C)\,
{\mathop{\rm det\/}\nolimits}(B)\,(B\cdot C)^{-1}=
({\mathop{\rm det\/}\nolimits}(C)\,C^{-1})\cdot
({\mathop{\rm det\/}\nolimits}(B)\,B^{-1})=
\widetilde{C}\cdot\widetilde{B}\/$$
est une matrice {\`a} coefficients dans 
$\widetilde{\Lambda}(2)\/$%
. D'autre part en substituant dans la matrice
g{\'e}n{\'e}rale 
$A\/$ 
les coefficients de la matrice produit 
$B\cdot C\/$
on obtient 
${\mathop{\rm det\/}\nolimits}(B\cdot C)\,
(B\cdot C)^{-1}=
\widetilde{B\cdot C}\/$%
. D'o{\`u}~:
$${\mathop{\rm det\/}\nolimits}(B\cdot C)\,
\widetilde{C}\cdot\widetilde{B}=
{\mathop{\rm det\/}\nolimits}(C)\,
{\mathop{\rm det\/}\nolimits}(B)\,%
\widetilde{B\cdot C}\/$$

Les polyn{\^o}mes 
${\mathop{\rm det\/}\nolimits}(C)\/$ 
et 
${\mathop{\rm det\/}\nolimits}(B)\/$
{\'e}tant premiers avec 
$\widetilde{C}\/$ et $\widetilde{B}\/$%
, leur produit
${\mathop{\rm det\/}\nolimits}(C)\,
{\mathop{\rm det\/}\nolimits}(B)\/$
divise 
${\mathop{\rm det\/}\nolimits}(B\cdot C)\/$ 
et, ces deux polyn{\^o}mes {\'e}tant homog{\`e}nes 
de m{\^e}me degr{\'e} par rapport aux
$b_{i,\,j}\/$ 
et aux 
$c_{r,\,s}\/$ 
et vallant 
$1\/$
si 
$B={\mathop{\rm Id\/}\nolimits}$
et 
$C={\mathop{\rm Id\/}\nolimits}$%
, ils sont {\'e}gaux.
D'o{\`u} la premi{\`e}re 
identit{\'e} remarquable.

La seconde s'obtient
de mani{\`e}re analogue 
 en consid{\'e}rant, si
$n=p+q$
le sous-anneau 
$\widetilde{\Lambda}_{n;\,p,\,q}\!=\!
{\Bbb Z\/}\,[u_{i,\,j}\!=\!a_{i,\, j}\, %
w_{i,\,k}\!=\!a_{i\,p+k},\,%
v_{r,\,s}\!=\!a_{p+r,\,p+s}]_%
{1\leq i,\,j\leq p, 1\leq k,\,r,\,s\leq q}\subset
{\Bbb Z\/}\,[a_{i,\,j}]_{1\leq i,\,j\leq n}\/$
 des polyn{\^o}mes {\`a}
coeficients entiers en les 
$p^2+pq+q^2\/$ 
variables
$u_{i,\,j},\, w_{i,\,k},\,v_{r,\, s}$
et\hfill\break
$\widetilde{K}_{n;\,p,\,q}={\Bbb Z\/}\,
(u_{i,\,j}\,w_{i,\,k}\,
v_{r,\,s})_{1\leq i,\,j\leq p, 1\leq r,\,s\leq q}\/$
son corps des fractions
et en utilisant que 
l'inverse de la matrice g{\'e}n{\'e}rale
triangulaire sup{\'e}rieure par blocs 
$p\!\times\!p, p\!\times\!q, q\!\times\!q$
g{\'e}n{\'e}rale est~:
$${\begin{pmatrix}
  U&W\\
  0&V
  \end{pmatrix}
}^{-1}=
{\begin{pmatrix}
    (u_{i,\,j})&(w_{i\, k})\\
        0&(v_{r,\,s})
\end{pmatrix}}^{-1}=
\begin{pmatrix}U^{-1}&-U^{-1}\cdot W\cdot V^{-1}\\
0&V^{-1}
\end{pmatrix}$$
}
\vskip-2mm
\rmc Remarque| D'apr{\`e}s le  Corollaire
de {\bf 2.2.1},
d{\'e}terminant et  comatrice de 
la matrice g{\'e}n{\'e}rale
$A\in
{\mathop{\rm M\/}\nolimits}_n%
(\widetilde{\Lambda})\/$ 
sont respectivement
${\mathop{\rm det\/}\nolimits} A=
\frac{b_n}{\delta}\/$ 
et
$\widetilde{A}=
\frac{B^n}{\delta}\/$ 
o{\`u}
$\delta\/$ 
est le p.g.c.d. 
v{\'e}rifiant 
$\delta\,({\rm Id\/})=1\/$
de 
$b_n\/$ 
et
$B^n\/$%
.
On aurait eu une expression analogue 
si on avait effectu{\'e} la m{\'e}thode
de Gauss g{\'e}n{\'e}rale en {\'e}crivant 
les colonnes dans un autre ordre\footnote{\small
Et en mettant les lignes dans le m{\^e}me ordre
pour ne pas changer la matrice
${\mathop{\rm Id\/}\nolimits}$
permettant de normaliser le d{\'e}terminant
g{\'e}n{\'e}ral par
${\mathop{\rm det\/}\nolimits}%
({\mathop{\rm Id\/}\nolimits})=1$%
.
}
.
\finc
Ainsi d'apr{\`e}s la fin du {\petcap Corollaire\/} 
de {\bf 2.2.1}  on a :

\Thc Formules de Crammer|
Si 
$(x_1,\,\ldots,\, x_n)^t\/$ 
est solution du syst{\`e}me g{\'e}n{\'e}ral 
$A\,X=Y\/$
alors pour 
$i=1,\ldots,\,n\/$ 
on a l'identit{\'e} :
$${\rm det\/}\,A\cdot x_i=
{\rm det\/}\,(A_{\bullet,\, i}=Y)
\leqno{(4)}$$
\finc
\vskip-2mm
{\small
La Remarque donne ausi que le 
d{\'e}terminant 
est lin{\'e}aire en chaque colonne
de variable.\hfil\break
Si
$1\!\leq\!i\ne\!j\!\leq\!n$
la matrice 
$A(A_{\bullet,\, i}\!=\!A_{\bullet,\,j})$
ayant deux colonnes {\'e}gales n'est pas
inversible donc~:
$${\mathop{\rm det\/}\nolimits}%
(A(A_{\bullet,\, i}=A_{\bullet,\,j})=0$$
ainsi la forme
$n$%
-lin{\'e}aire d{\'e}terminant g{\'e}n{\'e}ral
${\mathop{\rm det\/}\nolimits}:
{\mathop{M}\nolimits}_n%
(\widetilde{\Lambda})=
\Bigl(\widetilde{\Lambda}^n\Bigr)^n%
\rightarrow\widetilde{\Lambda}$
est altern{\'e}e d'o{\`u}~:
}

\Thc Expression du d{\'e}terminant|
 Le d{\'e}terminant g{\'e}n{\'e}ral est 
$n\/$%
-lin{\'e}aire altern{\'e}
en les colonnes de la matrice g{\'e}n{\'e}rale. 
V{\'e}rifiant
${\rm det\/}\,({\mathop{\rm Id\/}\nolimits})=1\/$%
, il vaut~: 
$${\rm det\/}\,A=\sum_{\sigma\in{\goth S\/}_n}
\epsilon\,(\sigma)a_{\sigma\,(1),\,1}\cdots%
 a_{\sigma\,(n),\,n} \leqno{(5)}$$
o{\`u}
$\epsilon\,(\sigma)\in\{-1,\,1\}\/$
est la signature de la permutation
$\sigma\in{\goth S}_n\/$ 
de 
$\{1,\ldots,\,n\}$%
.
\finc
Comme
$A\cdot\widetilde{A}=
{\mathop{\rm det\/}\nolimits}(A)\,
{\mathop{\rm Id\/}\nolimits}$
la solution du \og syst{\`e}me transpos{\'e}\fg
$X^t\,A=Y^t\/$
est
$X^t=Y^t\cdot\frac{\widetilde{A}}{\mathop{\rm det\/} A}$%
. Elle aurait aussi {\'e}t{\'e} obtenue 
par  m{\'e}thode de Gauss sur les colonnes de
$A$
au lieu de l'avoir effectu{\'e}e sur les lignes
comme en {\bf 2.2.1\/}.
Ainsi~:
\Thc D{\'e}terminant de la transpos{\'e}e|
 Le d{\'e}terminant g{\'e}n{\'e}ral est 
$n\/$%
-lin{\'e}aire altern{\'e} en les lignes
 et on a la relation~:
$${\rm det\/}\,(A^t)={\rm det\/}\,(A)
\leqno{(6)}$$
\finc

{\small
Si 
$A_{i,\,j}$ 
est la matrice 
$(n-1)\!\times\!(n-1)\/$ 
obtenue en effa\c cant la 
$\imath^{\hbox{i\`eme}}$
 ligne et la 
$\jmath^{\hbox{i\`eme}}$
 colonne de 
$A\/$%
, la matrice g{\'e}n{\'e}rale 
$n\!\times\!n$%
, les formules de Crammer, 
l'alter-multilin{\'e}arit{\'e} {\`a} la fois
en les colonnes et en les lignes et la 
multiplicativit{\'e} par forme triangulaire%
\footnote{\small
pour
$p=1, q=n-1$
et
$\begin{pmatrix}
  A&B\\
  0&D
\end{pmatrix}=
\begin{pmatrix}1&\ast\\
  0&A_{i,\,j}
  \end{pmatrix}$%
.
}
du d{\'e}terminant donnent :
}
\Thp Expression de la comatrice et 
d{\'e}vellopements du d{\'e}terminant.|
$$\widetilde{A}_{i,\,j}=
(-1)^{i+j}{\rm det\/}\,(A_{j,\,i})
\leqno{(7)}$$
$${\rm det\/}\,A=
\sum_{i=1}^n(-1)^{i+j}{\rm det\/}\,(A_{i,\,j})\,a_{i,\,j}
=
\sum_{j=1}^n(-1)^{i+j}a_{i,\,j}\,{\rm det\/}\,(A_{i,\,j})
\leqno{(8)}$$
\finp
%
\npage
\nsssecp
L'identit{\'e} remarquables de Cayley Hamilton.|
\Defn Le {\it polyn{\^o}me caract{\'e}ristique 
g{\'e}n{\'e}ral}\footnote{\small
ou {\it polyn{\^o}me caract{\'e}ristique} de la
matrice g{\'e}n{\'e}rale
$A\in{\mathop{\rm M\/}\nolimits}_n%
(\widetilde{\Lambda})$%
.}
est~:
$\chi_A(X)=
{\mathop{\rm det\/}\nolimits}%
(X{\mathop{\rm Id\/}\nolimits}-A)%
$
\rmc Remarque|Le terme constant de
$\chi_A(X)=X^n+X^{n-1}c_1+\cdots+c_n%
\in\widetilde{\Lambda}[X]$
est
$$c_n=\chi_A(0)=
{\mathop{\rm det\/}\nolimits}(-A)=
(-1)^n{\mathop{\rm det\/}\nolimits}(A)$$
.
\finc
\Thc Th{\'e}or{\`e}me| On a la relation de Cayley-Hamilton~:
$$\chi_A(A)=0$$
\finc

{\small
\Dem Les coefficients de la comatrice
$\widetilde{X{\mathop{\rm Id\/}\nolimits}-A}$
de 
$X{\mathop{\rm Id\/}\nolimits}-A$
{\'e}tant
$(\widetilde{X{\mathop{\rm Id\/}\nolimits}-A})_{i,\, j}\!=\!
{\mathop{\rm det\/}\nolimits}%
((X{\mathop{\rm Id\/}\nolimits}-A)_{j,\,i})$
de degr{\'e} au plus 
$n\!-\!1$%
,  
il y a des matrices
$S_i\!\!\in\!\!{\mathop{\rm M\/}\nolimits}_n%
(\widetilde{\Lambda}),%
i\!=\!0,\ldots, n-1$
tels que~:
$$\widetilde{X{\mathop{\rm Id\/}\nolimits}-A}%
=\sum_{i=0}^{n-1}S_iX^{n-1-i}=
S_0X^{n-1}+\cdots S_kX^{n-1-k}+\cdots+S_{n-1}$$
En explicitant la relation
$X^n+c_1X^{n-1}+\cdots+c_n=
{\mathop{\rm det\/}\nolimits}%
(X{\mathop{\rm Id\/}\nolimits}-A)=%
\widetilde{X{\mathop{\rm Id\/}\nolimits}-A}\cdot
(X{\mathop{\rm Id\/}\nolimits}-A)=%
(X^{n-1}S_0+X^{n-1}S_1+\cdots+S_{n-1})\cdot
(X{\mathop{\rm Id\/}\nolimits}-A)$
on obtient~:
$$S_0={\mathop{\rm Id\/}\nolimits},\
S_1- S_0\cdot A=c_1,\cdots,%
S_k- S_{k-1}\cdot A=c_k,\cdots,%
S_{n-1}- S_{n-2}\cdot A=c_{n-1},\
- S_{n-1}\cdot A=c_n$$
d'o{\`u}
$\displaystyle
\chi_A(A)=A^n+\sum_{i=1}^nc_iA^{n-1}=
S_0\cdot A^n+\bigl\{%
\sum_{k=1}^{n-1}(S_k- S_{k-1}\cdot A)\cdot A^{n-k}%
\bigr\}-
 S_{n-1}\cdot A
=S_0\cdot A^n+\bigl\{%
\sum_{k=1}^{n-1}S_k\cdot A^{n-k}-S_{k-1}\cdot A^{n-k+1}%
\bigr\}-
S_{n-1}\cdot A=
A^n+\sum_{k=1}^{n-1}
S_k\cdot A^{n-k}-S_{k-1}\cdot A^{n-(k-1)}- S_{n-1}\cdot A=0$%
.

}
\Thc Corollaire|
Si 
$\displaystyle\chi_A(X)\!=\!%
X^n+X^{n-1}c_1+\cdots+c_n\!=\!%
\sum_{i=0}^nc_iX^{n-i}%
\in\widetilde{\Lambda}[X]$
est le polyn{\^o}me caract{\'e}ristique 
g{\'e}n{\'e}ral,
alors la comatrice g{\'e}n{\'e}rale a l'expression~:
$$\widetilde{A}=(-1)^{n-1}%
\sum_{i=1}^{n}c_{i-1}A^{n-i}$$
\finc
\npage
{\small
\nsssecp {\bf Appendice\/}~: Th{\'e}or{\`e}me fondamental de l'arithm{\'e}tique dans
les anneaux de polyn{\^o}mes {\`a} coefficients entiers en un nombre fini de variables
$\Pi_n={\bf Z\/}[X_1,\ldots, X_n]$%
 .|
\Defns Pour tout entier naturel
$n\in{\bf N\/}$
l'anneau de polyn{\^o}me
$\Pi_n$
est
$\Pi_0={\bf Z\/}$
et, si
$n>0$
est  positif, d{\'e}fini par la relation de r{\'e}currence~:
$\Pi_n=\Pi_{n-1}[X_n]$%
. 

\Defn 
$P\in\Pi_n$
est  {\it positif,\/} not{\'e}
$P\!>\!0$%
, si soit
$n\!=\!0$
et
$P\!\in\!\Pi_0\!\setminus\!(-{\bf N\/})\!=\!{\bf Z\/}\!\setminus\!(-{\bf N\/})$
 est positif pour l'ordre usuel de
$\Pi_0\!=\!{\bf Z\/}$
et les r{\'e}cur\-rences~: si
$n\!>\!0, \displaystyle P\!=\!\sum_{i=0}^dP_iX_n^{d-i}; i\!=\!0,\ldots, d, P_i\!\in\!\Pi_{n-1}, P_0\!\ne\!0$
non nul de degr{\'e}
$d$
et {\it  (coefficient) dominant\/}
$P_0\!>\!0$
 positif\footnote{\small
d{\'e}j{\`a} d{\'e}fini par l'hypoth{\`e}se de r{\'e}currence.
}
 dans
$\Pi_{n-1}$%
.\hfill\break 
Il est dit {\it n{\'e}gatif,\/} not{\'e}
$P\!<\!0$%
, si son oppos{\'e}
$-P\!>\!0$
est positif. 
On note
$\Pi_n^+\!=\!\{P\in\Pi_n | P\!>\!0\}$
et
$\Pi_n^-(=-\Pi_n^+)\!=\!\{P\in\Pi_n | P\!<\!0\}$
les ensembles des {\'e}l{\'e}ments positifs et n{\'e}gatifs de
$\Pi_n$%
.
\Thc Proposition|
$({\romannumeral1})$
Si
$m\!<\!n$
alors
$\Pi_m\subset\Pi_n$
et si
$P\!\in\!\Pi_m^+$
on a~:
$P\!\in\!\Pi_n^+$%
.

$({\romannumeral2})$ 
$0\!\not\in\!\Pi_n^+\cup\Pi_n^-, \Pi_n^+\cap\Pi_n^-\!=\!\emptyset$
et on a la partition~: 
$\Pi_n\!=\!\Pi_n^+\amalg\{0\}\amalg\Pi_n^-$

$({\romannumeral3})$
L'ensemble
$\Pi_n^+$
est stable par multiplication
$\cdot$
et addition 
$+$%
.
\finc
{\small
\Dem 
$({\romannumeral1})$ 
Tout
$P\!\in\!\Pi_m\!\setminus\!\{0\}$%
,  consid{\'e}r{\'e} comme
$P\!\in\!\Pi_n$%
, est
$P\!=\!P_0$%
~:  de degr{\'e}
$0$
donc {\'e}gal {\`a} son coefficient dominant, d'o{\`u}
le r{\'e}sultat par r{\'e}currence sur
$n-m$%
.\hfill\carre

$({\romannumeral2})$ 
et
$({\romannumeral3})$ 
sont clairs si
$n\!=\!0$
(car
$\Pi_0\!=\!{\bf Z\/}$%
). On suppose
$n\!>\!0$
et ces r{\'e}sultats pour
$n-1$%
.\hfil\break 
Soit
$P\!=\!P_0X_n^e+\cdots+P_e, Q\!=\!Q_0X_n^f+\cdots+Q_f\!\in\!\Pi_n, P_0\!\ne\!0\!\ne\!Q_0$
deux polyn{\^o}mes non nuls,
les coefficients dominants de 
$P\cdot Q$
et
$P+Q$
sont respectivement
$C_0\!=\!P_0\cdot Q_0$
et
$S_0$
o{\`u}~:\hfill\break
 si
$e>f, S_0\!=\!P_0$%
, si
$e\!=\!f, S_0\!=\!P_0+Q_0$
et si
$e<f, S_0\!=\!Q_0$%
. Donc, si
$P_0, Q_0\!\in\!\Pi_{n-1}^+$%
, alors
$C_0, S_0\!\in\!\Pi_{n-1}^+$%
, soit si
$P, Q\!\in\!\Pi_n^+$
alors
$P\cdot Q, P+Q\!\in\!\Pi_n^+$%
~:
$({\romannumeral3})$
pour
$n$%
.
\hfill\carre\break
On a~:
$0\!\ne\!P_0\!\in\!\Pi_{n-1}\!\setminus\!\{0\}\!=\!\Pi_{n-1}^+\amalg\Pi_{n-1}^-$
ainsi
$P\!\in\!\Pi_n^+\amalg\Pi_n^-$
d'o{\`u}
$({\romannumeral2})$
pour
$n$%
.
\hfill\carre\findem

}

\thc Corollaire|L'anneau de polyn{\^o}mes
$\Pi_n$
{\`a} coefficients entiers, muni de 
$\geq$
d{\'e}finie par
$P\geq Q$
si 
$P-Q\in\Pi_n^+\cup\{0\}$
est un anneau totalement ordonn{\'e}.

En particulier l'anneau
$\Pi_n$
est int{\`e}gre donc, si
$n>0$
et
$P, Q\in\Pi_n\setminus\{0\}$%
, on a la relation~:
$${\mathop{\rm deg}\nolimits}_{X_n}(P\cdot Q)=%
{\mathop{\rm deg}\nolimits}_{X_n}(P)+{\mathop{\rm deg}\nolimits}_{X_n}(Q)$$
\finc
\vskip-3mm
{\small
\Dem Si
$P, Q, R\in\Pi_n$
on a~:
$P-P=0\in\Pi_n^+\cup\{0\}$
donc
$P\geq P$
et
$\geq$
est r{\'e}flexive.
Si
$P\geq Q\geq P$
alors
$P-Q\!=\!-(Q-P)\!\in\!(\Pi_n^+\cup\{0\})\cap(\Pi_n^-\cup\{0\})\!=\!\{0\}$%
, et
$\geq$
est antisym{\'e}trique.
Si
$P\!\geq\!Q\!\geq\!R$
on a
$P-R\!=\!(P-Q)+(Q-R)\in(\Pi_n^+\cup\{0\})+(\Pi_n^+\cup\{0\})\!=\!%
\Pi_n^+\cup\{0\}$
(par
$({\romannumeral3})$%
)
donc
$P\!\geq\!R$
et
$\geq$
est transitive et 
$\geq$
est un ordre sur
$\Pi_n^+$
qui, par
$({\romannumeral2})$%
, est total.\hfill\carre

Si
$P\geq Q$
et
$P+R-(Q+R)=P-Q\in\Pi_n^+\cup\{0\}$%
, donc
$P+R\geq Q+R$%
.\hfill\break
De plus si
$R\geq 0$%
, on a
$P\cdot R-(Q\cdot R)=(P-Q)\cdot R\in(\Pi_n^+\cup\{0\})\cdot(\Pi_n^+\cup\{0\})=\Pi_n^+\cup\{0\}$
d'o{\`u} par
$({\romannumeral3})$
$P\cdot R\geq Q\cdot R$
et l'ordre 
$\geq$
est compatible avec les op{\'e}rations de l'anneau
$\Pi$%
.\hfill\carre

Si
$\displaystyle P\!=\!\sum_{i=0}^{d={\mathop{\rm deg}\nolimits}P}\!\!a_iX_n^{d-i},%
 Q\!=\!\sum_{j=0}^{e={\mathop{\rm deg}\nolimits}Q}\!\!b_jX_n^{e-j}\in\Pi_n\!\setminus\!\{0\}$
alors
$a_0\!\cdot\!P,\,  b_0\!\cdot\!Q\!>\!0$
donc
$(a_0\cdot P)\dot(b_0\cdot Q)\!>\!0$
donc
$P\ne0, Q\ne0$
et
$\displaystyle P\cdot Q\!=\!\sum_{k=0}^{d+e}c_kX_n^{d+e-k}$
avec
$c_0\!=\!a_0\cdot b_0\ne0$%
. Ainsi
${\mathop{\rm deg}\nolimits}(P\cdot Q)\!=\!d+e$%
. \hfill\carre\findem
}

\rmc Remarque|
L'ensemble
$\Pi_n^+$
des polyn{\^o}mes positifs a
$1$
comme plus petit {\'e}l{\'e}ment, mais si
$n>0$%
, contrairement au cas
$n=0$
des entiers, la restriction de
$\geq$
{\`a} 
$\Pi_n^+$
n'est pas un bon ordre.
\finc
{\small
\Dem Soit
$P\in\Pi_n^+$
si 
${\mathop{\rm deg}\nolimits}P>0$
alors
${\mathop{\rm deg}\nolimits}(P-1)={\mathop{\rm deg}\nolimits}P>0$
et
$P\geq 1$%
.\hfill\break
Sinon, puisque
$1$
est le plus petit {\'e}l{\'e}ment de
$\Pi_0^+={\bf N\/}\setminus\{0\}$
on a aussi
$P\geq 1$%
.\hfill\carre

Comme, si
$n>0$%
, pour tout entier naturel
$m\in{\bf N}$
on a
$X-m\geq X-(m+1)$
la partie
$\{X-m | m\in{\bf N\/}\}\subset\Pi_1^+\subset\Pi_n^+$
est une partie non vide sans plus petit {\'e}l{\'e}ment de
$\Pi_n^+$%
.\hfill\findem
}
\thc Corollaire|
Les unit{\'e}s de
$\Pi_n$
sont
$\Pi_n^\ast=\{-1, 1\}$
\finc
{\small
\Dem Soit
$u\in\Pi_n^\ast$
une unit{\'e}, quitte {\`a} consid{\'e}rer
$-u$
on peut supposer 
$u\in\Pi_n^+$%
, donc
$u\geq 1$%
. Comme 
$0\geq v$
implique
$0=u\cdot 0\geq u\cdot v$%
. L'inverse
$v\in\Pi_n^\ast$%
, v{\'e}rifiant
$u\cdot v=1>0$%
, doit v{\'e}rifier
$v>0$%
, d'o{\`u}
$1=u\cdot v\geq u\cdot 1=u$
et par antisym{\'e}trie de
$\geq$%
, on a
$u=1$%
.\hfill\findem
}
{\small
\Defns Les applications {\it signe\/} et {\it valeur absolue\/}~:
$${\mathop{\rm sgn}\nolimits} : \Pi_n\rightarrow \{-1, 0, 1\},\quad
|\ | : \Pi_n\rightarrow \Pi_n^+\cup\{0\}$$
sont d{\'e}finies par
${\mathop{\rm sgn}\nolimits}(0)=0, |0|=0$%
, si
$P\in\Pi_n^+, {\mathop{\rm sgn}\nolimits}(P)=1, |P|=P$
et si
$P\in\Pi_n^-, {\mathop{\rm sgn}\nolimits}(P)=-1, |P|=-P$
et sont multiplicatives~:\hfill\break
 pour tout
$P, Q\in\Pi_n$
on a~:
${\mathop{\rm sgn}\nolimits}(P\cdot Q)={\mathop{\rm sgn}\nolimits}(P)\cdot{\mathop{\rm sgn}\nolimits}(Q), |P\cdot Q|=|P|\cdot|Q|$%
.
}
\Defn Un polyn{\^o}me
$P\in\Pi_n$
est dit {\it irr{\'e}ductible\/}
 si
$U, V\in\Pi_n$
v{\'e}rifient
$P\!=\!U\cdot V$
alors soit
$U\!\in\!\Pi_n^\ast, V\!=\!U^{-1}\cdot P$%
, soit
$V\!\in\!\Pi_n^\ast, U\!=\!V^{-1}\cdot P$%
.

\Defn Un polyn{\^o}me {\it premier\/}
est un polyn{\^o}me irr{\'e}ductible positif.\hfill\break
 On note
${\cal P\/}_n$
l'ensemble des polyn{\^o}mes premiers de 
$\Pi_n$%
.
\Thc Th{\'e}or{\`e}me fondamental de l'arithm{\'e}tique dans
$\Pi_n$%
|
Tout polyn{\^o}me non nul
a une unique {\it factorisation en premiers\/}~:
$$\hbox{\rm Si \/} M\in\Pi_n,\quad%
 M={\mathop{\rm sgn}\nolimits}(M)\prod_{p\in{\cal P}_n}p^{\nu_P(M)}$$
{\small
o{\`u} les  {\it multiplicit{\'es} 
$\nu_P(M)\in{\bf N\/}$
de 
$p$ 
dans 
$M$%
}
sont des entiers naturels presque tous nuls.
}
\finc
{\small
\Dem Si 
$n\!>\!0$
notons
$Y\!=\!X_n$
et
${\mathop{\rm deg}\nolimits}\!=\!{\mathop{\rm deg}\nolimits}_Y$%
. Si
$M\!=\!M_1\cdot M_2,  M_1, M_2\!\not\!\in\{-1, 1\}$
alors
$|M|\!>\!|M_1|, |M_2|$
et si
$n\!>\!0, {\mathop{\rm deg}\nolimits}(M)\!\geq\!%
{\mathop{\rm deg}\nolimits}(M_1), {\mathop{\rm deg}\nolimits}(M_2)$
avec {\'e}galit{\'e} pour un 
$M_i$
seulement si 
$M_j$
est de degr{\'e} 
$0$%
. L'ordre sur  valeurs absolues si
$n\!=\!0$
et  degr{\'e}s sinon {\'e}tant bon, 
 et par l'hypoth{\`e}se de r{\'e}currence pour le terme dominant si
$n>0$%
, on ne peut factoriser ind{\'e}finiment
$M_I\!=\!M_{i_1,\ldots,i_k}\!=\!M_{I,1}\cdot M_{I,2}$%
. A l'arr{\^e}t les facteurs
$M_I$
sont irr{\'e}ductibles
d'o{\`u} l'existence.\hfill\carre

L'unicit{\'e}~: On peut supposer le polyn\^ome
$M$
positif et raisonne par r\'ecurrence~:

Si
$n\!=\!0$
sur l'ordre usuel de
$M\in \Pi_{0}^{+}\setminus\{0\}\!=\!{\Bbb N\/}\setminus\{0\}$
et, si
$n\!>\!0$%
, sur l'ordre lexicographique des triplets d'entiers naturels
$(n, {\mathop{\rm deg}\nolimits}M, \sum_{p\!\in\!{\cal P}_n}\nu_p(M_0))$%
, ces  {\it tailles\/}
$t_{0}, t_{+}$
sont des bons ordres.

Soit
$p, q\!\in\!{\cal P}_n$
 premiers de deux d{\'e}compositions\footnote{\small
les produits
$u\!=\!p_2\cdots p_k, v\!=\!q_2\cdots q_l$
pouvant {\^e}tre vides (donc valoir 
$1$%
).
}~:
$p\cdot u\!=\!M\!=\!q\cdot v$
avec
$q\geq p$%
, alors soit

${\romannumeral 1})$
$p\!=\!q$
donc
$u\!=\!v$%
, de taille
$t_{\ast}(u)\!<\!t_{\ast}(M)$%
, d'o\`u le r{\'e}sultat par r{\'e}currence.\hfill\carre

soit
$q\!>\!p$%
. Si%
\footnote{\small
la preuve de Zermelo de factorialit\'e de
${\Bbb Z\/}$%
, initialis{\'e} par
$M\!=\!1$%
.
}
$n\!=\!0$
posons
$N\!=\!(q\!-\!p)\!\cdot\!v\!=\!q\!\cdot\!v\!-\!p\!\cdot\!v\!=\!p\!\cdot\!u\!-\!p\!\cdot\!v\!=\!p\!\cdot\!(u\!-\!v)$
et, si
$n\!>\!0, N\!=\!|p_{0}q-q_{0}pY^{e}|v\!=\!p\cdot |p_{0}u-q_{0}vY^{e}|$
o\`u
$e\!=\!{\mathop{\rm deg}\nolimits}q-{\mathop{\rm deg}\nolimits}p$%
. Dans les deux cas%
\footnote{\small
puisque
$(q-p)v\!<\!qv\!=\!M$
et
${\mathop{\rm deg}\nolimits}|p_{0}q-q_{0}pY^{e}|v\!<\!{\mathop{\rm deg}\nolimits}(qv)\!=\!{\mathop{\rm deg}\nolimits}M$%
.}
 $t{\ast}(N)\!<\!t_{\ast}(M)$
et, par r{\'e}currence,
$N$
se factorise  uniquement.

Si
$n\!=\!0$%
, ne divisant 
$q\!-\!p$%
, car il diviserait 
$q$
(
$\!\ne\!p$
et premier),
$p$
factorise
$v$%
. En ce cas~:

${\romannumeral 2})$
$u\!=\!q\cdot\frac{v}{p}\!<\!M$%
, d'o{\`u} le r{\'e}sultat car, par r{\'e}currence,
$u$
se factorise uniquement.\hfill\carre

Si
$n\!>\!0$
et
${\mathop{\rm deg}\nolimits}q = {\mathop{\rm deg}\nolimits}M$
alors
${\mathop{\rm deg}\nolimits}v=0$%
. Comme
$q$
est premier on ne peut avoir%
\footnote{\small
car on aurait
$q\!=\!pu$
avec, puisque
$p<q, p\ne q$
et, et
$p$
\'etant premier
$p\notin\Pi_{n}^{\ast}$%
.}
$v=1$%
.

Soit donc
$q'$
un premier de la seconde d\'ecomposition divisant
$v$%
, il divise
$p\cdot u\!=\!q\cdot v$%
.

Le lemme de Gau\stz%
\footnote{\small
{\sl Un anneau
$A$
est {\it factoriel\/} si le th\'eor\`eme fondamental de l'arithm\'etique a lieu dans
$A$%
.\hfill\break
Dans un anneau factoriel
$A$%
, le lemme de Gau\stz \
est~: {\sl Si
${\goth p\/}\in A$ est un premier d'un anneau factoriel ne divisant pas deux polyn{\^o}mes
$S, T\in A[Y]$
alors il ne divise pas leur produit
$ST$%
.}
}
\Dem 
Soit
$i\!=\!\min\{i ; {\goth p\/}\not| S_{i}\}, j\!=\!\min\{j ; {\goth p\/}\not| T_{j}\}$%
, le c\oe ficient d'indice
$i+j$
de
$ST$%
~:
$$\displaystyle (ST)_{i+j}\!=\!\sum_{k\!=\!0}^{i-1}S_{k}T_{i-k+j}+S_{i}T_{j}+\sum_{l\!=\!j-1}^{0}S_{i+j-l}T_{l}$$
Par d\'efinition de
$i, j, {\goth p\/}$
divise les
$S_{k}, T_{l}$%
, donc les sommes extr{\^e}mes et, comme
$A$
est factoriel, 
${\goth p\/}$ ne divise pas le terme central
$S_{i}T_{j}$%
, donc ni le coefficient
$(ST)_{i+j}$%
, ni le produit
$ST$%
.\findem
}
appliqu\'e \`a
${\goth p\/}\!=\!q'$
et l'anneau
$A\!=\!\Pi_{n-1}$%
,  factoriel par r\'ecurrence, 
donne~:

Soit
$q'$
divise
$u$
et l'unicit\'e de la factorisation de
$p\cdot\frac{u}{q'}\!=\!\frac{M}{q'}\!=\!q\cdot\frac{v}{q'}$
conclut.\hfill\carre

Soit
$q'$
divise
$p$
donc, puisque
$p$
est premier,
$q'\!=\!p$
et
${\romannumeral 2\/})$
conclut.\hfill\carre

Reste le cas
$n\!>\!0, {\mathop{\rm deg}\nolimits}q\!<\!{\mathop{\rm deg}\nolimits}M$%
. Comme
$N$
a une unique factorisation en premiers,
$p$
divise

Soit
$v$%
, auquel cas
${\romannumeral 2\/})$
encore conclut.\hfill\carre

Soit
$|p_{0}q-q_{0}pY^{e}|$%
, donc
$p_{0}q$%
.\hfill\break
Comme
${\mathop{\rm deg}\nolimits}(p_{0}q)\!=\!{\mathop{\rm deg}\nolimits}q\!<\!{\mathop{\rm deg}\nolimits}M$
on a
$t_{\ast}(p_{0}q)\!<\!t_{\ast}(M)$
et, par r\'ecurrence,
$p$
divise~:

Soit
$q$%
, donc
$p\!=\!q$%
, cas trait\'e en
${\romannumeral 1\/})$%
.\hfill\carre

Soit
$p_{0}$%
, donc
$p\!=\!p_{0}$
et, par le lemme de Gau\stz,
$p$
divise soit
$q$
donc encore
$p\!=\!q$
, soit
$v$%
.\hfill\break
Ces deux cas \'etant  trait\'es en
${\romannumeral 1\/})$
et
${\romannumeral 2\/})$%
, le th\'eor\`eme  est \'etabli.\hfill\carre\findem

}

}
\npage 
\nssecp D{\'e}terminants sur un anneau commutatif 
et applications.|
\parc
On rappelle qu'un anneau de polyn{\^o}mes
{\`a} coefficients entiers
${\Bbb Z\/}[X_i]_{i\in I}$
est un {\it anneau commutatif libre\/}
{\it sur}
l'ensem\-ble
$I$
(ou {\it sur\/} ses variables
$X_i, {i\!\in\!I}$%
), {\it c. a d.\/} pour toute
application
$f : I\rightarrow A$
de l'ensemble 
$I$
vers un anneau commutatif
$A$
il y a un unique morphisme d'anneau
$F : {\Bbb Z}[X_i]_{i\in I}\rightarrow A$
tel que pour tout
$i\in I$
on ait
$F(X_i)=f(x_i)$%
.
\finc
On suppose dans cette sous-section  que l'anneau 
$\Lambda\/$ 
est commutatif.

{\small
Pour
toute 
$M=
(m_{i,\,j})_{1\leq i,\,j\leq n}\in%
{\rm M\/}_n\,(\Lambda)\/$
il y a un unique morphisme d'anneau,
$\rho_M : \widetilde{\Lambda}_n\rightarrow \Lambda\/$%
, dit 
{\it sp{\'e}cialisation en \/} 
$M\/$ 
et d{\'e}fini par
$\rho_M\,(a_{i,\,j})=m_{i,\,j}\/$ 
pour 
$i,\,j=1,\,\ldots,\,n\/$%
.
}
\Defn Le 
{\it d{\'e}terminant\/} 
de la matrice 
$M\in{\rm M\/}_n\,(\Lambda)\/$
est l'image 
$\rho_M\,({\rm det\/}\,(A))\/$ 
du d{\'e}terminant g{\'e}n{\'e}ral
par la sp{\'e}cialisation en 
$M\/$%
. On le note 
${\rm det\/}\,(M)\/$%
.

{\small
Les neuf relations 
$(0)\/$ 
{\`a} 
$(8)\/$ 
de {\bf 2.2.2\/} et le th{\'e}or{\`e}me de {\bf 2.2.3\/}
sont alors valables si les matrices 
$A, B, U, V, W\/$ 
qui y interviennent sont {\`a} coefficients
dans l'anneau commutatif 
$\Lambda\/$%
.
}
\Thc Th{\'e}or{\`e}me| Soit 
$X\/$ 
et 
$Y\/$ 
des modules libres sur l'anneau commutatif 
$\Lambda\/$ 
ayant des bases finies 
${\cal B\/}=(b_1,\,\ldots,\,b_n)\/$ 
et
${\cal C\/}=(c_1,\,\ldots,\,c_n)\/$ 
avec le m{\^e}me nombre 
$n\/$
d'{\'e}l{\'e}ments et
$f : X\rightarrow Y\/$ 
un morphisme de matrice 
$M\/$ 
dans les bases 
${\cal B\/}$ et ${\cal C\/}$%
.

Alors on a les deux \og cercles\fg d'{\'e}quivalence~:

({\romannumeral1}) D'une part entre~:

\quad ({\romannumeral1} 1) L'application 
$f\/$ 
est injective. 

\quad ({\romannumeral1} 2) Le d{\'e}terminant 
${\rm det\/}\,(M)\/$ 
de
$M$
n'est pas  diviseur de z{\'e}ro dans 
$\Lambda\/$.

(\romannumeral2) D'autre part entre~:

\quad ({\romannumeral2} 1) L'application 
$f\/$ 
est surjective.

\quad ({\romannumeral2} 2) L'application 
$f\/$ 
est un isomorphisme.

\quad ({\romannumeral2} 3) Le d{\'e}terminant 
${\rm det\/}\,(M)\!\in\!\Lambda^{\bullet}\/$
de
$M$
 est une unit{\'e} de l'anneau 
$\Lambda\/$%
. 

{\small
En ce cas l'isomorphisme inverse 
$f^{-1} : Y\rightarrow X\/$
est, dans les bases 
${\cal C\/}$ et ${\cal B\/}$%
, de matrice~:
}
$$M^{-1}=({\rm det\/}(M))^{-1}\,\widetilde{M}\/$$

\finc
{\small
\Dem (\romannumeral1) Si
$f$ 
n'est pas injective il y a
$0\ne X=(x_1,\ldots, x_n)\in\Lambda^n$
tel que
$M\cdot X=0$
d'o{\`u}
${\mathop{\rm det\/}\nolimits}(M) X\!=\!%
(\tilde{M}\cdot M)\cdot X\!=\!%
\tilde{M}\cdot (M\cdot X)\!=\!0$%
, donc pour
$i\!=\!1,\ldots, n$
on a
${\mathop{\rm det\/}\nolimits}(M) x_i\!=\!0$%
. Ainsi 
${\mathop{\rm det\/}\nolimits}(M)$
est diviseur de z{\'e}ro, car
$X$
{\'e}tant non nul, une de ses composantes
$x_i\!\ne\!0$
est non nulle.

Si 
${\mathop{\rm det\/}\nolimits}(M)$
n'est pas diviseur de z{\'e}ro alors pour tout
$0\ne X\in\Lambda^n$
on a~:\hfill\break
$0\!\ne\!{\mathop{\rm det\/}\nolimits}(M) X\!=\!%
\tilde{M}\cdot (M\cdot X)$%
, donc
$M\cdot X\!\ne\!0$
et, d'apr{\`e}s la Proposition 1 de {\bf 1.3\/},
$f$
est injective.\hfill\carre

({\romannumeral2}) L'implication
$({\romannumeral2}\ 2)\Rightarrow({\romannumeral2}\ 1)$
est claire. Si
$f$ est surjective,
d'apr{\`e}s le Corrolaire de {\bf 1.4\/},\hfill\break
il y a 
$g : Y\rightarrow X$
tel que
$g\circ f={\mathop{\rm Id}\nolimits}_Y$%
. Si
$N$
est sa matrice dans les bases
${\cal C\/}$
et
${\cal D\/}$
on a
$N\cdot M={\mathop{\rm Id\/}\nolimits}$
donc
${\mathop{\rm det\/}\nolimits}(N)
{\mathop{\rm det\/}\nolimits}(M)=1$
et
${\mathop{\rm det\/}\nolimits}(M)$
est une unit{\'e} de
$\Lambda$%
, do{\`u} l'implication
$({\romannumeral2}\ 1)\Rightarrow({\romannumeral2}\ 3)$%
.\hfill\break
 L'implication
$({\romannumeral2}\ 3)\Rightarrow({\romannumeral2}\ 2)$
suit de ce qu'alors
$M\cdot
({\mathop{\rm det\/}\nolimits}(M)^{-1}\tilde{M})=
{\mathop{\rm Id\/}\nolimits}=
({\mathop{\rm det\/}\nolimits}(M)^{-1}\tilde{M})%
\cdot M$%
~: \hfill\break
l'application
$g : Y\rightarrow X$
de matrice
$({\mathop{\rm det\/}\nolimits}(M)^{-1}\tilde{M})$
est isomorphisme r{\'e}ciproque de
$f$%
.\hfill\carre\findem
}
\Thc Corollaire|
Soit
$m, n$
deux entiers naturels et
$f :\Lambda^n\rightarrow\Lambda^m$
un morphisme surjectif.
Alors, si l'anneau
$\Lambda\ne\{0\}$
n'est pas l'anneau nul, on a~:
$m\leq n$%
.
\finc
{\small
\Dem Si 
$m>n$
le morphisme
$\overline{f} : \Lambda^m\rightarrow\Lambda^m$
qui
 envoit le
$k^{\hbox{i\`eme}}$
 {\'e}l'{\'e}ment
$e_k$
de la base canonique,
 si
$n<k\leq m$
 sur z{\'e}ro
(et sur
$f(e_k)$
sinon)
est toujours surjectif, de matrice 
de d{\'e}terminant nul\footnote{\small
car ayant au moins une une colonne nulle.
}, contredisant\footnote{\small
car pour tout
$a\in\Lambda, 0a=0\ne1$
car l'anneau 
$\Lambda$
est non nul.
 } 
les {\'e}quivalences ({\romannumeral2})
du th{\'e}or{\`e}me.\hfill\findem
}
\Thc co-Corrollaire| Si l'anneau
$\Lambda\ne\{0\}$
n'est pas l'anneau nul et
$m, n$
sont deux entiers naturels tels que
$\Lambda^m$
et
$\Lambda^n$
sont isomorphes alors
$m= n$%
.\hfill\findem
\finc
\vskip-3mm
\rmc Remarques |({\romannumeral1}) Si l'anneau 
$\Lambda$
n'est pas commutatif, le co-Corrollaire
(donc le Corollaire aussi) n'a plus lieu~:
Soit
$M$
un groupe ab{\'e}lien non nul muni d'un isomorphisme
$\Phi : M\oplus M\rightarrow M$%
, d{\'e}fini par
$\phi_i\!=\!\Phi\circ{\mathop{\rm i\/}\nolimits}_i :%
 M\rightarrow M, i\!=\!1, 2$
(par exemple
$M\!=\!{\bf Z\/}^{({\bf N\/})}, 
\phi_1(e_{n})\!=\!e_{2n}, \phi_2(e_{n})\!=\!e_{2n+1}$%
).

Alors si
$\Lambda\!=\!%
{\mathop{\rm End\/}\nolimits}(M)$
on a~: 
$\Lambda_d=\phi_1\Lambda\oplus\phi_2\Lambda$%
. D'o{\`u} pour tout entier
$n$
positif~:\hfill\break
$\displaystyle\Lambda_d=
\phi_1^n\Lambda\oplus_{i=0}^{n-1}%
\phi_1^i\circ \phi_2\Lambda$ 
et le module monog{\'e}ne
$\Lambda_d$
contient le sous-module libre {\`a} un nombre
infini de g{\'e}n{\'e}rateurs
$\displaystyle 
\oplus_{n\in{\bf N\/}}\phi_1^n\circ \phi_2\Lambda$%
. Par contre~:

({\romannumeral2})%
\footnote{\small voir la section {\bf 3\/} suivante.
}%
Si 
$\Lambda$
est un {\it corps\/}%
\footnote{\small
{\it c. a d.\/} anneau, non n{'e}c{\'e}ssairement
commutatif, non nul (%
$1\ne 0$%
) et tout {\'e}l{\'e}ment
$0\ne \lambda\in\Lambda$
non nul est {\it inversible\/} 
[{\it c. a d.\/} il y a
$\mu\in\Lambda$
(not{\'e}
$\mu=\lambda^{-1}$%
)
tel que
$\lambda \mu=1=\mu\lambda$%
].
}
ce Corrollaire (donc le co-Corollaire) a lieu.

({\romannumeral3}) Sans supposer l'anneau
$\Lambda$
commutatif le Corrollaire
(et donc le co-Corrollaire) restent vrais pour les
modules libres sur des ensembles infinis.
Pus pr{\'e}cis{\'e}ment~:
\finc
\vskip-3mm
\Thc Proposition|Si l'anneau
$\Lambda\ne\{0\}$
est non nul et
$f :\Lambda^{(J)}\rightarrow\Lambda^{(I)}$
est un morphisme surjectif alors,
si
$I$
est infini, son cardinal 
est au plus celui de
$J$%
~:
$${\mathop{\rm Card\/}\nolimits}(I)\leq%
{\mathop{\rm Card\/}\nolimits}(J)$$
\finc
\vskip-5mm
{\small
\Dem Soit 
$(a_{i,\, j})_{(i,\,j)\in I\times J}$
la matrice de
${\mathop{\rm i\/}\nolimits}_I\circ%
f :\Lambda^{(J)}\rightarrow\Lambda^{I}$
alors pour tout
$j\in J$
l'ensemble
$I_{f,\, j}=\{i\in I | a_{i,\,j}\ne0\}$
est fini, la r{\'e}union
$\displaystyle I_f=\cup_{j\in J}I_{f,\, j}$
est donc finie ou de cardinal
${\mathop{\rm Card\/}\nolimits}(I_f)\leq
{\mathop{\rm Card\/}\nolimits}(J)$
au plus celui de
$J$%
. L'image de 
$f$
{\'e}tant incluse dans 
$\Lambda^{(I_f)}\subset\Lambda^{(I)}$
qui serait  sous-module strict si
${\mathop{\rm Card\/}\nolimits}(I)>%
{\mathop{\rm Card\/}\nolimits}(J)$%
, le morphisme 
$f$
ne serait pas surjectif,
en ce cas.\hfill\findem
}
\Thc Compl{\'e}ment|
Soit
$M$
un module de type fini sur un anneau commutatif
$\Lambda$%
. Alors tout endomorphisme 
$f : M\rightarrow M$
surjectif de
$M$
est un isomorphisme.\hfill\break
{\small
Plus pr{\'e}cis{\'e}ment si
$M$
est engendr{\'e} par
$m_1,\ldots, m_n\in M$
et 
$A=(a_{i,\, j})_{1\leq i, j\leq n}\!\in\!%
{\mathop{\rm M\/}\nolimits}_n(\Lambda)$
une matrice telle que pour
$i=1,\ldots, n$
on ait 
$\displaystyle m_i\!=\!f(\sum_{j=1}^na_{i,\, j}m_j)%
\!=\!\sum_{j=1}^na_{i,\, j}f(m_j))$%
.\hfill\break
 Alors 
$f$
est d'inverse
$f^{-1}=P(f)$
o{\`u}\footnote{\small
l'intervention de
$\chi_A$
ici est uniquement pour d{\'e}crire la formule
avec un invariant connu.\hfill\break
 La preuve n'utilise que
$XP(X)\!=\!1-{\mathop{\rm det\/}\nolimits}%
({\mathop{\rm Id\/}\nolimits}-XA)$%
, et  pas l'identit{\'e} de Cayley-Hamilton.
}
le polyn{\^o}me
$P(X)\!\in\!\Lambda[X]$
est
$P=X^{-1}-X^{n-1}\chi_A(X^{-1})$%
.
}
\finc
{\small
\Dem
Dans le sous-anneau (commutatif!)
$\Lambda[f]\!\subset\!%
{\mathop{\rm End\/}\nolimits}_\Lambda(M)$
engendr{\'e} par
$f$
il suit des rela\-tions
$({\mathop{\rm det\/}\nolimits}%
({\mathop{\rm Id\/}\nolimits}-XA)){\mathop{\rm Id\/}\nolimits}\!=\!%
\widetilde{%
({\mathop{\rm Id\/}\nolimits}-XA)}
({\mathop{\rm Id\/}\nolimits}-XA),
1-XP(X)\!=\!1-P(X)X\!=\!%
{\mathop{\rm det\/}\nolimits}%
({\mathop{\rm Id\/}\nolimits}-AX)$%
~:
$$({\mathop{\rm Id\/}\nolimits}-f\circ P(f)){\mathop{\rm Id\/}\nolimits}\!=\!%
({\mathop{\rm Id\/}\nolimits}-P(f)\circ f){\mathop{\rm Id\/}\nolimits}\!=\!%
\widetilde{%
({\mathop{\rm Id\/}\nolimits}-Af)}%
({\mathop{\rm Id\/}\nolimits}-Af)%
\!=\!0$$
puisque pour
$i\!=\!1,\ldots, n$
on a
$\displaystyle
({\mathop{\rm Id\/}\nolimits}-Af)\cdot^t\!\!(m_1,\ldots,m_j,\ldots,m_n)_i%
=m_i-\sum_{i=1}^na_{i,\, j}f(m_j)\!=\!0$%
.\hfill\findem

}
\npage
\nsecp Espaces vectoriels et dimension.|
\vskip5mm
Dans cette section
$K$
d{\'e}signe un corps (non suppos{\'e} commutatif).
\Defns Un {\it espace vectoriel\/}
(abr{\'e}g{\'e} en {\it e.v.\/}) {\it {\`a} gauche\/} 
(resp. {\it {\`a} droite\/})
$V$
{\it sur\/} le corps
$K$
 [dit aussi {\it
$K$%
-espace vectoriel\/}
(abr{\'e}g{\'e} en {\it 
$K$%
-e.v.\/})
{\it {\`a} gauche\/}
(resp. {\it {\`a} droite\/})]
est un
$K$%
-module {\`a} gauche (resp. {\`a} droite).\hfill\break
Un {\'e}l{\'e}ment
$v\in V$
d'un 
$K$%
-espace vectoriel est appel{\'e}
{\it vecteur de\/}
$V$%
.

{\small
Suivant l'usage, sauf pour~:

les syst{\`e}mes lin{\'e}aires en {\bf 3.4.2\/},

les  structures rationelles sur un
sous-corps en {\bf 3.5\/} et leurs applications~:

quelques
propri{\'e}t{\'e}s des corps gauches en {\bf 3.5.2\/}.

la dimension du dual d'un espace vectoriel de dimension
infinie en {\bf 3.5.3\/},\hfill\break
on ne consid{\`e}rera que des structures {\`a} gauche
et abr{\'e}gera l'expression\hfill\break
{\it espace vectoriel {\`a} gauche\/}
 en 
{\it espace vectoriel\/}.

De m{\^e}me avant {\bf 3.5\/},
ne consid{\'e}rant que le seul corps
$K$%
,
 on ommetra de pr{\'e}ciser
{\it sur
$K$%
.
}

}

%
%

%
\nssecp Combinaisons lin{\'e}aires.|
{\small
On rappelle (Cf. {\bf 1.4\/})
que
$K^{(I)}$
d{\'e}signe l'espace vectoriel des familles
$(\lambda_i)_{\in I}\!\in\!K^{(I)}$
d'{\'e}l{\'e}\-ments
$\lambda_i\!\in\!K$
du corps
$K$
qui sont presque tous nuls
{\small[%
{\it c. a d.\/} de {\it supports\/}
$\{i\in I | \lambda_i\ne 0\}$
 finis].
}%
}

\Defns Soit
$\underline{v}=(v_i)_{\in I}$
une famille indic{\'e}e par un ensemble
$I$
d'{\'e}l{\'e}ments
$v_i\in V$
d'un espace vectoriel
$V$%
. On consid{\`e}re l'application lin{\'e}aire~:
$${\mathop{cl_{\underline{v}}}\nolimits} :%
K^{(I)}\rightarrow V,\
\underline{\lambda}=(\lambda_i)_{i\in I}\mapsto%
{\mathop{cl_{\underline{v}}}\nolimits}%
(\underline{\lambda})=%
\sum_{i\in I}\lambda_i v_i$$
qui {\`a}
$\underline{\lambda}=(\lambda_i)_{\in I}\!\in\!K^{(I)}$
associe la
{\it combinaison lin{\'e}aire de\/}
$\underline{v}$
{\it {\`a} coefficients\/}
$\underline{\lambda}$%
.\hfill\break
La famille
$\underline{v}\!=\!(v_i)_{i\in I}$
est 
{\it g{\'e}n{\'e}ratrice 
{\small
(sur
$K$%
, dans
$V$%
)\/}
}
si le morphisme
${\mathop{cl_{\underline{v}}}\nolimits}$
est surjectif~:\hfill\break
tout vecteur
$v\in V$
est
$\displaystyle v=\sum_{i\in I}\lambda_i v_i$%
, combinaison lin{\'e}aire de cette famille.

L'ensemble
des combinaison lin{\'e}aires de la famille
$\underline{v}$%
, image 
${\mathop{\rm Im\/}\nolimits}(c_{\underline{v}})$
du morphisme 
${\mathop{cl_{\underline{v}}}\nolimits}$%
, est le {\it sous-espace vectoriel engendr{\'e} par\/}
la famille
$\underline{v}$
et est not{\'e}
$<\underline{v}>=<v_i, i\in I>$
ou, si
$I=\{1,\ldots,n\}$
est fini,
$<v_1,\ldots, v_n>$%
.

\Defn Une {\it relation lin{\'e}aire\/} de la famille
$\underline{v}\!=\!(v_i)_{i\in I}, v_i\!\in\!V$
est une combinaison lin{\'e}aire nulle~:
$\displaystyle\sum_{i\in I}\lambda_i v_i\!=\!0\!\in\!V$%
. 
Elle est dite
{\it non triviale\/}
si la famille de ses coefficients 
$(\lambda_i)_{\in I}\!\ne\!0\!\in\!K^{(I)}$
est non nulle {\small
({\it c. a d.\/} il y a un
$i\in I$
avec
$\lambda_i\ne0$%
)%
}%
.

\Defn La famille de vecteurs
$\underline{v}\!=\!(v_i)_{i\in I}, v_i\!\in\!V$
est
{\it libre\/}
[ou {\it lin{\'e}airement ind{\'e}pendante\/}]
{\small
({\it sur
$K$%
, dans
$V$%
)}
}
si le morphisme
${\mathop{cl_{\underline{v}}}\nolimits}$
est injectif, ce qui d'apr{\`e}s  la 
 Proposition 1
de {\bf 1.3\/} {\'e}quivaut {\`a}
{\sl sa seule relation lin{\'e}aire est la
relation triviale~: 
$\displaystyle\sum_{i\in I}\lambda_i v_i\!=\!0$
si et seul pour tout
$i\!\in\!I$
on a 
$\lambda_i\!=\!0$%
.}
Une famille non libre est dite {\it li{\'e}e.\/}
\rmc Remarques|
1) {\^E}tre g{\'e}n{\'e}ratrice pour une famille
$\underline{v}=(v_i)_{i\in I}, v_i\in V$
de vecteurs d'un espace vectoriel
$V$
 ne d{\'e}pend que de sa partie image
$\{v_i | i\in I\}\subset V$%
.

2) Une famille 
$\underline{v}$
non injective n'est\footnote{\small
 puisque si pour des
$i\ne j\in I$
on a
$v_i=v_j$
alors
$1 v_i+(-1) v_j=0$%
.
}
pas libre.
\finc
\vskip-3mm
\Defns Une partie
$X\subset V$
d'un e.v.
$V$
est {\it g{\'e}n{\'e}ratrice\/}, resp. {\it libre\/},
si la famille
$\underline{X}=(x)_{x\in X}$
de ses {\'e}l{\'e}ments indic{\'e}s par eux-m{\^e}me est
g{\'e}n{\'e}ratrice, resp. libre.
\Defns Une famille
$\underline{b}=(b_i)_{i\in I}, b_i\in V$%
, resp. partie
$B\subset V$
d'un espace vectoriel
$V$
est une {\it base de\/} 
$V$
si elle est libre et g{\'e}n{\'e}ratirice
{\it c. a d.\/} l'application combinaison
${\mathop{cl_{\underline{b}}}\nolimits} :%
 K^{(I)}\rightarrow V,
\underline{\lambda}\mapsto%
{\mathop{cl_{\underline{b}}}\nolimits}%
(\underline{\lambda})=%
\sum_{i\in I}\lambda_i b_i$
est un isomorphisme.

 L'isomorphisme inverse
${\mathop{co_{\underline{b}}}\nolimits}=%
{\mathop{cl_{\underline{b}}}\nolimits}^{-1} %
: V\rightarrow K^{(I)}$
est dite {\it application coordonn{\'e}es
de\/} l'espace vectoriel 
$V$
{\it dans la base
$\underline{b}$%
.\/}
\rmc Exemples|
$({\romannumeral1})$
La famille
$(e_i)_{i\in I}\!\in\!K^{(I)}\!=\!L_K$
canonique de base\footnote{\small
(C.f. {\bf 1.4\/})
$(e_i)_{i\in I}\in K^{(I)}=L_I,
e_i=(e_{i_j})_{j\in I}; e_{i_i}=1$
et si
$j\ne i, e_{i_j}=0$
} 
de  l'espace vectoriel
libre sur le corps
$K$ 
est g{\'e}n{\'e}ratrice et libre, la
{\it base canonique\/} de
$K^{(I)}$%
.

$({\romannumeral2})$ 
Si  tout {\'e}l{\'e}ment
$x_i\in X_i\subset V, i\in I$
d'une partie
$X_i$ 
d'un espace vectoriel
$V$
est combinaison lin{\'e}aire
$x_i=\sum_{j\in J}\lambda_{i,\, j} y_j$%
d'une partie
$Y=\{y_j | j\in J\}$
alors toute combinaison lin{\'e}aire
$c=\sum_{i\in I}\mu_i x_i$
des
$x_i, i\!\in\!I$
est
$c\!=\!\sum_{i\in I}\mu_i\sum_{j\in J}\lambda_{i,\, j}y_j\!=\!%
\sum_{j\in J}\bigl(%
\sum_{i\in I}\mu_i\lambda_{i,\, j}\bigr)y_j$%
,\hfill\break
combinaison lin{\'e}aire des
$y_j, j\!\in\!J$%
. En particulier~:

{\sl  Une partie engendrant une partie
g{\'e}n{\'e}ratrice est g{\'e}n{\'e}ratrice.
}

$({\romannumeral3})$ 
Une famille 
$\underline{v}=(v_i)_{i\in\emptyset}$
(ou partie
$X=\emptyset\subset V$%
) vide est libre.

$({\romannumeral4})$
Une famille
$\underline{v}=(v_i)_{i\in\{i_0\}}$
(ou partie
$X=\{v_{i_0}\}$%
)
{\`a} un seul {\'e}l{\'e}ment est libre si
et seulement si cet {\'e}l{\'e}ment
$v_{i_0}\ne0$
est non nul.

$({\romannumeral5})$
Une famille
$\underline{v}\!=\!(v_i)_{i\in\{i_0, i_1\}}$
(ou partie
$X\!=\!\{v_{i_0}, v_{i_1}\}$%
)
{\`a} deux {\'e}l{\'e}ments est li{\'e}e
si et seulement si soit
$v_{i_0}\!=\!0$
est nul, soit il y a
$\lambda\!\in\!K$
tel que
$v_{i_1}\!=\!\lambda v_{i_0}$%
.
\finc
\Thc Proposition| Soit
$X\!\subset\!V$
une partie g{\'e}n{\'e}ratrice d'un e.v.
$V$
alors une partie
$Y\subset X$
maximale, pour
$\subset$
 dans l'ensemble des parties libres de
$X$
est g{\'e}n{\'e}ratrice.
\finc
{\small
\Dem Soit
$x\in X$
alors soit
$x=y_x\in Y$%
,
soit la partie
$Z=\{x\}\cup Y$
est li{\'e}e.\hfill\break
Ainsi il y a\footnote{\small
 dans le premier cas~: 
$\lambda_x=1,\lambda_{y_x}=-1$
et si
$y\ne y_s, \lambda_y=0$%
, par d{\'e}finition dans le second.
}
$\lambda_x\in K$
et, pour
$y\in Y$
des
$\lambda_y\in K$
non tous nuls tels que
$0=\lambda_x x+\sum_{y\in Y}\lambda_y y\in V$%
. \hfill\break
Comme la partie
$Y$
est libre,
$\lambda_x\ne0$%
, et 
$\displaystyle x=\sum_{y\in Y}(-{\lambda_x}^{-1}\lambda_y)y$
est combinaison lin{\'e}aire des
$y\in Y$%
. L'espace
$V$
{\'e}tant engendr{\'e} par les
$x\in X$
l'est donc par les
$y\in Y$%
. Ainsi
$Y$
est g{\'e}n{\'e}ratrice.\hfill\findem 
}

\Thc Corollaire|
Tout espace vectoriel 
$V$
poss{\`e}de une base.
Plus pr{\'e}cis{\'e}ment~:
 
Soit
$Y, X\!\subset\!V$
deux parties libres avec
$X$
de plus g{\'e}n{\'e}ratrice de l'espace
$V$%
,\hfill\break alors 
$V$ 
a une base
$B$
v{\'e}rifiant
$Y\subset B\subset Y\cup X$%
~:

{\small
\og De toute partie partie g{\'e}n{\'e}ratrice,
on compl{\`e}te en une base
toute partie libre\fgf.
}
\finc
{\small
\Dem L'inclusion {\'e}tant inductive sur
les parties libres
$Y',\ Y\!\subset\!Y'\!\subset\!Y\!\cup\!X$
il y en a, par le lemme de Zorn, une
$Y'=B$
qui est maximale, par la Proposition elle
est g{\'e}n{\'e}ratrice.
\hfill\findem

}
\Thc co-Corollaire|
Tout sous-espace
$U\subset V$
d'un espace vectoriel
$V$
poss{\`e}de un suppl{\'e}mentaire
$W\subset V=U\oplus W$%
.
\finc
{\small
\Dem Soit
$(b_k)_{k\in I\amalg J}, b_k\in V$
une base compl{\'e}tant une base
$(b_i)_{i\in I}, b_i\in U$
du sous-espace
$U$
de 
$V$%
. Il suffit de prendre
$W=\oplus_{j\in J}K b_j$%
.\hfill\findem
}

\rmc Remarque|Les d{\'e}finition de partie
g{\'e}n{\'e}ratrice,
libre, base ont un sens pour les modules sur
un anneau quelconque,
mais les propri{\'e}t{\'e}s des exemples 
$({\romannumeral4}), ({\romannumeral5})$
et la proposition n'ont plus lieu~:

a) Dans le
${\bf Z\/}$%
-module
${\bf Z\/}/6{\bf Z\/}$
la partie
$\{\overline{2}\}$
est, car
$3\overline{2}=\overline{6}=0$%
, non libre {\`a} un {\'e}l{\'e}ment.

b) Dans le
${\bf Z\/}$%
-module
${\bf Z\/}$
la partie
$\{2\}$
est libre maximale, mais non g{\'e}n{\'e}ratrice.
\finc
\nssecp Dimension d'un espace vectoriel.|
\Thc Lemme|
Soit
$f : K^{(J)}\rightarrow K^{(I)}$
un morphisme surjectif entre des modules libres
sur des ensembles
$I, J$%
, alors~:
$${\mathop{\rm Card\/}\nolimits}(I)\leq%
{\mathop{\rm Card\/}\nolimits}(J)$$
\finc
{\small
\Dem Dans le cas o{\`u} l'ensemble
$I$
est infini, cela est m{\^e}me vrai pour
les modules libres
sur un anneau
$\Lambda$
 quelconque et a {\'e}t{\'e} {\'e}tabli dans
la Proposition en fin de {\bf 2.3\/}.

Le r{\'e}sulat est aussi clair si
$I=\emptyset$%
, on suppose donc
$I=\{1,\ldots,m\}$
et\footnote{\small
Comme
$K^{(\emptyset)}=\{0\}$%
, si
$I \ne\emptyset,
K^{(\emptyset)}\rightarrow K^{(I)}\ne\{0\}$
ne peut {\^e}tre surjectif.
}
$ J=\{1,\ldots, n\}$
finis d'ordre
$n, m>0$
positifs et proc{\`e}de par r{\'e}currence 
(initialis{\'e}e par 
$I\!=\!\emptyset, m\!=\!0$
ci-dessus) sur
$m$%
~:

 Soit pour
$1\leq j\leq n, v_j=f(e_j)=%
(a_{1,\, j},\ldots, a_{m,\, j})$%
, les vecteurs images de la base
$e_1,\ldots, e_n$
canonique de
$K^{(J)}$%
. L'une des derni{\`e}res
coordonn{\'e}es
$a_{m,\, j}\!\ne\!0$
est (car
$f$
est  surjectif) non nulle. Quitte {\`a} renum{\'e}roter
les
$e_j$%
,
on peut supposer
$a_{m,\, n}\ne0$%
. 

Soit pour
$1\leq j<n, w_j=v_j-(a_{m,\, j}a_{m,\, n}^{-1})v_n=%
(b_{1,\,j},\ldots, b_{m-1,\, j}, 0);\
w_n=v_n$
et
$$g : K^n=K^{(J)}\rightarrow K^{(I)}=K^m,\ g(e_j)=w_j$$

 Pour tout
$(x_1,\ldots, x_n)\!\in\!K^n$%
, si pour
$1\!\leq\!j<\!n, y_j\!=\!x_j$
et
$\displaystyle y_n=%
\sum_{j=1}^nx_j a_{m,\, j}a_{m,\,n}^{-1}$
 on a~:
$\displaystyle\sum_{j=1}^nx_j v_j=\sum_{j=1}^ny_j w_j$%
. Ainsi le morphisme
$g$
est aussi surjectif, envoit
$K^{n-1}=K^{(J\setminus\{n\})}$
dans
$K^{m-1}=K^{(J\setminus\{m\})}$
et sa restriction
$h=g_| : K^{n-1}=K^{(J\setminus\{n\})}\rightarrow%
K^{(J\setminus\{m\})}=K^{m-1}$%
, puisque si
$\displaystyle\sum_{j=1}^ny_jw_j=%
(z_1\ldots, z_{m-1}, 0)\in K^{m-1}=K^{(J\setminus\{m\})}$
alors
$y_m=0$%
, est aussi surjective. D'o{\`u} par
l'hypoth{\`e}se de r{\'e}currence
$m-1\leq n-1$
et donc
${\mathop{\rm Card\/}\nolimits}(I)=m-1+1\leq%
n-1+1={\mathop{\rm Card\/}\nolimits}(J)$%
.\hfill\findem
}
\Thc Corollaire|
Soit
$V, W$
des espaces vectoriels munis de bases
\hfill\break
$(b_i)_{i\in I}, b_i\in V;\ %
(c_j)_{j\in J}, c_j\in W$
et 
$f : W\rightarrow V$
un morphisme alors~:

$\quad({\romannumeral1} 1)$ Si
$f$ 
est surjectif
${\mathop{\rm Card\/}\nolimits}(J)\geq%
{\mathop{\rm Card\/}\nolimits}(I)$%
.

$\quad({\romannumeral1} 2)$ Si
$f$ 
est isomorphisme
${\mathop{\rm Card\/}\nolimits}(J)=%
{\mathop{\rm Card\/}\nolimits}(I)$%
.

$\quad({\romannumeral1} 3)$ Si
$f$ 
est injectif
${\mathop{\rm Card\/}\nolimits}(J)\leq%
{\mathop{\rm Card\/}\nolimits}(I)$%
.

$({\romannumeral2})$ Si
$(u_l)_{l\in L}, u_l\in V$
est une famille g{\'e}n{\'e}ratrice alors
${\mathop{\rm Card\/}\nolimits}(I)\leq
{\mathop{\rm Card\/}\nolimits}(L)$%

$({\romannumeral3})$ Si
$(v_k)_{k\in K}, v_k\in V$
est une famille libre alors
${\mathop{\rm Card\/}\nolimits}(K)\leq
{\mathop{\rm Card\/}\nolimits}(I)$%
.
\finc

\Defn Le cardinal  d'une base 
$(b_i)_{i\in I}, b_i\!\in\!V$
d'un espace vectoriel
$V$
sur le corps
$K$
ne d{\'e}pend pas du choix de la base,
c'est sa {\it dimension~:\/}
$${\mathop{\rm dim\/}\nolimits}(V)%
\build{=}_{}^{\hbox{\rm Syn}}%
{\mathop{\rm dim\/}_K\nolimits}(V)%
\build{=}_{}^{\hbox{\rm Syn}}%
[V:K]\build{=}_{}^{\hbox{\rm D{\'e}f}}
{\mathop{\rm Card\/}\nolimits}(I)$$

{\small\sl
Deux espaces vectoriels 
$V, W$
sont isomorphes
si et seulement si ils ont m{\^e}me dimension~:
$${\mathop{\rm dim\/}\nolimits}(V)=%
{\mathop{\rm dim\/}\nolimits}(W)$$
}
\vskip-5mm
\Thc co-Corollaire|
Soit
$U\subset V$
un sous-espace d'un espace vectoriel
$V$%
.
Alors
$${\mathop{\rm dim\/}\nolimits}(U)\leq 
{\mathop{\rm dim\/}\nolimits}(V)$$

De plus si
$V$
est de dimension finie, il y a {\'e}galit{\'e} si et
seulement si
$U=V$%
.
\finc
\rmc Exemples|
$({\romannumeral1})$ Un espace vectoriel
$V=\{0\}$
est l'espace nul si et seulement si
${\mathop{\rm dim\/}\nolimits}(V)=0$%
.

$({\romannumeral2})$ La dimension de
$K^n$
est
${\mathop{\rm dim\/}\nolimits}(K^n)=n$%
.

$({\romannumeral3})$ Si
$V^\ast$
est le dual d'un espace vectoriel 
$V$
de dimension finie\footnote{\small
voir {\bf 3.5.2 \uppercase\expandafter{\romannumeral3}}
pour le cas de dimension infinie.
}
alors\footnote{\small
d'apr{\`e}s l'exemple et la Proposition 1
de {\bf 1.5\/}.%
}~:
$${\mathop{\rm dim\/}\nolimits}(V^\ast)=%
{\mathop{\rm dim\/}\nolimits}(V)$$
\finc
\nssecp Le th{\'e}or{\`e}me du rang.|
\Thc Proposition|
Soit
$$0\rightarrow U%
\build{\rightarrow}_{}^{\iota}%
V%
\build{\rightarrow}_{}^{\pi}%
W\rightarrow0$$
une suite exacte courte d'espaces vectoriels.
Alors
$${\mathop{\rm dim\/}\nolimits}(V)=
{\mathop{\rm dim\/}\nolimits}(U)+%
{\mathop{\rm dim\/}\nolimits}(W)$$
{\small
Plus pr{\'e}cis{\'e}ment si
$\underline{b}\!=\!(b_i)_{i\in I}, b_i\!\in\!U, 
\underline{d}\!=\!(d_j)_{j\in J}, d_j\!\in\!W$
sont des bases des espaces extr{\^e}mes 
$U, W$
de la suite exacte et si, pour
$i\!\in\!I, c_i\!=\!\iota(b_i)\!\in\!V$
et pour
$j\!\in\!J, c_j\!\in\!\pi^{-1}(d_j)\subset V$%
,\hfill\break
 la famille
$\underline{c}\!=\!(c_k)_{k\in K=I\amalg J}, c_k\!\in\!V$
est une base de l'espace central
$V$
de la suite exacte.
}
\finc
{\small
\Dem Soit
$\displaystyle\sum_{k\in K}\lambda_k c_k\!=\!0$
une relation lin{\'e}aire de la famille
$\underline{c}$%
. 

Alors, comme pour
$i\!\in\!I, \pi(c_i)\!=\!\pi(\iota(b_i))\!=\!\pi\circ\iota(b_i)\!=\!0$
et, pour
$j\!\in\!J, \pi(c_j)\!=\!d_j$%
, on a
$\displaystyle0=\pi(\sum_{k\in K}\lambda_k c_k)\!=\!%
\pi(\sum_{i\in I}\lambda_i c_i)+%
\pi(\sum_{j\in J}\lambda_j b_j)=%
\sum_{i\in I}\lambda_i \pi(c_i)+%
\sum_{j\in J}\lambda_j \pi(c_j)\!=\!%
\sum_{j\in J}\lambda_j d_j$
et donc,\hfill\break
 puisque
$\underline{d}$
est libre pour
$j\!\in\!J$
on a 
$\lambda_j\!=\!0$%
. Ainsi
$\displaystyle
0\!=\!\sum_{i\in I}\lambda_i c_i\!=\!%
\sum_{i\in I}\lambda_i \iota(b_i)\!=\!%
\iota\bigl(\sum_{i\in I}\lambda_i b_i\bigr)$
et, puisque
$\iota$
est injective
$\displaystyle\sum_{i\in I}\lambda_i b_i\!=\!0$%
, d'o{\`u} puisque
$\underline{b}$
est libre, pour tout
$i\!\in\!I, \lambda_i\!=\!0$%
.

Ainsi pour tout
$k\!\in\!K\!=\!I\amalg J, \lambda_k\!=\!0$
et la famille
$\underline{c}$
est libre.\hfill\carre

Soit
$v\!\in\!V$ un vecteur de l'espace central.
Comme la famille
$\underline{d}$
est g{\'e}n{\'e}ratrice il y a 
$(\lambda_j)_{j\in J}\!\in\!K^{(J)}$
telle que
$\displaystyle\pi(v)\!=\!\sum_{j\in J}\lambda_j d_j$
donc le vecteur
$\displaystyle v'\!=\!v-\sum_{j\in J}\lambda_j c_j$%
, satisfaisant {\`a}~:\hfill\break
$\displaystyle\pi(v')\!=\!%
\pi(v)-\sum_{j\in J}\lambda_j \pi(c_j)\!=\!%
\pi(v)-\sum_{j\in J}\lambda_j d_j\!=\!0$
est dans le noyau
${\mathop{\rm ker\/}\nolimits}(\pi)\!=\!%
{\mathop{\rm Im\/}\nolimits}(\iota)$
et donc il y a
$u\!\in\!U$
tel que
$\iota(u)\!=\!v'$%
. Comme la famille
$\underline{b}$
est g{\'e}n{\'e}ratrice, il y a
$(\lambda_i)_{i\in I}\!\in\!K^{(I)}$
telle que
$\displaystyle u\!=\!\sum_{i\in I}\lambda_i b_i$%
, donc
$\displaystyle v'\!=\!\iota(u)\!=\!%
\sum_{i\in I}\lambda_i \iota(b_i)\!=\!%
\sum_{i\in I}\lambda_i c_i$%
. 

Ainsi
$\displaystyle v\!=\!v'+\sum_{j\in J}\lambda_j c_j\!=\!%
\sum_{k\in I\amalg J}\lambda_j c_j$
et la famille
$\underline{c}$
est aussi g{\'e}n{\'e}ratrice.\hfill\carre\findem
}
\Thc Corollaire|
Soit
$f : V\rightarrow W$
un morphisme d'espaces vectoriels. Alors~:
$${\mathop{\rm dim\/}\nolimits}(V)=%
{\mathop{\rm dim\/}\nolimits}(%
{\mathop{\rm ker\/}\nolimits}(f))+
{\mathop{\rm dim\/}\nolimits}(%
{\mathop{\rm Im\/}\nolimits}(f))$$
\finc
{\small
\Dem La suite
$0\rightarrow{\mathop{\rm ker\/}\nolimits}(f)%
\build{\rightarrow}_{}^{\subset}V%
\build{\rightarrow}_{}^{\overline{f}}%
{\mathop{\rm Im\/}\nolimits}(f)%
\rightarrow0$
est exacte.\hfill\findem
}
\Defn Le {\it rang\/} d'une application lin{\'e}aire
est la dimension de son image~:
$${\mathop{\rm rg\/}\nolimits}(f)=%
{\mathop{\rm dim\/}\nolimits}(%
{\mathop{\rm Im\/}\nolimits}(f))$$

\nssecp Application aux syst{\'e}mes lin{\'e}aires.|
\nsssecp Syst{\`e}mes lin{\'e}aires
\uppercase\expandafter{\romannumeral1}~:
le langage.| 
\Defn Soit
$(v_j)_{j\in J}, v_j\in V$
une famille indic{\'e}e par un ensemble
$J$
de vecteurs d'un espace vectoriel
et
$v\in V$
un autre vecteur de
$V$%
. {\it R{\'e}soudre le syst{\`e}me lin{\'e}aire\/}
$$\sum_{j\in J}X_j v_j=v$$
consiste {\`a} d{\'e}terminer les familles
$(x_j)_{j\in J}\in K^{(J)}$
d'{\'e}l{\'e}ments 
$x_j, j\in J$
presque tous nuls du corps
$K$
telles que l'on a l'{\'e}galit{\'e} de vecteurs
$\displaystyle\sum_{j\in J} x_j v_j=v$%
.

{\small
Avec les notations de {\bf 3.1\/}, il s'agit 
de d{\'e}terminer la pr{\'e}image de
$v$
par le  morphis\-me
${\mathop{cl_{\underline{v}}}\nolimits} :%
K^{(J)}\rightarrow V,\
\underline{\lambda}\!=\!(\lambda_i)_{j\in J}\mapsto%
{\mathop{cl_{\underline{v}}}\nolimits}%
(\underline{\lambda})\!=\!%
\sum_{j\in J}\lambda_j v_j$%
. 
C'est  l'ensemble vide si 
$v$
n'est pas dans le sous-espace engendr{\'e}
par la famille,
  une classe modulo le
sous-espace vectoriel 
${\mathop{\rm ker\/}\nolimits}(%
{\mathop{cl_{\underline{v}}}\nolimits})$
de l'espace vectoriel libre
$K^{(J)}$
sur
$J$
sinon. Si
$J$
est fini et donn{\'e}, la dimension de ce noyau, selon
le th{\'e}or{\`e}me du rang
(Corollaire de {\bf 3.3\/})
d{\'e}termine
et est d{\'e}termin{\'e}e par le rang de
${\mathop{cl_{\underline{v}}}\nolimits}$
.
}

Hormi cette pr{\'e}sentation on se limitera au cas
o{\`u} l'espace vectoriel
$V$
est de dimension finie, puis par choix d'une de
ses bases
$(b_i)_{i=1,\ldots,m}$
{\`a}
$V=K^m$
et o{\`u} l'ensemble
$J$
est fini puis, quitte {\`a} le num{\'e}roter
$J=\{1,\ldots, n\}$%
. 

En ce cas, en {\'e}crivant les coordonn{\'e}es
des vecteurs en colonnes\hfill\break
$v_j=(a_{1,\,j},\ldots,a_{m,\,j})^t,\ %
j=1,\ldots,m; v=(b_1,\ldots, b_m)^t\in K^m$
c'est le syst{\`e}me
$$\left\{%
\begin{matrix}X_1 a_{1,\, 1}&+\ldots+&X_k a_{1,\, k}&%
+\ldots+&X_n a_{1,\, n}&=&b_1\\
\vdots& &\vdots& &\vdots&&\vdots\\
X_1 a_{i,\, 1}&+\ldots+&X_k a_{i,\, k}&%
+\ldots+&X_n a_{i,\, n}&=&b_i\cr\vdots& &\vdots& &\vdots& &\vdots\\
X_1 a_{m,\, 1}&+\ldots+&X_k a_{m,\, k}&%
+\ldots+&X_n a_{m,\, n}&=&b_m
\end{matrix}
\right.\leqno{(\Sigma)}$$
\Defn Le {\it rang\/} d'une famille
$\underline{v}=(v_i)_{i\in I}, v_i\in V$
est 
${\mathop{\rm rg\/}\nolimits}(\underline{v})=%
{\mathop{\rm dim\/}\nolimits}(<\underline{v}>)$%
, \hfill\break
la dimension du sous-espace engendr{\'e}.
\Thc Proposition|Le syst{\`e}me
$(\Sigma)$
a une solution si et seulement si
$${\mathop{\rm rg\/}\nolimits}(v_1,\ldots, v_n)=%
{\mathop{\rm rg\/}\nolimits}(v_1,\ldots, v_n, v)$$
Ses solutions forment alors une classe
modulo un sous-espace de dimension~:
$$n-{\mathop{\rm rg\/}\nolimits}(v_1,\ldots, v_n)$$
\finc
\vfill\eject
\nsssecp Syst{\`e}mes lin{\'e}aires
\uppercase\expandafter{\romannumeral2}~:
rang du syst{\`e}me transpos{\'e}.| 
Le syst{\`e}me
$(\Sigma)$
est associ{\'e} {\`a} la matrice
$A=\begin{pmatrix}
a_{1,\, 1}&\cdots&a_{1,\, }\\
\vdots&&\vdots\\
a_{m,\, 1}&\cdots&a_{m,\, n}
\end{pmatrix}%
\in{\mathop{\rm M\/}\nolimits}_{m,\, n}$%
.
{\small
\Defns Les {\it rangs de colonne\/}
(resp. {\it ligne\/}) {\it  {\`a} droite\/}
(resp. {\it gauche\/}) de
$A$
sont les dimensions~: 
$${\mathop{\rm rg\/}\nolimits}_d^c(A)=%
{\mathop{\rm dim\/}\nolimits}(<c_1,\ldots, c_n>_d),\quad
{\mathop{\rm rg\/}\nolimits}_s^l(A)=%
{\mathop{\rm dim\/}\nolimits}(<l_1,\ldots, l_m>_s)$$
$${\mathop{\rm rg\/}\nolimits}_s^c(A)=%
{\mathop{\rm dim\/}\nolimits}(<c_1,\ldots, c_n>_s),\quad
{\mathop{\rm rg\/}\nolimits}_d^l(A)=%
{\mathop{\rm dim\/}\nolimits}(<l_1,\ldots, l_m>_s)$$
des sous-espace vectoriels {\`a} droite 
(resp. {\`a} gauche) de
$K^m$
(resp. 
$K^n$%
)
engendr{\'e} par\hfill\break
 les colonnes
$c_j\!=\!(a_{1,\,j},\ldots,a_{m,\,j})^t\in K^m$
(resp. lignes
$l_i=(a_{i,\,1},\ldots,a_{i,\,n})\in K^n$%
)
.
}

\Thc Proposition|
On a les {\'e}galit{\'e}s~:
$${\mathop{\rm rg\/}\nolimits}_d^c(A)=%
{\mathop{\rm rg\/}\nolimits}_s^l(A),\quad
{\mathop{\rm rg\/}\nolimits}_s^c(A)=%
{\mathop{\rm rg\/}\nolimits}_d^l(A)$$
\finc
{\small
\Dem Posons
$p={\mathop{\rm rg\/}\nolimits}_s^l(A)$%
. Quitte {\`a} renum{\'e}roter les lignes on peut
supposer que l'espace 
$<l_1,\ldots,l_m>_s=<l_1,\ldots, l_{p}>_s$
des lignes est engendr{\'e} par les
$p$
premi{\`e}res lignes.\hfill\break
Cette renum{\'e}rotation correspondant {\`a} un 
isomorphisme de permutation
 des coordonn{\'e}es de l'espace
$K^m$
des colonnes, ne change pas le rang {\`a} droite
${\mathop{\rm rg\/}\nolimits}_d^c(A)$
de colonne.

Soit
$A_{p,\,n}\in%
{\mathop{\rm M\/}\nolimits}_{p,\,n}$
la matrice des 
$p$
premi{\`e}res lignes de
$A$%
.

Comme les autres lignes 
$\displaystyle
l_k\!=\!\sum_{i=1}^{p}c_{k,\,i}l_i, p\!<k\!\leq\!m$
sont combinaisons lin{\'e}aires
des premi{\`e}res, toute relation lin{\'e}aire
({\`a} droite)
$\displaystyle\sum_{j=1}^n a_{i\,j} x_j\!=\!0, %
 1\!\leq\!i\!\leq\!p$
entre les colonnes de
$A_{p,\, n}$
donne pour
$\displaystyle p<k\leq m,\ 
\sum_{j=1}^na_{k,\, j} x_j=\sum_{j=1}^n%
\bigl(\sum_{i=1}^{p}c_{k\, i} a_{i,\,j}\bigr) x_j=%
\sum_{i=1}^{p}c_{k\, i}\bigl(%
\sum_{j=1}^n a_{i,\,j} x_j\bigr)=0$%
, la m{\^e}me relation entre les colonnes compl{\`e}tes de
$A$%
. Ainsi puisque\footnote{\small
d'apr{\`e}s l'exemple de {\bf 3.2\/}.
}
${\mathop{\rm rg\/}\nolimits}_d^c(A_{p,\, n})%
\leq p$
on a~:
${\mathop{\rm rg\/}\nolimits}_d^c(A)\leq p=%
{\mathop{\rm rg\/}\nolimits}_s^l(A)$%
.

Soit
$A^t=(a_{j,\, i})_{\begin{matrix}
                    1\!\leq\!j\!\leq\!n\\
                    1\!\leq\!i\!\leq\!m
                \end{matrix}
                 }\in%
{\mathop{\rm M\/}\nolimits}_{n,\,m}$
la matrice dont les lignes (resp. colonnes)
sont les colonnes (resp. lignes) de
$A$%
. L'in{\'e}galit{\'e} pr{\'e}c{\'e}dante appliqu{\'e}e
une fois {\`a}
$A$%
, puis {\`a}
$A^t$
et au corps oppos{\'e}
$K^\circ$
donne~:
${\mathop{\rm rg\/}\nolimits}_d^c(A)\leq%
{\mathop{\rm rg\/}\nolimits}_s^l(A)=%
{\mathop{\rm rg\/}\nolimits}_s^c(A^t)\leq%
{\mathop{\rm rg\/}\nolimits}_d^l(A^t)=%
{\mathop{\rm rg\/}\nolimits}_d^c(A)$%
, d'o{\`u} la premi{\`e}re {\'e}galit{\'e}
${\mathop{\rm rg\/}\nolimits}_d^c(A)=%
{\mathop{\rm rg\/}\nolimits}_s^l(A)$%
. La seconde suit en appliquant la premi{\`e}re
dans le corps oppos{\'e}
$K^\circ$%
.\hfill\findem
}
\Thc Corollaire|
Soit
$f : U\rightarrow V$
un morphisme entre deux espaces vectoriels
de dimension finie et
$f^\ast : V^\ast\rightarrow U^\ast$
son morphisme transpos{\'e} alors~:
$${\mathop{\rm rg\/}\nolimits}(f)=%
{\mathop{\rm rg\/}\nolimits}(f^\ast)$$
\finc
\rmc Remarque|
Il n'y a pas {\'e}galit{\'e} entre les quatre rangs
d'une matrice~: si
$a, b\in K, ab\ne ba$
la matrice
$A\!=\!\begin{pmatrix}
            1&b\\
            a&ab
            \end{pmatrix}$
v{\'e}rifie\footnote{\small
puisque
$a\begin{pmatrix}
           1&b
  \end{pmatrix}%
\!=\!\begin{pmatrix}
          a&ab
      \end{pmatrix}$%
, mais 
$\begin{pmatrix}
  1&b
\end{pmatrix}
a\!=\!\begin{pmatrix}
             a&ba
\end{pmatrix}%
\ne%
\begin{pmatrix}
        a&ab
\end{pmatrix}=a\begin{pmatrix}
                 1&b
                \end{pmatrix}$%
.
}~:
${\mathop{\rm rg\/}\nolimits}_s^l(A)\!=\!1\!\ne2\!=\!%
{\mathop{\rm rg\/}\nolimits}_d^l(A)$%
.
\finc
\rmc Exemple|
Soit
$H=
\{z_{a, b}=\begin{pmatrix}
              a&b\\
              -\overline{b}&\overline{a}
              \end{pmatrix} | %
 a, b\in{\bf C\/}\}\subset%
{\mathop{\rm M\/}\nolimits}_2({\bf C\/})$%
, c'est un sous-anneau non commutatif\footnote{\small
$a=\begin{pmatrix}
          0&1\\
         -1&0
\end{pmatrix}%
, b=\begin{pmatrix}
       0&i\\
       i&0
       \end{pmatrix}$ 
v{\'e}rifient
$ab=\begin{pmatrix}
          i&0\\
          0&-i
\end{pmatrix}%
=-ba\ne ab$%
.
}
qui, puisque la transconjugaison
$\sigma : \begin{pmatrix}
              a&b\\
              -\overline{b}&\overline{a}
\end{pmatrix}%
\mapsto
\begin{pmatrix}
  \overline{a}&-b\\
  \overline{b}&a
  \end{pmatrix}$
est un anti-automorphisme\footnote{\small
{\it c. a d.} pour tout
$x, y\in H, \sigma(x-y)=\sigma(x)-\sigma(y)$
et
$\sigma(xy)=\sigma(y)\sigma(x)$%
.
}
v{\'e}rifiant pour tout
$z=z_{a,\, b}\in H, z\sigma(x)=|a|^2+|b|^2$%
, est un corps, le {\it corps des quaternions\/},
 l'inverse d'un {\'e}l{\'e}ment
$z=z_{a,\,b}\ne0$
non nul {\'e}tant
$z^{-1}=(|a|^2+|b|^2)^{-1}\sigma(z)$%
.
\finc
\vskip3mm
{\small
\centerline{\small\petcap
 Commentaires bibliographiques\/}
\vskip3mm

Corps et espaces vectoriels sont tr{\`e}s
efficacement pr{\'e}sent{\'e}s dans le chap.
\uppercase\expandafter{\romannumeral1}\footnote{\small
d'o{\`u} (Theorem 4) le pr{\'e}sent  {\bf 3.4.2\/}
a {\'e}t{\'e} pirat{\'e}!
} de {\bf [Ar2]}.

Pour r{\'e}soudre effectivement les syst{\'e}mes lin{\'e}aires et la
\og m{\'e}thode de Gauss sur un corps\fgf, outre C et E du
chap. \uppercase\expandafter{\romannumeral1} de {\bf [Ar2]},
voir A2 de {\bf [Ga]\/}
(en particulier l'{\'e}nonc{\'e} d'unicit{\'e} de {\bf 2.7\/}.)

Voir au \S3 du chap. 7 de {\bf [He2]\/} une preuve du th{\'e}or{\`e}me
de Fr{\"o}benius {\'e}tablissant que
{\sl le corps
$H$
de l'exemple est {\`a} isomorphisme pr{\`e}s le seul
corps non commutatif contenant les r{\'e}els comme sous-corps,
en {\'e}tant\footnote{\small
voir {\bf 3.5.1\/} et {\bf 3.5.2\/} ci-dessous pour
cette alg{\`e}bre lin{\'e}aire \og {\`a} plusieurs corps\fgf.
} 
de dimension
finie sur les r{\'e}els.\/}

Plus de d{\'e}tails sur matries et  espaces vectoriels
sont dans ~:

les \S14 {\`a} \S20 de {\bf [Go]\/}.

le \S3 (pour la 
$1^\circ$
 ({\oldstyle 1947}) {\'e}dition,
\S7 apr{\`e}s la  troisi{\`e}me ({\oldstyle 1959}) {\'e}dition 
du Chap. \uppercase\expandafter{\romannumeral2} de {\bf [Bo]\/}.\hfill\break
la dualit{\'e} est au \S4  (%
$1^\circ$%
) ({\oldstyle 1947}) {\'e}dition,
\S2 apr{\'e}s la troisi{\`e}me ({\oldstyle 1959})
 {\'e}dition 
du Chap. \uppercase\expandafter{\romannumeral2} de {\bf [Bo]\/}.

Kap. 4 de {\bf [vW]\/}, chap. \uppercase\expandafter{\romannumeral1}
de vol. \uppercase\expandafter{\romannumeral2} de {\bf [Ja]\/}
ou \S 5 de
chap. \uppercase\expandafter{\romannumeral3} de {\bf [Ld]}.

}

\npage
\nssecp Structure rationnelle sur un sous-corps.|
\nsssecp 
D{\'e}finitions exemples et
propri{\'e}t{\'e}s.|
{\small
Dans cette sous-section on consid{\`e}re
$k\subset K$
un sous-corps du corps
$K$%
.
Ainsi tout
$K$%
-espace vectoriel
$U$
est un
$k$%
-espace vectoriel, en particulier le corps
$K$
est un 
$k$%
-espace vectoriel. 
}

\Defns 
Un 
$k$%
-sous-espace  vectoriel
$U_k\subset V$
d'un
$K$%
-espace vectoriel
$V$
est
{\it
$(K,k)$%
-libre\/}
si toute famille
$\underline{u}=(u_i)_{i\in I}, u_i\in U_k$
libre sur
$k$
est, consid{\'e}r{\'e}e comme famille
$(u_i)_{i\in I}, u_i\in V$%
, libre sur 
$K$%
. Le sous-%
$K$%
-espace vectoriel~:
$$U=U_K=<U_k>_K\subset V$$
engendr{\'e} (sur %
$K$%
) par
$U_k$
est dit {\it port{\'e} par le sous-%
$k$%
-espace-%
$(K,k)$%
-libre\/}
$U_k\subset U$%
. 

Si
$U_K=V$
on dit que l'espace 
$V$
est {\it muni de la 
$k$%
-structure\/}
$V_k=U_k\subset V$%
.

Un sous-espace
$U\subset V$
est dit {\it 
$k$%
-rationnel\/}
(pour la
$k$%
-structure
$V_k\subset V$%
) 
 si il est engendr{\'e} sur
$K$
par le sous-%
$k$%
-espace vectoriel
$U_k=U\cap V_k\subset V_k$
de la
$k$%
-structure.
Ce sous-espace
$U_k$
est dit
{\it espace des
$k$%
-points du\/} sous-espace
$k$%
-rationnel
$U$%
.

\rmc Remarques et exemples|
a) Un sous-%
$k$%
-espace
$U_k\!\subset\!V$
est
$(K, k)$%
-libre si et seulement si pour toute famille
$(u_i)_{i\in I}, u_i\!\in\!U_k$%
, il est {\'e}quivalent qu'elle soit, comme famille de
$U_k$%
, li{\'e}e sur
$k$
et que, consid{\'e}r{\'e}e comme
famille 
$(u_i)_{i\in I}, u_i\!\in\!V$
de
$V$%
,  elle soit li{\'e}e sur
$K$%
.

b) Un sous-%
$k$%
-espace
$W_k\!\subset\!U_k\!\subset\!V$
d'un sous-espace
$(K, k)$%
-libre est
$(K, k)$%
-libre,\hfill\break
 en particulier l'espace des
$k$%
-points d'un sous-espace
$k$%
-rationnel est
$(K, k)$%
-libre.

 c) Si 
$(b_i)_{i\in I}, b_i\!\!\in\!\!V$
est
$K$%
-base de
$V\!$%
, alors
$V_k\!\!=\!\!\oplus_{i\in I}k b_i\!\!\subset\!\!%
\oplus_{i\in I}K b_i\!\!=\!\!V\!\!$
est 
$\null\!\!k$%
-structure de 
$\!V\!$%
.
 R{\'e}ciproquement\footnote{\small
car une
$k$%
-base
d'une 
$k$%
-structure
$(c_i)_{i\in I}, c_i\!\in\!V_k\!\subset\!V$
engendre
$V\!=\!(V_k)_K$
et est libre sur
$K$%
.
}
toute 
$k$%
-structure sur un
$K$%
-espace vectoriel
$V$
est obtenue de cette mani{\`e}re~:

c1) {\sl Si 
$V_k\subset V$
est une
$k$%
-structure sur un
$K$%
-espace vectoriel
$V$
alors pour tout 
$K$-espace vectoriel
$W$
la restriction
${\mathop{\rm Hom\/}\nolimits}_K(V, W)%
\rightarrow%
{\mathop{\rm Hom\/}\nolimits}_k(V_k, W)$
est un isomorphisme du groupe des applications
$K$%
-lin{\'e}aires de
$V$
dans
$W$
dans celui des applications
$k$%
-lin{\'e}aires de
$V_k$
dans
$W$%
.
}
{\small
\Dem Soit
$(b_i)_{i\in I}$
une base de
$V_k$%
. Ainsi 
$V_k$
et
$V$
sont respectivement
$k$%
-espace
et
$K$%
-espace libre sur cette m{\^e}me base.
 Donc, d'apr{\`e}s leur
propri{\'e}t{\'e} universelle
(Proposition 3 de {\bf 1.4\/}),
les restrictions {\`a} cette base commune
${\mathop{\rm Hom\/}\nolimits}_K(V, W)\ni f%
\mapsto f_{|\{b_i | i\in I\}}\in W^I$
et
${\mathop{\rm Hom\/}\nolimits}_k(V_k, W)\ni g%
\mapsto g_{|\{b_i | i\in I\}}\in W^I$
sont des bijections.\hfill\findem
}

c2) {\sl Si
$U\!\subset\!V$
est un sous-espace
$k$%
-rationnel pour la 
$k$%
-structure
$V_k\!\subset\!V$
alors
$V_k/U\!\cap\!V_k\!\subset\!V/U$
est une
$k$%
-structure, dite
{\it
$k$%
-structure quotient\/}
sur le quotient
$V/U$%
.
}
{\small
\Dem Soit
$\pi : V\rightarrow V/U$
l'application quotient et 
$(b_l)_{i\in L=I\amalg J}, b_l\in V_k$
une base de
$V_k$
compl{\'e}tant une base
$(b_i)_{i\in I}, b_i\!\in\!U_k$
de
$U_k$
alors,
les restrictions de
$\pi$
aux suppl{\'e}mentaires
$W_k\!=\!\oplus_{j\in J}k b_j\!\subset\!V_k\!=\!%
U_k\oplus W_k$
et
$W\!=\!\oplus_{j\in J}K b_j\subset V=U\oplus W$
de
$U_k$
et
$U$
dans
$V_k$
et
$V$
induisent des isomorphismes sur les espaces quotients
$V_k/U_k\!=\!V_k/U\cap V_k\!\subset\!V/U$
et
$V/U$
respectivement. 

Ainsi
$V_k/U_k\!=\!V_k/U_k\cap V=\pi(\oplus_{l\in L} k b_l)%
\!=\!\pi(\oplus_{j\in J}k b_j)\subset%
\pi(\oplus_{j\in J}K b_j)\!=\!%
\pi(\oplus_{l\in L}K b_l)\!=\!V/U$
est la 
$k$%
-structure du quotient
$V/U$
associ{\'e} {\`a} sa base
$(\pi(b_j))_{j\in J}, \pi(b_j)\in V/U$%
.\hfill\findem
}

\finc
\Thc Proposition|
Soit
$\underline{b}\!=\!(b_j)_{j\in J}, b_j\!\in\!V$
une base d'un 
$K$%
-espace vectoriel
$V$
\hfill\break
et
$\underline{\beta}\!=\!(\beta_i)_{i\in I}, \beta_i\!\in\!K$
une base du
$k$%
-espace vectoriel (aussi {\`a} gauche)
$K$%
. 

Alors
$\underline{\beta b}\!=\!%
(\beta_i b_j)_{(i,\, j)\in I\times J}, %
\beta_i b_j\!\!\in\!\!V$
est une base du
$k$%
-espace vectoriel
$V$%
, ainsi~:
$${\mathop{\rm dim\/}\nolimits}_k(V)=%
{\mathop{\rm dim\/}\nolimits}_k(K)\cdot%
{\mathop{\rm dim\/}\nolimits}_K(V)$$
\finc

On se donne
$V_k\!\subset\! V, W_k\!\subset\!W$
deux
$K$%
-espaces vectoriels munis de
$k$%
-structure.
\Defn 
 Un morphisme
$f\in{\mathop{\rm Hom\/}\nolimits}(V, W)$
est dit {\it d{\'e}fini sur 
$k$%
\/} si 
$f(V_k)\subset W_k$%
. On note alors
$f_k\in{\mathop{\rm Hom\/}\nolimits}(V_k, W_k)$
le
{\it 
$k$%
-morphisme d{\'e}finissant
$f$%
.\/}
\Thc Th{\'e}or{\`e}me|
Soit
$f : V\rightarrow W$
un morphisme d{\'e}fini sur
$k$%
. Alors

$({\romannumeral1})\quad f(V)\cap W_k=f(V_k)$%
.

$({\romannumeral2})\quad%
{\mathop{\rm ker\/}\nolimits}(f)=%
<{\mathop{\rm ker\/}\nolimits}(f_k)>_K$%
.

En particulier 
$f$
et
$f_k$
 ont m{\^e}me rang et
m{\^e}me dimension de noyau~:
$${\mathop{\rm rg\/}\nolimits}_K(f)=%
{\mathop{\rm rg\/}\nolimits}_k(f_k),\quad%
{\mathop{\rm dim\/}\nolimits}_K(%
{\mathop{\rm ker\/}\nolimits}(f))=%
{\mathop{\rm dim\/}\nolimits}_k(%
{\mathop{\rm ker\/}\nolimits}(f_k))$$
.
\finc
{\small
\Dem Soit
$(b_i)_{i\in I}, b_i\!\in\!V_k$
une base de
$V_k$%
. Un vecteur
$w\!\in\!W_k$
est dans 
$f(V)$%
, l'image de
$f$
si et seulement si la famille
$(w_j)_{j\in I\amalg\{w\}}, w_i\!=\!%
b_i, w_w\!=\!w\!\in\!W$
est li{\'e}e dans
$W$%
. 

Comme  famille du sous-espace
$(K, k)$%
-libre
$W_k$%
, cela arrive si et seulement si cette famille est
li{\'e}e dans
$W_k$%
,{\it c. a d.\/} si et seulement si
$w$
est dans l'image de
$f_k$%
. d'o{\`u}
$({\romannumeral1)}$%
.\hfill\carre

L'inclusion
$<\!{\mathop{\rm ker\/}\nolimits}(f_k)\!>_K%
\!\subset\!{\mathop{\rm ker\/}\nolimits}(f)$
{\'e}tant imm{\'e}diate, 
 il suffit\footnote{\small
quitte {\`a}
changer
$V, V_k$
en
$ V/\!%
<\!{\mathop{\rm ker\/}\nolimits}(f_k)\!>_K,%
 V_k/%
{\mathop{\rm ker\/}\nolimits}(f_k)$
et
$f_k$
en
$\overline{f_k} : V_k/%
{\mathop{\rm ker\/}\nolimits}(f_k)\rightarrow W_k$%
.} de consid{\'e}rer le cas
${\mathop{\rm ker\/}\nolimits}(f_k)\!=\!0$%
. Ainsi
$f_k$
est injective et, si
$(b_i)_{i\in I}, b_i\!\in\!V_k$
est une base de
$V_k$
alors la famille
$(f(b_i))_{i\in I}, f(b_i)\!\in\!W_k$
est libre dans
$W_k$
donc dans
$W$%
, puisque
$W_k$
est
$(K, k)$%
-libre dans
$W$
et
$f$
est injective.
\hfill\carre\findem
}
\Thc Corollaire|
Un syst{\`e}me lin{\'e}aire de
$m$
{\'e}quations {\`a}
$n$
inconnues {\`a} coefficients
et second membre dans un sous-corps
$k\subset K$
d'un corps
$K$
a une solution dans 
$k^n$
si et seulement si il a une solution dans
$K^n$%
. En ce cas ces ensembles de solutions sont
des classes de sous-%
$k$%
(resp. 
$K$)%
-espaces de
$k^n$%
(resp. 
$K^n$%
)
de m{\^e}me dimension.
\finc
\vskip10mm
{\small
\centerline{\small\petcap
 Commentaires bibliographiques\/}
\vskip3mm

Pour plus de d{\'e}tails voir le  \S3 de a $1^\circ$
 ({\oldstyle 1947}) {\'e}dition 
du Chap. \uppercase\expandafter{\romannumeral2} de {\bf [Bo]\/},
qui se base sur la tr{\`e}s concr{\`e}te
notion de {\it  solution primordiale\/}\footnote{\small
{\it c. a d.\/} solution dont l'ensemble des indices des coordonn{\'e}es 
non nulles 
est minimal et dont (au moins) une des coordonn{\'e}e vaut
$1$%
. Les r{\'e}sultat sont~:\hfill\break
{\sl L'espace des solutions est engendr{\'e} par les solutions primordiales.\/}\hfill\break
{\sl Les solutions primordiales d'un syst{\`e}me {\`a} coefficients dans 
$k\!\subset\!K$%
, sont {\`a} coefficients dans
$k$%
.
}
} %
 d'un syst{\`e}me lin{\'e}aire
homog{\`e}ne.
Notion abandonn{\'e}e dans les {\'e}ditions 
post{\'e}rieures {\`a } la trois{\`e}me qui, faisant
usage du produit tensoriel, ont un traitement 
des structures rationnelles (aux \S5 et \S8) moins {\'e}l{\'e}mentaire
qu'aux pr{\'e}c{\'e}dentes.

}
\npage
\nsssecp Applications~:
{\uppercase\expandafter{\romannumeral1}}
Quelques propri{\'e}t{\'e}s des corps gauches.|
\Defn Le {\it centre\/} d'un corps
$K$
est l'ensemble~: 
$${\mathop{\rm Z\/}\nolimits}(K)=%
\{z\in K\ |\ %
\forall \lambda\in K, z\lambda=\lambda z\}$$
des {\'e}l{\'e}ments de
$K$
commutant {\`a} tout {\'e}l{\'e}ment de
$K$%
.

 Plus g{\'e}n{\'e}ralement, le
{\it commutant dans
$K$%
\/} d'une partie
$L\subset K$
est~:
$${\mathop{\rm C\/}\nolimits}_K(L)=
\{z\in K\ |\ %
\forall \lambda\in L, z\lambda=\lambda z\}$$
On a
${\mathop{\rm Z\/}\nolimits}(K)\!=\!%
{\mathop{\rm C\/}\nolimits}_K(K)$%
, et ces
${\mathop{\rm C\/}\nolimits}_K(L)$
 sont\footnote{\small
Car 
$0\!\in\!{\mathop{\rm C\/}\nolimits}(L)$
et si 
$z_1, z_2, z\!\in\!{\mathop{\rm C\/}\nolimits}(L)$
avec 
$z\!\ne\!0$%
, pour tout
$\lambda\!\in\!K$
on a
$\lambda(z_1-z_2)\!=\!\lambda z_1-\lambda z_2\!=\!%
z_1\lambda-z_2\lambda\!=\!(z_1-z_2)\lambda,
\lambda(z_2z_1)\!=\!(\lambda z_2)z_1\!=\!%
(z_2\lambda)z_1\!=\!z_2(\lambda z_1)\!=\!%
z_2(z_1\lambda)\!=\!(z_2z_1)\lambda$
et 
$ z^{-1}\lambda\!=\!z^{-1}\lambda(zz^{-1})\!=\!%
z^{-1}(\lambda z)z^{-1}\!=\!%
z^{-1}(z\lambda )z^{-1}\!=\!%
(z^{-1}z)\lambda z^{-1}\!=\!\lambda z^{-1}$%
, donc
$z_1-z_2, z_2z_1, z^{-1}\!\in\!%
{\mathop{\rm C\/}\nolimits}(L)$%
.
}
 des sous-corps de
$K$%
.

Soit
$K$
un corps de centre
$k={\mathop{\rm Z\/}\nolimits}(K)$%
. 

{\small
Consid{\'e}rons le groupe des endomorphismes 
${\mathop{\rm End\/}\nolimits}_k(K_s)$
du
$k$%
-espace vectoriel {\`a} gauche
$K$%
. La multiplication {\`a} droite
${\mathop{\rm End\/}\nolimits}_k(K)\!\times\!\!K\!\ni\!%
(f, b)\!\mapsto\!f b : x\!\mapsto\!f(x)b$
en fait un
$K$%
-espace vectoriel {\`a} droite.
Il contient, 
$k$
{\'e}tant dans le centre de
$K$%
,
le sous-%
$k$%
-espace vectoriel 
$\sigma_K$
(isomorphe {\`a}
$K$%
) des multiplications {\`a} gauche
(lin{\'e}aires {\`a} droite)
$\sigma_a : K\rightarrow K, x\mapsto ax$
par les {\'e}l{\'e}ments 
$a\in K$%
. 
}
\Thc Lemme 1|
Ce sous-espace vectoriel
$\sigma_K$
des multiplications {\`a} gauche
par les {\'e}l{\'e}ments de
$K$
est 
$(K, k)$%
-libre dans le
$K$%
-espace vectoriel {\`a} droite
${\mathop{\rm End\/}\nolimits}_k(K_s)$
des
$k$%
-endomorphismes lin{\'e}aires ({\`a} gauche) de
$K$%
.
\finc
{\small
\Dem Soit
$\underline{\sigma}\!=\!(\sigma_{a_i})_{i\in I}, %
\sigma_{a_i}\in\sigma_K,a_i\!\in\!K$
une famille li{\'e}e de
$\sigma_K$
et
$\displaystyle\sum_{j=1}^n\sigma_{a_{i_j}}b_{i_j}\!=\!0$
une relation lin{\'e}aire non triviale de
longueur
$n$
minimale.
Quitte {\`a} remplacer 
$b_{i_j}$
par 
$b_{i_j}b_{i_n}^{-1}$%
, on peut supposer
$b_{i_n}\!=\!1$%
. Alors pour tout
$c\!\in\!K$
et
$x\in K$
on a
$\displaystyle0\!=\!%
\sum_{j=1}^n\sigma_{a_{i_j}}(xc)b_{i_j}\!=\!%
\sum_{j=1}^n\sigma_{a_{i_j}}(x)cb_{i_j}$
d'o{\`u}
$\displaystyle0\!=\!%
\sum_{j=1}^n\sigma_{a_{i_j}}cb_{i_j}-%
\sum_{j=1}^n\sigma_{a_{i_j}}b_{i_j}c\!=\!%
\sum_{j=1}^{n}\sigma_{a_{i_j}}(x)(cb_{i_j}-b_{i_j}c)%
\!=\!%
\sum_{j=1}^{n-1}\sigma_{a_{i_j}}(x)(cb_{i_j}-b_{i_j}c)$
et par minimalit{\'e} de
$n$
on a pour tout
$c\!\in\!K$
et
$j\!=\!1,\ldots, n, cb_{i_j}-b_{i_j}c\!=\!0$
et donc
$b_{i_j}\!\in\!{\mathop{\rm Z\/}\nolimits}(K)\!=\!k$%
.\hfill\findem
}
\Thc Compl{\'e}ment|
Soit
$(a_i)_{i\in I}, a_i\!\in\!K$
une 
$k$%
-base et 
$(b_i)_{i\in I}, b_i\!\in\!K$
une famille d'{\'e}l{\'e}ments presque tous nuls de
$K$%
.

Alors, si
$k\subset L\subset K$
est un sous corps de
$K$
contenant
$k$%
, l'endomorphisme
$$\displaystyle\sum_{i\in I}\sigma_{a_i}b_i\!\in\!%
{\mathop{\rm End\/}\nolimits}_k(K_s)$$
est
$L$%
-lin{\'e}aire [{\it c. a d.\/}
$\displaystyle\sum_{i\in I}\sigma_{a_i}b_i\!\in\!%
{\mathop{\rm End\/}\nolimits}_L(K_s)\!\subset\!%
{\mathop{\rm End\/}\nolimits}_k(K_s)$%
] si et seulement si~:\hfill\break
 pour tout
$i\in I$
le coefficient
$b_i\!\in\!{\mathop{\rm C\/}\nolimits}_K(L)$
est dans le commutant de
$L$
dans 
$K$%
.
\finc
{\small
\Dem 
$\displaystyle\sum_{i\in I}\sigma_{a_i}b_i\!\in\!%
{\mathop{\rm End\/}\nolimits}_L(K)$
si et seulement si pour tout
$\lambda\!\in\!L$
et
$x\!\in\!K$
on a
$\displaystyle%
\sum_{i\in I} \sigma_{a_i}(x)b_i\lambda\!=\!%
\bigl(\sum_{i\in I}\sigma_{a_i}(x)b_i\bigr)%
\lambda\!=\!%
\sum_{i\in I}\sigma_{a_i}(x\lambda)b_i\!=\!%
\sum_{i\in I}a_i(x\lambda)b_i\!=\!%
\sum_{i\in I}(a_ix)\lambda b_i\!=\!%
\sum_{i\in I}\sigma_{a_i}(x)\lambda b_i$%
, donc
$\displaystyle%
\sum_{i\!\in\!I}\sigma_{a_i}b_i\lambda\!=\!%
\sum_{i\!\in\!I}\sigma_{a_i}\lambda b_i$
et, puisque d'apr{\`e}s le Lemme la famille
$(\sigma_{a_i})_{i\in I}$
est libre sur
$K$%
, pour tout
$i\!\in\!I$%
et
pour tout
$\lambda\!\in\!L$%
, on a
$b_i\lambda=\lambda b_i$%
, {\it c. a d.\/} pour tout
$i\!\in\!I, b_i\!\in\!%
{\mathop{\rm C\/}\nolimits}_K(L)$%
.\hfill\findem
}
\vskip-1mm
On suppose d{\'e}sormais le corps
$K$
de dimension\footnote{\small
Comme
$k={\mathop{\rm Z\/}\nolimits}(K)$
les structure 
$K_s, K_d$
de
$k$%
-espace vectoriel {\`a} gauche et {\`a} droite
sur
$K$
co{\"\i}ncident et cette dimension sur
$k$
ne d{\'e}pend
de ce \og choix de c{\^o}t{\'e} des scalaires\fgf.
\hfill\break
De m{\^e}me si
$L\!\subset\!K$
est un sous-corps commutatif de
$K$
les structures {\`a} gauche et {\`a} droite de
$k\cap L$%
-espaces vectoriels sur
$k$
et
$L$
co{\"\i}ncident et les relations 
de la Proposition 
{\bf 3.5.1\/}
$${\mathop{\rm dim\/}\nolimits}_{k\cap L}(K)\!=\!%
{\mathop{\rm dim\/}\nolimits}_{k\cap L}(L)\cdot%
{\mathop{\rm dim\/}\nolimits}_{L}(K)\!=\!%
{\mathop{\rm dim\/}\nolimits}_{k\cap L}(k)\cdot%
{\mathop{\rm dim\/}\nolimits}_{k}(K)$$
donnent que les dimensions sur
$L$
ne d{\'e}pendent pas non plus
du choix de c{\^o}t{\'e}.
}
 finie 
$d$
sur son centre
$k$%
.\hfill\break
 En ce cas le Lemme et son compl{\'e}ment se
pr{\'e}cisent en la~:
\vskip-1mm
\Thc Proposition 1|
Soit
$(a_i)_{1\leq i\leq d}, a_i\!\in\!K$
une base du
$k$%
-e.v.
$K$%
. Alors~:

$({\romannumeral1)}$ la famille
$(\sigma_{a_i})_{1\!\leq\!i\!\leq\!d}\!\in\!%
{\mathop{\rm End\/}\nolimits}_k(K_s)$
est une 
$K$%
-base.

$({\romannumeral2)}$ Si
$k\!\subset\!L\!\subset\!K$
est un sous-corps de
$K$
contenant le centre
$k$%
, sa dimension et celle de son commutant
${\mathop{\rm C\/}\nolimits}_K(L)$
dans
$K$
sont li{\'e}es par la relation~:
$${\mathop{\rm dim\/}\nolimits}_k(%
{\mathop{\rm C\/}\nolimits}_K(L))\!\cdot\!%
{\mathop{\rm dim\/}\nolimits}_k(L)\!=\!%
{\mathop{\rm dim\/}\nolimits}_k(K)$$

En particulier
$d=n^2$%
, la dimension de
$K$
sur
$k$
est le carr{\'e} de la dimension
$n={\mathop{\rm dim\/}\nolimits}_k(L)$
de tout sous-corps commutatif maximal
$L\!\subset\!K$%
.
\finc
\vskip-1mm
{\small
\Dem Le sous-espace
$<\!\sigma_K\!>_K\subset\!%
{\mathop{\rm End\/}\nolimits}_k(K))$
est, d'apr{\`e}s la Proposition de 
{\bf 3.5.1\/} de dimension~:
${\mathop{\rm dim\/}\nolimits}_k(<\sigma_K>_K)\!=\!%
{\mathop{\rm dim\/}\nolimits}_k(K)\!\cdot\!%
{\mathop{\rm dim\/}\nolimits}_K(<\sigma_K>_K)\!=\!%
d^2\!=\!%
{\mathop{\rm dim\/}\nolimits}_k(%
{\mathop{\rm End\/}\nolimits}_k(K))$
et donc d'apr{\'e}s le co-Corollaire de {\bf 3.2\/}
on a l'{\'e}galit{\'e}
$<\sigma_K>_K\!=\!%
{\mathop{\rm End\/}\nolimits}_k(K))$%
.\hfill\carre

Soit
$c\!=\!{\mathop{\rm dim\/}\nolimits}_L(K)$%
. On a
${\mathop{\rm dim\/}\nolimits}_k(%
{\mathop{\rm End\/}\nolimits}_L(K))\!=\!%
{\mathop{\rm dim\/}\nolimits}_k(%
{\mathop{\rm M\/}\nolimits}_{c,\, c}(L))$%
. D'apr{\`e}s la Proposition de {\bf 3.5.1\/} on a
$d\!=\!{\mathop{\rm dim\/}\nolimits}_k(K)\!=\!%
{\mathop{\rm dim\/}\nolimits}_k(L)\cdot c$
et
${\mathop{\rm dim\/}\nolimits}_k(%
{\mathop{\rm End\/}\nolimits}_L(K))\!=\!%
{\mathop{\rm dim\/}\nolimits}_k(L)\!\cdot\!%
{\mathop{\rm dim\/}\nolimits}_L(%
{\mathop{\rm M\/}\nolimits}_{c,\, c}(L))\!=\!%
{\mathop{\rm dim\/}\nolimits}_k(L)\!\cdot\!c^2$.

D'autre part d'apr{\`e}s le Compl{\'e}ment du Lemme
on a
${\mathop{\rm dim\/}\nolimits}_k(%
{\mathop{\rm End\/}\nolimits}_L(K))\!=\!%
{\mathop{\rm dim\/}\nolimits}_k(<\!\sigma_K\!>_%
{{\mathop{\rm C\/}\nolimits}_K(L)}\!)$
qui, toujours d'apr{\`e}s
la Proposition de {\bf 3.5.1\/} est
${\mathop{\rm dim\/}\nolimits}_k(%
{\mathop{\rm C\/}\nolimits}_K(L))\!\cdot\!%
{\mathop{\rm dim\/}\nolimits}_k(K)$%
. 

D'o{\`u}
$\bigl[{\mathop{\rm dim\/}\nolimits}_k(%
{\mathop{\rm C\/}\nolimits}_K(L))\!\cdot\!%
({\mathop{\rm dim\/}\nolimits}_k(L)\bigr]\cdot c\!=\!%
{\mathop{\rm dim\/}\nolimits}_k(%
{\mathop{\rm C\/}\nolimits}_K(L))\!\cdot\!%
\bigl[{\mathop{\rm dim\/}\nolimits}_k(L)%
\!\cdot\!c\bigr]\!=\!%
{\mathop{\rm dim\/}\nolimits}_k(%
{\mathop{\rm C\/}\nolimits}_K(L))\!\cdot\!d\!=\!%
{\mathop{\rm dim\/}\nolimits}_k(%
{\mathop{\rm End\/}\nolimits}_L(K))\!=\!%
\!=\!%
{\mathop{\rm dim\/}\nolimits}_k(L)\!\cdot\!c^2\!=\!%
\bigl[{\mathop{\rm dim\/}\nolimits}_k(L)%
\!\cdot\!c\bigr]\!\cdot\!c\!=\!%
{\mathop{\rm dim\/}\nolimits}_k(K)\!\cdot\!c$
et, en simplifiant cette relation entre entiers
positifs par
$c$%
, la relation de
$({\romannumeral2)}$%
.\hfill\carre

Enfin,
$k$
{\'e}tant commutatif, le corps 
$K$
contient des sous-corps commutatifs contenant
$k$%
, 
 de dimension sur
$k$
au plus
$d$%
, chacun est inclus dans un sous-corps commutatif
maximal
$L$%
.

 Pour tout
$x\!\in\!{\mathop{\rm C\/}\nolimits}_K(L))$
le sous-corps 
$L(x)\!\subset\!H$
engendr{\'e} par
$L$
et
$x$
contient
$L$
et
$x$
et est commutatif.
$L$
{\'e}tant suppos{\'e} maximal on a
$L(x)\!=\!L$
et donc
$x\!\in\!L$%
, {\it c. a d.\/}
${\mathop{\rm C\/}\nolimits}_K(L))\!=\!L$%
. 

La relation de 
$({\romannumeral2)}$
donne alors
${\mathop{\rm dim\/}\nolimits}_k(K)\!=\!%
{\mathop{\rm dim\/}\nolimits}_k(L)^2$%
.\hfill\carre\findem
}
\vskip-1mm
\Thc Corollaire 1|
Soit
$\alpha : K\rightarrow K$
un
$k$%
-automorphisme de corps.

Alors il y a
$a\in K$
tel que pour tout
$x\in K$
on a
$\alpha(x)=a^{-1}xa$%
.
\finc
\vskip-1mm
{\small
\Dem Comme
$\alpha(1)\!=\!1\!\ne\!0$
et avec les notations ci-dessus, il y a
$b_1,\cdots, b_d\!\in\!K$%
, non tous nuls, t.q. pour tout
$y\!\in\!K,\
\displaystyle\alpha(y)\!=\!\sum_{i=1}^d\sigma_{a_i}(y)b_i$%
. Donc pour tout
$x\!\in\!K, 0\!=\!\alpha(yx)-\alpha(y)\alpha(x)\!=\!$
$$\displaystyle\!=\!%
\sum_{i=1}^d\sigma_{a_i}(yx)b_i-\sum_{i=1}^d\sigma_{a_i}(y)b_i\alpha(x)\!=\!%
\sum_{i=1}^d\sigma_{a_i}(y)xb_i-\sum_{i=1}^d\sigma_{a_i}(y)b_i\alpha(x)\!=\!%
\sum_{i=1}^d\sigma_{a_i}(y)(xb_i-b_i\alpha(x))$$
et, par le Lemme 1, pour
$i\!=\!1,\ldots,\!d,\ xb_i=b_i\alpha(x)$%
. Si
$a\!=\!b_i\!\ne\!0$
on a alors~:
$\alpha(x)\!=\!a^{-1}xa$%
.\hfill\findem
}
\rmc Remarque et d{\'e}finitions|
Si
$a\!\in\!K$
la sous-
$k$%
-alg{\`e}bre
$k[a]\!\subset\!K$
engendr{\'e} par
$a$
est  un sous-corps commutatif de
$K$%
.
Le {\it degr{\'e} 
$d_a$
de
$a$
sur
$k$%
\/}
est sa dimension sur
$k$%
, 
le plus petit entier positif
$p$
tel que
$<1=a^0, a,\ldots, a^p>=<a^0,\ldots, a^{p-1}>$%
. Ainsi
$\underline{a}\!=\!(a^j)_{0\leq j\leq d_a-1}$
est base de
$k[a]$
et il y a 
$a_0,\ldots,a_{d_a-1}\!\in\!k$
tels que
$\displaystyle a^{d_a}\!=\!\sum_{j=0}^{d_a-1}a_ia^i$%
.\hfill\break
 La matrice\footnote{\small
op{\'e}rant par multiplication {\`a} droite
sur les lignes de coordon{\'e}es puisque l'on
consid{\`e}re ici le sous-corps
$k[a]$
comme
$k$%
-espace vectoriel {\`a} gauche.
} 
dans la base
$\underline{a}$
de la restriction
${\sigma_a}_{|k[a]}$
est donc
 la {\it matrice compagnon\/}~:
$$
 \begin{pmatrix}
   0&1\!\!&\cdots0\cdots&\!\!0\!\!&\cdots0\cdots&0\\
        \vdots&\vdots\!\!&\ddots&\!\!\!\!&\vdots&0\\
        0&0\!\!&\cdots0\cdots&\!\!1&\!\!\cdots0\cdots&0\\
        \vdots&0\!\!&\vdots&\!\!\!\!&\ddots&0\\
        0&0\!\!&\cdots0\cdots&\!\!0\!\!&\cdots0\cdots&1\\
        a_0&a_1\!\!&\cdots&\!\!a_k\!\!&\cdots&a_{d_a-1}
\end{pmatrix}
$$
du
{\it polyn{\^o}me minimal de
$ a$%
\/}~:
$\displaystyle P_a\!=\!\chi_{{\sigma_a}_{|k[a]}}\!=\!%
X^{d_a}-\sum_{j=0}^{d_a-1}a_jX^j\!\in\!k[X]$%
.
\finc
{\small
 \Dem
Puisque 
$k\!=\!{\mathop{\rm Z\/}\nolimits}(K)\!\subset\!%
{\mathop{\rm C\/}\nolimits}_K(a)$
la sous-alg{\`e}bre
$k[a]$
est commutative.\hfill\break
 Comme sous-espace de 
$K$%
,
elle est un
$k$%
-espace vectoriel de dimension finie.\hfill\break
 Pour tout
$0\!\ne\!c\!\in\!k[a]$
la restriction 
${\sigma_c}_{|k[a]}$
est un endomorphisme injectif du
$k$%
-espace vectoriel de dimension finie
$k[a]$%
, donc surjectif  et
$c'$
tel que
$\sigma_c(c')\!=\!1$
est inverse de
$c$
dans
$k[a]$%
.\hfill\findem
}
\Thc Corollaire 2|
Le degr{\'e} 
$d_a$
d'un {\'e}l{\'e}ment
$a\in K$
divise 
$n=d_ae_a$%
.\hfill\break
{\small
Le quotient
$e_a\!=\!{\mathop{\rm deg\/}\nolimits}_{k[a]}(L)$
{\'e}tant  degr{\'e} sur
$k[a]$
d'un
sous-corps commutatif maximal 
$L\!\ni\!a$%
.
}
\finc
\Thc Lemme 2|
Soit
$L\subset K$
un sous-corps commutati maximal.\hfill\break
Alors pour tout
$a\in K$
le polyn{\^o}me caract{\'e}ristique de
$\sigma_a\in{\mathop{\rm End\/}\nolimits}_L(K)$
v{\'e}rifie~:
$$\chi_{\sigma_a}=P_a^{e_a}$$
en particulier 
$\chi_{\sigma_a}\in k[X]$%
, il est {\`a} coefficients dans
$k$%
.
\finc
{\small
\Dem Soit
$(\alpha_j)_{1\leq j\leq e_an}$
une base du
$k[a]$%
-espace vectoriel {\`a} gauche
$K$%
. Le
$k$%
-espace vectoriel
$\sigma_K\subset<\sigma_K>_L\!=\!%
{\mathop{\rm End\/}\nolimits}_L(K)$
a donc la base
$(a^i\alpha_j), 0\leq i<d_a, 
                       1\leq j\leq e_an$%
, base dans laquelle la matrice du
$L$%
-endomorphisme {\`a} droite
$${{\mathop{\rm End\/}\nolimits}_L}_{\sigma_a} %
 :\ <\sigma_K>_L\rightarrow<\sigma_K>_L, %
[x\mapsto\sum a^i\alpha_jxb_{i,\, j}]\mapsto%
[x\mapsto a\sum a^i\alpha_jxb_{i,\, j}]$$
est une somme directe de
$e_an$
transpos{\'e}es de la matrice de
$\sigma_a : k[a]_s\rightarrow k[a]_s$%
. Son polyn{\^o}me caract{\'e}ristique est donc
$\chi_{{{\mathop{\rm End\/}\nolimits}_L}_{\sigma_a}}\!=\!%
{P_a}^{e_an}$%
. D'autre-part, un choix d'une
$L$%
-base de
$K\!=\!\oplus_{l=1}^nc_lL$
permet d'{\'e}crire
${\mathop{\rm End\/}\nolimits}_L(K_s)\!=\!K_d^n%
, f\mapsto(f(b_l))_{1\leq l\leq n}$
comme somme directe de
$n$
exemplaires du
$L$%
-espace {\`a} droite
$K$%
. Ainsi on a aussi
$\chi_{{{\mathop{\rm End\/}\nolimits}_L}_{\sigma_a}}\!=\!%
\chi_{\sigma_a}^n$%
. Les deux polyn{\^o}mes
$\chi_{\sigma_a}, P_a^{e_a}$
{\'e}tant unitaires, {\`a} coefficient dans
le corps
$L\supset k$
et ayant m{\^e}me puissance
$n^{\hbox{i\`eme}}$
sont\footnote{\small
Car leur quotient, racine
$n^{\hbox{i\`eme}}$
de l'unit{\'e} du corps
$L(X)$
des fractions rationnelles {\`a} coefficients dans
$L$%
, est dans
$L$%
, donc
$1$
car les deux polyn{\^o}mes sont unitaires.%
\findem\hfill\break
Ceci ne se g{\'e}n{\'e}ralise pas
{\`a}  un anneau commutatif non int{\`e}gre~:
dans
${\bf Z}/4{\bf Z}[X]$
on a
$X^{2}=(X+2)^{2}$%
. 
}
{\'e}gaux.\hfill\findem

}
\Thc Proposition 2|
Soit
$L\subset K$
un sous-corps commutatif maximal l'application
$$N : K\rightarrow L, %
N(a)={\mathop{\rm det\/}\nolimits}_L(\sigma_a)$$
 est polynomiale\footnote{\small
{\it c. a d.\/} si
$(a_i)_{i=1,\ldots, d\!=\!n^2}, a_i\in K$
est une 
$k$%
-base de
$K$%
, il y un polyn{\^o}me
$P\in k[X_1,\ldots, X_{d}]$
homog{\`e}ne de degr{\'e}
$n$
tel que pour tout
$\underline{x}\!=\!(x_1,\ldots, x_d)\!\in\!k^d$
on a
$N(\displaystyle\sum_{i=1}^dx_ia_i)\!=\!%
P(x_1,\ldots, x_d))$
}
homog{\`e}ne sur 
$k$
de degr{\'e}
$n\!=\!{\mathop{\rm dim\/}\nolimits}_k(L)$
sur le 
$k$%
-espace vectoriel
$K$%
.
\finc
\vskip-2mm
{\small
\Dem D{\'e}composons
$\displaystyle N(a)\!=\!%
N_1(a)+\sum_{i=2}^nN_i(a)\beta_i$
dans une
$k$%
-base de
$L$
qui compl{\'e}te
$\beta_1\!=\!1$%
. Comme
$N$
est polynomiale sur
$L$
homog{\`e}ne de degr{\'e}
$n$%
, les 
$N_i$
sont polyn{\^o}miales sur
$k$
homog{\`e}nes de degr{\'e}
$n$%
. Le Lemme 2, 
$N(a)\!=\!(-1)^n\chi_{\sigma_a}(0)\!\in\!k$
donne 
$N_1=N$%
.\hfill\findem
}
Comme
$N(a)\!=\!0$ 
ne s'annule que pour
$a\!=\!0$%
, le th{\'e}eor{\`e}me de Chevalley-Warning\footnote{\small\sl
  un polyn{\^o}me homog{\'e}ne
$P\!\in\!k[X_1,\ldots,X_N]$
 sur un corps
fini sans z{\'e}ros non trivial a son degr{\'e}
$d\geq N$
au moins {\'e}gal {\`a} son nombre de variables. Warning pr{\'e}cise~:
le nombre de z{\'e}ros d'un polyn{\^o}me 
$Q\!\in\!k[X_1,\ldots,X_N]$
de degr{\'e}
$d<N$ est nul modulo la caract{\'e}ristique 
$p$
de 
$k$%
.
}
donne, puisque
$n\geq n^2$
si et seulement si
${\mathop{\rm dim\/}\nolimits}_k(K)\!=\!n\!=\!1$%
, le~:

\Thc Th{\'e}or{\`e}me de Maclagan Wedderburn|
Tout corps fini est commutatif.
\finc
\vskip2mm
{\small
\centerline{\small\petcap
 Commentaires bibliographiques\/}
\vskip1mm

Pour les propri{\'e}t{\'e}s
(inutiles ici)
 des
 corps commutatifs finis et la preuve {\bf [Wa]\/} de Warning
du  th{\'e}or{\`e}me de Chevalley  
voir {\bf 1.7\/} de {\bf [Sa2]\/}, chap.
\uppercase\expandafter{\romannumeral1} \S1 et \S2 de {\bf [Se]}
ou \S1 du chap. 7 de {\bf [He2]\/}.

Les propri{\'e}t{\'e}s des corps gauches sont usuellement
pr{\'e}sent{\'e}es comme partie de la th{\'e}orie des modules et
anneaux semi-simple 
(Chap. \uppercase\expandafter{\romannumeral8} de {\bf [Bo]\/}) ou
Kap 12, 13 de {\bf [vW]\/},
 pour un expos{\'e}
 concis voir
le chap. \uppercase\expandafter{\romannumeral9} de {\bf [Wl]\/} et/ou
(avec plus de d{\'e}tails et exemples) {\bf [Bl]\/} et {\bf [He3]\/}.

Que le th{\'e}or{\`e}me de Wedderburn suive de celui de Chevalley remonte {\`a}
Artin, voir {\bf [Ch]\/}.

Pour le th{\'e}or{\`e}me de Wedderburn voir aussi la populaire preuve
en une page\footnote{\small
qui ouvre aussi {\bf [Wl]\/} avant un tr{\`e}s
bref expos{\'e} des propri{\'e}t{\'e}s des corps finis.
}
 de {\bf [Wt]\/},\hfill\break les d{\'e}monstrations originelles {\bf [MW]\/}, 
celle\footnote{\small
qui contient une belle arithm{\'e}tique des 
\og polyn{\^o}mes {\`a} coefficients dans un corps gauche\fgf.
}
 de {\bf [Ar1]\/}
et l'expos{\'e} \og pour sophomore\footnote{\small
{\it c. a d.\/} l'{\'e}quivalent {\`a} l'universit{\'e} de Cornell en
{\oldstyle 1961} du L2 du LMD Valentinois de {\oldstyle 2010}.
}%
\fg
de {\bf [He1]\/} qui est repris, avec une comparaison 
des autres preuves en fin du manuel {\bf [He2]\/}.
}
\npage
\nsssecp Applications~:
{\uppercase\expandafter{\romannumeral2}}
Tours de sous-corps d'un corps %
non d{\'e}nombrable.|
\Thc Lemme|
Si
$M\subset G$
est une partie infinie d'un groupe
$G$%
.

Alors le sous-groupe
$G_M\subset G$
engendr{\'e} par
$M$
dans 
$G$
est de m{\^e}me cardinal~:
$${\mathop{\rm Card\/}\nolimits}(G_M)=
{\mathop{\rm Card\/}\nolimits}(G)$$
\finc
{\small
\Dem
Quitte {\`a} changer
$M$
en
$\{e_G, m, m^{-1} | m\!\in\!M\}$
(de m{\^e}me cardinal que
$M$%
) on peut supposer 
$G_M$
image par la composition de
l'ensemble
$M^{({\Bbb N})}$
des suites de longueur finie dans
$M$
dont le cardinal est celui de 
$M$%
, d'o{\`u}
${\mathop{\rm Card\/}\nolimits} M\!\leq\!%
{\mathop{\rm Card\/}\nolimits} G_M\!\leq\!%
{\mathop{\rm Card\/}\nolimits} M^{({\Bbb N})}\!=\!%
{\mathop{\rm Card\/}\nolimits} M$%
.\hfill\findem
}
\Thc Corollaire|
Soit
$M\subset K$
une partie infinie d'un corps
$K$%
.

Alors le sous-corps
$K_M\subset K$
engendr{\'e} par
$M$
est de m{\^e}me cardinal~:
$${\mathop{\rm Card\/}\nolimits}(K_M)=%
{\mathop{\rm Card\/}\nolimits}(M)$$

\finc
\vskip-3mm
{\small
\Dem
Soit
$A_n\!\subset\!K, n\!\in\!{\Bbb N\/}$
les parties d{\'e}finies par
$A_0\!=\!M$
et les relations de r{\'e}cur\-rence
$A_{2n+1}\!=\!%
{K^\times}_{\!\!A_{2n}\setminus\{0\}}\!\subset\!K^\times$%
, sous-groupe du groupe multiplicatif de
$K$
engendr{\'e} par
$A_{2n}\!\setminus\!\{0\}$
et
$A_{2n+2}\!=\!{K^+}_{A_{2n+1}}\!\subset\!K^+$%
, sous-groupe du groupe additif de
$K$
engendr{\'e} par
$A_{2n+1}$%
. D'apr{\`e}s le  Lemme elles ont m{\^e}me
cardinal que
$M$%
, il en est de m{\^e}me de leur union
$\displaystyle K_M\!=\!%
\cup_{n\in{\Bbb N}}A_n$%
.\hfill\findem

Si
$B\!\subset\!K$
est une
$k$%
-base d'un surcorps
$K\!\supset\!k$
d'un corps
$k$
on a
$B\!\subset\!K\!=\!K_{k\amalg B}$
d'o{\`u} le~:
}
\Thc co-Corollaire|
Le cardinal d'un surcorps
$K\!\supset\!k$
d'un corps
$k$
est~:
$${\mathop{\rm Card\/}\nolimits}(K)=%
{\mathop{\rm Max\/}\nolimits}\bigl(%
{\mathop{\rm Card\/}\nolimits}(k), %
{\mathop{\rm dim\/}\nolimits}_k(K)\bigr)
$$
\finc
\Thc Proposition|
Soit
$K$
un corps de cardinal sup{\'e}rieur au d{\'e}nombrable~:
$${\mathop{\rm Card}\nolimits}(K)>%
{\mathop{\rm Card}\nolimits}({\Bbb N\/})$$
et
$I$
le plus petit ordinal de cardinal
${\mathop{\rm Card}\nolimits}(I)=%
{\mathop{\rm Card}\nolimits}(K)$
celui de
$K$%
.

Alors il y a 
$(k_i)_{i\in]0, I[}, k_i\!\subset\!K, %
{\mathop{\rm card\/}\nolimits}(k_{i})\!=\!
{\mathop{\rm Max\/}\nolimits}\bigl(%
{\mathop{\rm Card}\nolimits}(i),\,
{\mathop{\rm Card}\nolimits}({\Bbb N\/})\bigr)$%
, une famille de sous-corps infinis de
$K$
qui, si
$k_0=K_{\{1\}}$
le {\it sous-corps premier\/}\footnote{\small
{\it c. a d.\/} engendr{\'e} dans
$K$
par la partie
$\{1\}\subset K$%
.
}
de
$K$%
, est \og infiniment croisssante\fgf, {\it c. a d.\/}
 pour tout
$j\!<\!i\!\in\![0, I[$
on a~:
$$k_j\subset k_i\quad\hbox{\rm et\/}\quad
{\mathop{\rm dim\/}\nolimits}_{k_i}(k_{i+1})\!\geq\!
{\mathop{\rm Card}\nolimits}({\Bbb N\/})$$
\finc
{\small
\Dem
Soit
$${\cal K\/}^\infty_{<\!<\!K}\!=\!%
\bigl\{k\subset K\,|\, %
{\mathop{\rm Car\/}\nolimits}({\bf N\/})\!\leq\!%
{\mathop{\rm Car\/}\nolimits}(k)\!<\!%
{\mathop{\rm Car\/}\nolimits}(K)\bigr\}$$
l'ensemble des sous-corps infinis de 
$K$
de cardinal inf{\'e}rieur {\`a} celui de
$K$%
.

Pour tout
$k\!\in\!{\cal K\/}^\infty_{<\!<\!K}$
le co-Corollaire donne
${\mathop{\rm dim\/}\nolimits}_{k}(K)\!=\!%
{\mathop{\rm Card\/}\nolimits}(K)$%
. 

Ainsi, l'axiome du choix donne
 une partie libre d{\'e}nombrable
$J_k\!\subset\!K$
et donc, toujours d'apr{\'e}s le co-Corollaire,
le sous-corps de
$K$
engendr{\'e} par cette partie,
$k(J_k)$
v{\'e}rifie~:\hfill\break
$k\!\subset\!k(J_k),\  %
{\mathop{\rm Card\/}\nolimits}(k(J_k))\!=\!%
{\mathop{\rm Card\/}\nolimits}(k)$%
, donc
$k(J_k)\!\in\!{\cal K\/}^\infty_{<\!<\!K}$
et 
${\mathop{\rm dim\/}\nolimits}_k(k(J_k))\!\geq\!%
{\mathop{\rm Card\/}\nolimits}({\bf N\/})$%
.

De m{\^e}me, le corps
premier 
$k_0$
{\'e}tant fini ou d{\'e}nombrable, 
la dimension sur 
$k_0$
du corps
$K$
est  
${\mathop{\rm dim\/}\nolimits}_{k_0}(K)\!=\!%
{\mathop{\rm Card\/}\nolimits}(K)$%
, le cardinal de 
$K$%
. Il y a donc une partie 
$k_0$%
-libre d{\'e}nombrable
$J_0\!\subset\!K$%
.

La suite
$k_i\!\in\!{\cal K\/}^\infty_{<\!<\!K}, i\!\in\!]0, I[$
se d{\'e}finit par
$k_1\!=\!K_{J_0}$
le sous-corps de
$K$
engendr{\'e} par
$J_0\subset K$
et les relations de r{\'e}currence (transfinie)
$k_{i+1}\!=\!k_i(J_{k_i})$%
, corps engendr{\'e} dans
$K$
sur
$k_i$
par
$J_{k_i}$%
, et si
$j\!\in\!]0, I[$
n'a pas de pr{\'e}c{\'e}dent 
$k_i\!=\!\cup_{j<i}k_i\!=\!k_0(\cup_{j<i}J_j)$%
. Ainsi
$k_i$
est le corps engendr{\'e} dans 
$K$
par 
$I_i\!=\!\cup_{j<i}J_j$%
, partie de cardinal 
${\mathop{\rm card\/}\nolimits}(I)\!=\!
{\mathop{\rm Max\/}\nolimits}\bigl(%
{\mathop{\rm Card}\nolimits}(i),\,
{\mathop{\rm Card}\nolimits}({\Bbb N\/})\bigr)$%
, et libre sur le corps premier.
Ce sous-corps est, d'apr{{\'e}s} le co-Corollaire,
de cardinal
${\mathop{\rm card\/}\nolimits}(k_i)\!=\!
{\mathop{\rm Max\/}\nolimits}\bigl(%
{\mathop{\rm Card}\nolimits}(i),\,
{\mathop{\rm Card}\nolimits}({\Bbb N\/})\bigr)$%
.\hfill\findem
}
\rmc Remarque|
Le corps
${\bf C\/}(X)$
des fractions rationnelles {\`a} coefficients complexes
est engendr{\'e} sur
${\bf C\/}$
par\footnote{\small
il est m{\^e}me engendr{\'e} par le seul {\'e}l{\'e}ment
$X$%
.
} 
la famille d{\'e}nombrable libre
$(X^n)_{n\in{\bf N\/}}$%
.
Par d{\'e}composition en {\'e}l{\'e}ments simples,
cette famille se compl{\`e}te par
$(\frac{1}{(X-a)^{m+1}})_{(a, m)\in{\bf R\/}\times{\bf N\/}}$
en une base. Ainsi~:
$${\mathop{\rm dim\/}\nolimits}_{\bf C\/}({\bf C\/}(X))\!=\!%
{\mathop{\rm Card\/}\nolimits}({\bf N\/}\cup({\bf C\/}\times{\bf N\/}))\!=\!
{\mathop{\rm Card\/}\nolimits}({\bf C\/})\!>\!%
{\mathop{\rm Card\/}\nolimits}({\bf N\/})$$
et dans la proposition l'in{\'e}galit{\'e}
${\mathop{\rm dim\/}\nolimits}_{k_i}(k_{i+1})\!\geq\!
{\mathop{\rm Card}\nolimits}({\Bbb N\/})$
peut {\^e}tre stricte.
\finc%

\nsssecp Applications~:
{\uppercase\expandafter{\romannumeral3}}
Dimension du dual d'un espace de dimension infinie.|
{\small
Le dual d'un 
$K$%
-espace vectoriel de dimension finie
$V$
est, d'apr{\`e}s la Proposition 1 de {\bf 1.5\/}~:
$${\mathop{\rm dim\/}\nolimits}_K(V^\ast)\!=\!%
{\mathop{\rm dim\/}\nolimits}_k(V)$$
Le propos de ce sous-sous\footnote{\small
pour ne pas dire \og fin-saoul\fgf!
}%
-paragraphe est
d'{\'e}tablir le~:
}
\Thc Th{\'e}or{\'e}me d'Erd{\"o}s-Kaplanski|
Le dual 
$V^\ast={\mathop{\rm Hom\/}\nolimits}_K(V, K)$
\hfill\break
 d'un 
$K$%
-espace vectoriel de dimension infinie
$V$
est de dimension~:
$${\mathop{\rm dim\/}\nolimits}_K(V^\ast)\!=\!%
{\mathop{\rm Card\/}\nolimits}(K)%
^{{\mathop{\rm dim\/}\nolimits}_K(V)}%
$$
En particulier pour
$V$
de  dimension infinie, le dual 
$V^\ast$
 n'est pas isomorphe {\`a}
$V$%
.
\finc
\Thc Proposition|
Soit
$V_k\!\subset\!V$
une
$k$%
-structure sur un
$K$%
-espace vectoriel {\`a} gauche. Alors
${\mathop{\rm Hom\/}\nolimits}_k(V_k, k)\!\subset\!%
{\mathop{\rm Hom\/}\nolimits}_k(V_k, K)\!=\!%
{\mathop{\rm Hom\/}\nolimits}_K(V, K)$
est
$(K, k)$%
-libre.\hfill\break
Si
$V_k$
est de dimension finie, c'est une 
$k$%
-structure sur le dual
$V^\ast={\mathop{\rm Hom\/}\nolimits}_K(V, W)$%
.
\finc

{\small
\Dem Soit 
une famille
$\underline{f}=(f_i)_{i\in I}, f_i\in%
{\mathop{\rm Hom\/}\nolimits}_k(V_k, k)$
libre sur
$k$%
. Une relation lin{\'e}aire 
$\displaystyle\sum_{i\in I} f_ib_i\!=\!0, b_i\!\in\!K$
sur
$K$
donne, si pour chaque
$\displaystyle i\in I,\ b_i\!=\!%
\sum_{j\in J}b_{i,\,j}\beta_j$
est l'{\'e}criture du coefficient
$b_i\in K$
dans la base 
$(\beta_i)_{i\in I}, \beta_i\in K_s$%
, les relations lin{\'e}aires
$\displaystyle\sum_{i\in I}f_ib_{i,\,j}\!=\!%
0, j\!\in\!J, %
b_{i,\,j}\in k$
qui, puisque la famille
$\underline{f}$
est 
$k$%
-libre, donnent pour tout 
$j\in J$
et
$i\in I, b_{i,\,j}\!=\!0$%
, et donc pour tout
$\displaystyle i\!\in\!I, %
b_i\!=\!\sum_{j\in J}b_{i,\,j}\beta_j\!=\!0$
et la famille
$\underline{f}$
est
$K$%
-libre.\hfill\carre

Ainsi
${\mathop{\rm dim\/}\nolimits}_K(%
<\!{\mathop{\rm Hom\/}\nolimits}_k(V_k, k)\!>_K)\!=\!%
{\mathop{\rm dim\/}\nolimits}_k(%
{\mathop{\rm Hom\/}\nolimits}_k(V_k, k))$
et, dans le cas o{\`u}
$V_k$
est de dimension finie, le sous-
$K$%
-espace engendr{\'e} par
${\mathop{\rm Hom\/}\nolimits}_k(V_k, k)$
de dimension
${\mathop{\rm dim\/}\nolimits}_k(V_k)%
{\mathop{\rm dim\/}\nolimits}_k(k)$
$\!=\!%
{\mathop{\rm dim\/}\nolimits}_K(V)%
{\mathop{\rm dim\/}\nolimits}_K(K)\!=\!%
{\mathop{\rm dim\/}\nolimits}_K(%
{\mathop{\rm Hom\/}\nolimits}_K(V, K))$%
, d'o{\`u}
$<\!{\mathop{\rm Hom\/}\nolimits}_k(V_k, k)\!>_K=%
{\mathop{\rm Hom\/}\nolimits}_K(V, K)$%
.\hfill\findem
}

\Thc Compl{\'e}ment|
Si
$V_k$
et
$K$
sont de dimension infinie, on a l'inclusion stricte~:
$$<{\mathop{\rm Hom\/}\nolimits}_k(V_k, k)>_K\subset%
\bigl\{f\!\in\!{\mathop{\rm Hom\/}\nolimits}_K(V, K)%
\, |\,%
{\mathop{\rm dim\/}\nolimits}_k(f(V_k))<\infty\bigr\}%
\build{\subset}_{\ne}^{}
{\mathop{\rm Hom\/}\nolimits}_K(V, K)$$
\finc
{\small
\Dem L'image 
$\displaystyle f(V)\!\subset\!\sum_{i=1}^nk b_i$
de
$\displaystyle f\!=\!\sum_{i=1}^nf_ib_i\in%
<\!{\mathop{\rm Hom\/}\nolimits}_k(V_k, k)>_K%
\!\subset\!
{\mathop{\rm Hom\/}\nolimits}_K(V, K)$
est de dimension finie sur
$k$%
, alors que si
$(x_i)_{i\in I}, x_i\!\in\!K$
est une famille d'ensemble d'indice
celui d'une base
$(e_i)_{i\in I}, e_i\!\in\!V_k$%
, la forme lin{\'e}aire 
$f\!\in\!{\mathop{\rm Hom\/}\nolimits}_K(V, K)$
d{\'e}termin{\'e}e par
$f(e_i)\!=\!x_i$
est d'image
$<x_i, i\in I>_k\subset K$%
, dont la dimension peut {\^e}tre  
${\mathop{\rm Min\/}\nolimits}\bigl(%
{\mathop{\rm dim\/}\nolimits}_K(V), %
{\mathop{\rm dim\/}\nolimits}_k(K)\bigr)$%
.\hfill\findem
}

Le th{\'e}or{\`e}me d{\'e}coule 
des trois Affirmations~:
\Thc Affirmation 1|
${\mathop{\rm Card\/}\nolimits}(V^\ast)=%
{\mathop{\rm Card\/}\nolimits}(k)^%
{{\mathop{\rm dim\/}\nolimits}_k(V)}$
\finc
{\small
\Dem
Cela suit de l'exemple et la Proposition 1
 de {\bf 1.5\/}.\hfill\findem
}

\Thc Affirmation 2|
Si l'espace vectoriel
$V$
est de dimension infinie on a~:
$${\mathop{\rm dim\/}\nolimits}_K(V^\ast)\geq%
{\mathop{\rm Card\/}\nolimits}(K)$$
\finc
{\small
\Dem 
Le dual de
$V^\ast$
est de dimension infinie car, d'apr{\'e}s 
$({\romannumeral1})$
de la Proposition 3 de {\bf 1.5\/}
il contient le sous-espace de dimension infinie
$\iota_V(V)$%
, ainsi
$V^\ast$
est de dimension infinie.\hfill\break 
Ayant le r{\'e}sultat pour
$k$
au plus d{\'e}nombrable, on suppose 
${\mathop{\rm Card\/}\nolimits}(k)>%
{\mathop{\rm Card\/}\nolimits}({\bf N\/})$%
. 

On choisit une
$K$%
-base
$(e_i)_{i<I}, e_i\!\in\!V$
de
$V$%
. Pour tout
$k\!\subset\!K$%
, on note
$\displaystyle V_k\!=\!\oplus_{i<I}ke_i\!\subset V$%
, c'est une famille coissante 
(si
$k\!\subset\!k'\!\subset\!K%
, V_{k}\!\subset\!V_{k'}$%
) de
$k$%
-structures pour
$k\!\in\!{\cal K\/}_<(K)$
sur
$V$%
.

 D'o{\`u} une famille croissante
$(V_{k_i})_{0<i<I}$
de structures rationnelles pour les sous-corps
$k_i\!\subset\!K, 0\!<\!i\!<\!$
donn{\'e} par la Proposition de {\bf 3.5.3\/},

D'apr{\`e}s  le Compl{\'e}ment de la Proposition
pour 
$0\leq i<I$
on a l'inclusion stricte~:
$$<{\mathop{\rm Hom\/}\nolimits}%
_{k_i}(V_{k_i},k_i>%
_{k_{i+1}}%
\build{\subset}_{\ne}^{}
{\mathop{\rm Hom\/}\nolimits}_{k_i}(V_{k_i},k_{i+1})$$
Il y a donc,
d'apr{\`e}s l'axiome du choix, une application
$g : [0, I[\rightarrow V^\ast$
telle si
$0\!\leq\!i\!<\!I$
on a~:
$$g(i)\!\in\!%
{\mathop{\rm Hom\/}\nolimits}_%
{\!\!k_i}(V_{k_i},k_{i+1})%
\setminus\!%
<\!{\mathop{\rm Hom\/}\nolimits}_%
{\!\!k_i}(V_{k_i},k_i)\!>%
_{k_{i+1}}\subset{\mathop{\rm Hom\/}\nolimits}_K(V, K)$$
La famille
$(g(i)_{0<i<I}$
est libre et de cardinal
${\mathop{\rm Card\/}\nolimits}(I)\!=\!
{\mathop{\rm Card\/}\nolimits}(K)$%
.\hfill\findem
}
\Thc Affirmation 3|
Si
$W$
est un
$k$%
-espace vectoriel on a~:
$${\mathop{\rm Card\/}\nolimits}(W)=%
{\mathop{\rm Max\/}\nolimits}\bigl(%
{\mathop{\rm dim\/}\nolimits}_k(X), %
{\mathop{\rm Card\/}\nolimits}(k))\bigr)$$
\finc
{\small
Si
$(w_j)_{j\in J}, w_j\in W$
est une base du dual on a
${\mathop{\rm Card\/}\nolimits}(J)\!=\!%
{\mathop{\rm dim\/}\nolimits}_k(W)$
infini et~:
$$\displaystyle W\!=\!\oplus_{j\in J}w_j k\!=\!%
\cup_{n=0}^{\infty}\cup_{J'\subset J,
  {\mathop{\rm Card\/}\nolimits}(J')=n}%
\oplus_{j'\in J'}w_{j'}k$$
d'o{\'u} le r{\'e}sultat, puisque
${\mathop{\rm Card\/}\nolimits}\bigl(\{J'\subset J | %
{\mathop{\rm Card\/}\nolimits}(J')\bigr)\!=\!n\}\!=\!%
{\mathop{\rm Card\/}\nolimits}(J)$
et
${\mathop{\rm Card\/}\nolimits}\bigl(%
\oplus_{i=1}^nw_{j'_i}k\bigr)\!=\!%
{\mathop{\rm Card\/}\nolimits}(k^n)$
qui est soit un cardinal fini, soit celui de
$k$%
.\hfill\findem
}
\vskip5mm
{\small
\centerline{\small\petcap
 Commentaires bibliographiques\/}
\vskip1mm

Pour une preuve interm{\'e}diaire entre l'originale
(\S3 du chap. \uppercase\expandafter{\romannumeral9} du
vol. \uppercase\expandafter{\romannumeral2} de {\bf [Ja]\/})
et la pr{\'e}sente, faire
les exercices 2 et 3 du \S7 du chap.
\uppercase\expandafter{\romannumeral2} des
{\'e}ditions apr{\`e}s la troisi{\`e}me ({\oldstyle 1959}) de
{[Bo]\/}.

}
\npage
\centerline{\petcap Compl{\'e}ments (non exhaustifs!) de bibliographie}
\vskip 5mm
 
{\hangindent=1cm\hangafter=1\noindent{\bf   [Ar1]}\
{\petcap E. Artin\/}\pointir 
{\sl {\"U}ber einen Satz von Herrn J.H. Maclagan Wederburn\/},\hfill\break
Hamb. Abb.  {\bf 5} ({\oldstyle 1928}) pp. 245-250,
Collected 
Papers 
{\oldstyle 1965} pp. 301-306.
\par
}

{\hangindent=1cm\hangafter=1\noindent{\bf   [Ar2]}\
{\petcap E. Artin\/}\pointir 
{\sl Galois Theory\/},\hfill\break
Notre Dame Math. Lect. {\bf 2},
Chap {\uppercase\expandafter{\romannumeral 1}}
Notre Dame, {\oldstyle 1942}.\par}

{\hangindent=1,5cm\hangafter=1\noindent
{\bf [Bl]}\ {\petcap A. Blanchard\/}\pointir
{\it Les corps non commutatifs\/},
{P.U.F}, {\oldstyle 1972}.\par}

{\hangindent=1cm\hangafter=1\noindent{\bf   [Bo]}\
{\petcap N. Bourbaki\/}\pointir 
{\sl Alg{\`e}bre\/},
Hermann, {\oldstyle 1947}, {\oldstyle 1959},
CCLS {\oldstyle 1970}.\par}

{\hangindent=1cm\hangafter=1\noindent{\bf   [Ch]}\
{\petcap C. Chevalley\/}\pointir 
{\sl D{\'e}monstration d'une hypoth{\`e}se de M. Artin\/},\hfill\break
Hamb.  {\bf 11\/} {\oldstyle 1936}
pp. 73-75.\par}

{\hangindent=1cm\hangafter=1\noindent{\bf   [Co]}\
{\petcap P. Colmez\/}\pointir 
{\sl {\'E}l{\'e}ments d'analyse et d'alg{\`e}bre
(et de th{\'e}orie des nombres)\/},
Ed. Ec. Polytechnique, {\oldstyle 2009}.\par}

{\hangindent=1cm\hangafter=1\noindent{\bf   [Ga]}\
{\petcap P. Gabriel}\pointir 
{\sl Matrices, g{\'e}om{\'e}trie, alg{\`e}bre
lin{\'e}aire\/},
Cassini, {\oldstyle 2001}.\par}

{\hangindent=1cm\hangafter=1\noindent{\bf   [Go]}\
{\petcap R. Godement\/}\pointir 
{\sl Cours d'alg{\`e}bre\/},
Hermann, {\oldstyle 1966}.\par}

{\hangindent=1cm\hangafter=1\noindent{\bf   [Gu]}\
{\petcap L. Guillou\/}\pointir
{\sl S{\'e}ries de Fourier et ensemble d'unicit{\'e}s\/}\footnote{\small
voir~:
http://www-fourier.ujf-grenoble.fr/%
$_{\tilde{}}\,$%
marin/Manuscripts/bili-ens/Reignier.pdf
},
{\oldstyle 2010}.
\par}

{\hangindent=1cm\hangafter=1\noindent{\bf   [Ha]}\
{\petcap P. R. Halmos\/}\pointir 
{\sl Naive Set Theory\/},
Litton Ed. Pub., {\oldstyle 1960}.\par}

{\hangindent=1cm\hangafter=1\noindent{\bf   [He1]}\
{\petcap I.N. Herstein\/}\pointir 
{\sl Wederburn's Theorem and a Theorem of Jacobson.\/},\hfill\break
Am. Math. Month.  {\bf 68} 3 ({\oldstyle 1961}) pp. 249-251.\par}

{\hangindent=1cm\hangafter=1\noindent{\bf   [He2]}\
{\petcap I.N. Herstein\/}\pointir 
{\sl Topics in Algebra\/},
Blaisdell {\oldstyle 1964}.\par}

{\hangindent=1cm\hangafter=1\noindent{\bf   [He3]}\
{\petcap I.N. Herstein\/}\pointir 
{\sl Non commutative rings\/},
Carus Math. Mon. {\oldstyle 1968}.\par}

{\hangindent=1cm\hangafter=1\noindent{\bf   [Ja]}\
{\petcap N. Jacobson\/}\pointir 
{\sl Lectures in Abstract Algebra\/},
van Nostrand {\oldstyle 1951}, {\oldstyle 53}, {\oldstyle 64}.\par}

{\hangindent=1cm\hangafter=1\noindent{\bf   [Kr]}\
{\petcap J.-P. Krivine\/}\pointir 
{\sl Th{\'e}orie axiomatique des ensembles\/},
P.U.F., {\oldstyle 1969}.\par}

{\hangindent=1cm\hangafter=1\noindent{\bf   [Ld]}\
{\petcap E. Landau\/}\pointir 
{\sl Grundlagen der Analysis\/},
Leipzig, {\oldstyle 1930}, Trad.  Chelsea {\oldstyle 1951}.\par}

{\hangindent=1cm\hangafter=1\noindent{\bf   [Lg]}\
{\petcap S. Lang\/}\pointir 
{\sl Algebra\/},
Addison-Wesley, {\oldstyle 1965}.\par}

{\hangindent=1cm\hangafter=1\noindent{\bf   [MW]}\
{\petcap J.H. Maclagan Wederburn\/}\pointir 
{\sl A theorem on finite algebras.\/},\hfill\break
T.A.M.S.  {\bf 6} ({\oldstyle 1905}) pp. 349-352.\par}

{\hangindent=1,5cm\hangafter=1\noindent
{\bf [Sa1]}\ {\petcap P. Samuel\/}\pointir
{\it Anneaux factoriels\/},
{Publ. Soc. mat. S{\~ a}o\  Paulo}, {\oldstyle 1964}.\par}

{\hangindent=1,5cm\hangafter=1\noindent
{\bf [Sa2]}\ {\petcap P. Samuel\/}\pointir
{\it Th{\'e}orie alg{\'e}brique des nombres\/},
{Hermann}, {\oldstyle 1967}.\par}

{\hangindent=1,5cm\hangafter=1\noindent
{\bf [Se]}\ {\petcap J.-P. Serre\/}\pointir
{\it Cours d'arithm\'etique},
{P.U.F}, {\oldstyle 1970}.\par}

{\hangindent=1,5cm\hangafter=1\noindent
{\bf [vW]}\ {\petcap B. L. van der Waerden\/}\pointir
{\it Algebra},
{Springer}, {\oldstyle 1936}.\par}

{\hangindent=1cm\hangafter=1\noindent{\bf   [Wa]}\
{\petcap E. Warning\/}\pointir 
{\sl Bemerkung zur vorstehendenArbeit Herrn Chevalley\/},\hfill\break
Hamb. Abh. Math.  {\bf 11\/} {\oldstyle 1936}
pp. 76-83.\par}

{\hangindent=1,5cm\hangafter=1\noindent
{\bf [Wb]}\ {\petcap H. Weber\/}\pointir
{\it Lehrbuch der Algebra  \uppercase\expandafter{\romannumeral1}\/}.
{Braunschweig, F. Vieweg}, {\oldstyle 1894}.\par}

{\hangindent=1,5cm\hangafter=1\noindent
{\bf [Wl]}\ {\petcap A. Weil\/}\pointir
{\it Basic Number Theory\/},
{Springer}, {\oldstyle 1967}.\par}

{\hangindent=1,5cm\hangafter=1\noindent
{\bf [Wt]}\ {\petcap E. Witt\/}\pointir
{\it {\"U}ber die Kommutativit{\"a}t
endlischer Schiefk{\"o}rper\/},\hfill\break
Hamb. Abb.  {\bf 8} ({\oldstyle 1931}) pp. 413,
Ges.
Abh.
, Springer {\oldstyle 1998} pp. 395.
\par}
\vfill

\npage
\null\vskip-10mm


 %
%
\vskip3mm
\centerline{\petcap 1 Modules}
\vskip1mm
Introduction et d{\'e}finitions. \page 2

{\bf 1.1 Applications lin{\'e}aires, morphismes et isomorphismes.\/}
 \page 3

{\bf 1.2 Sous-modules et modules quotients.\/}
 \page 4

{\bf 1.3 Noyau, image, factorisation canonique et suites exactes.\/}
 \page 5

{\bf 1.4 Produit, somme directe et module libre sur un ensemble.\/}
 \page 6

{\bf 1.5 Dualit{\'e}, transposition et application canonique  vers le bidual.\/}\nobreak
 \page 9

\vskip3mm
\centerline{\petcap 2 Calcul matriciel et d{\'e}terminants.\/}
\vskip1mm
{\bf 2.1 Calcul matriciel.\/}
 \page 10

{\bf 2.2 Quelques identit{\'e}s polynomiales.\/}\hfill
 \page 12
 
\null\quad{\bf 2.2.1\/} {\sl M{\'e}tode de Gauss pour le \og syst{\`e}me %
$n\!\times\!n$
g{\'e}n{\'e}ral\fgf.\/}
 \page 12
 
\null\quad{\bf 2.2.2\/} {\sl Les identit{\'e}s remarquables de Crammer %
et du d{\'e}terminant.\/}
 \page 14
 
\null\quad{\bf 2.2.3\/} {\sl L'identit{\'e} remarquable de  %
Cayley-Hamilton.\/}
\page 17

\null\quad{\bf 2.2.4\/} {\bf Appendice~:\/} 

\quad{\sl Th{\'e}or{\`e}me fondamental de l'arithm{\'e}tique dans
${\bf Z\/}[X_1,\ldots, X_n]$%
.\/}
 \page 18

{\bf 2.3 D{\'e}terminant sur un anneau commutatif et applications.\/}
 \page 21
 
\vskip3mm
\centerline{\petcap 3 Espaces vectoriels et dimension.\/}
Introduction et d{\'e}finition. \page 23

{\bf 3.1 Combinaisons lin{\'e}aires.\/}
 \page 23

{\bf 3.2 Dimension d'un espace vectoriel.\/}
 \page 25

{\bf 3.3 Le th{\'e}or{\`e}me du rang.\/}
 \page 26

{\bf 3.4 Applications aux syst{\'e}mes lin{\'e}aires.\/}
 \page 28
 \vskip-1mm
 
\null\quad{\bf 3.4.1\/} {\sl Syst{\`e}mes lin{\'e}aires~: %
\uppercase\expandafter{\romannumeral1} le langage.\/}
 \page 27
 \vskip-1mm
 
\null\quad{\bf 3.4.1\/} {\sl Syst{\`e}mes lin{\'e}aires~: %
\uppercase\expandafter{\romannumeral2} %
rang du syst{\`e}me transpos{\'e}.\/}
 \page 28

{\bf 3.5 Structure rationnelle sur un sous-corps.\/}
 \page 30
 \vskip-1mm
 
\null\quad{\bf 3.5.1\/} {\sl D{\'e}finitions exemples et propri{\'e}t{\'e}s.\/}
 \page 30
 
\null\quad{\bf 3.5.2\/} {\sl Applications~: %
\uppercase\expandafter{\romannumeral1} Quelques propri{\'e}t{\'e}s %
des corps gauches.\/}
 \page 32
 \vskip-1mm
 
\null\quad{\bf 3.5.3\/} {\sl 
\quad\uppercase\expandafter{\romannumeral2} Tours de sous-corps d'un corps %
non d{\'e}nombrable.\/}
 \page 36
 \vskip-1mm

\null\quad{\bf 3.5.4\/} {\sl
\quad\uppercase\expandafter{\romannumeral3} Dimension du dual d'un %
espace de dimension infinie.\/}
 \page 37
\vskip2mm 
{\bf Compl{\'e}ments (non exhaustifs!) de bibliographie.\/}
\page 39

{\bf Table des mati\`eres.\/}
\page 40

{\bf Contact et affiliation.\/}
\page 40
\parc
\noindent{
Metteur en sc\`ene et secr\'etaire~: Alexis Marin Bozonat\hfill\break
\null\ courriel~: alexis.charles.marin@gmail.com}

\noindent{
Univ. Grenoble Alpes, CNRS, IF, 38100, Grenoble, France}
\finc

\end{document}